\documentclass[11pt]{amsart} 

   \usepackage[a4paper, margin=100pt]{geometry}     
     \normalsize 
    \addtocontents{toc}{\setcounter{tocdepth}{2}} 
  
    \usepackage{xcolor}
    \usepackage{hyperref}
         \hypersetup{colorlinks, citecolor=blue, filecolor=blue, linkcolor=black, urlcolor=blue}
    \usepackage{soul}
    \usepackage[numbers]{natbib}    
    \usepackage{amssymb}
    \usepackage{tikz-cd}

\usepackage{marginnote}

\newcommand{\gen}[1]{\left\langle{#1}\right\rangle}

\DeclareMathOperator{\GL}{\rm{GL}}
\DeclareMathOperator{\SL}{\rm{SL}}
\DeclareMathOperator{\Aff}{\rm{Aff}}
\DeclareMathOperator{\Hol}{\rm{Hol}}

\DeclareMathOperator{\dist}{\rm{dist}}

\DeclareMathOperator{\diam}{\rm{diam}}

\DeclareMathOperator{\diff}{\rm{Diff}}

\newcommand{\R}{\mathbb{R}}\newcommand{\N}{\mathbb{N}}

\newcommand{\T}{\mathbb{T}}\newcommand{\C}{\mathbb{C}}

\newcommand{\id}{\mathrm{Id}}

\newcommand{\dimHD}{\operatorname{dim_{HD}}}

\newcommand{\ua}{{\underline{a}}}
\newcommand{\ub}{{\underline{b}}}

\newcommand{\cA}{\mathcal{A}}

\newcommand{\cC}{\mathcal{C}}
\newcommand{\cD}{\mathcal{D}}
\newcommand{\cE}{\mathcal{E}}
\newcommand{\cF}{\mathcal{F}}

\newcommand{\cH}{\mathcal{H}}

\newcommand{\cL}{\mathcal{L}}
\newcommand{\cM}{\mathcal{M}}

\newcommand{\cP}{\mathcal{P}}
\newcommand{\cQ}{\mathcal{Q}}
\newcommand{\cR}{\mathcal{R}}

\newcommand{\cT}{\mathcal{T}}
\newcommand{\cU}{\mathcal{U}}

\newcommand{\cW}{\mathcal{W}}

\newcommand{\bC}{\mathbb{C}}

\newcommand{\bF}{\mathbb{F}}

\newcommand{\bR}{\mathbb{R}}

\numberwithin{equation}{section}

\newtheorem{theorem}{Theorem}[section] 
\newtheorem{corollary}[theorem]{Corollary}
\newtheorem{lemma}[theorem]{Lemma}
\newtheorem*{lemma*}{Lemma}
\newtheorem{proposition}[theorem]{Proposition}
\newtheorem*{proposition*}{Proposition}

\newtheorem{problem}[theorem]{Problem}

\newtheorem*{question*}{Question}
\newtheorem*{theorem*}{Theorem}
\newtheorem*{claim*}{Claim}

\newtheorem{theoremain}{Theorem}

\theoremstyle{definition}
\newtheorem{definition}[theorem]{Definition}

\theoremstyle{remark}
\newtheorem{remark}[theorem]{Remark}

    \usepackage{bigints}
    \usepackage{xfrac}
    \usepackage{mathrsfs}
    \usepackage{comment}
    \usepackage{verbatim}
    \usepackage{adjustbox}
    \usepackage{graphicx}
    \usepackage{import}

    \newcommand{\diffloc}[1]{\diff^{#1}_{\rm loc}}

    \newcommand{\Sigmafin}{\Sigma^{\rm fin}}

\addtocontents{toc}{\setcounter{tocdepth}{1}}

\begin{document}
\title[Cantor sets in higher dimensions]{Cantor sets in higher dimensions ~ I:\\
criterion for stable intersections
}

\author{Meysam Nassiri}
\address{School of Mathematics, Institute for Research in Fundamental Sciences (IPM), P.O. Box 19395-5746, Tehran, Iran}
\email{nassiri@ipm.ir}
\author{Mojtaba Zareh Bidaki}
\address{School of Mathematics, Institute for Research in Fundamental Sciences (IPM), P.O. Box 19395-5746, Tehran, Iran}
\email{mojtabazare@ipm.ir}

\thanks {This work is partially supported by INSF grant no.  4001845.}
\date{}

\begin{abstract}
    We study the geometry of dynamically defined Cantor sets in arbitrary dimensions, introducing a
    criterion for $\cC^{1+\alpha}$ stable intersections of such Cantor sets, under a mild bunching condition. This condition is naturally satisfied for perturbations of conformal Cantor sets and, in particular, always holds in dimension one. Our work extends the celebrated \emph{recurrent compact set criterion} of Moreira and Yoccoz for stable intersection of Cantor sets in the real line to higher-dimensional spaces. Based on this criterion, we develop a method  
    for constructing explicit examples of stably intersecting Cantor sets in any dimension. This construction operates in the most fragile and critical regimes, where the Hausdorff dimension of one of the Cantor sets is arbitrarily small and both Cantor sets are nearly homothetical.
    All results and examples are provided in both real and complex settings. 
\end{abstract}

\maketitle
\tableofcontents
   
\section{Introduction} 

Cantor sets appear naturally in many dynamical systems as  invariant sets. Despite all being homeomorphic, Cantor sets may exhibit drastically different geometric features in the ambient space. 

Over the past five decades, a deep connection between the fractal geometry of regular  Cantor sets and the bifurcation theory of diffeomorphisms, primarily in dimension two, has been progressively uncovered; see, for instance,   
\cite{Palis-Takens_book93,Palis-Viana-94,PY94,Buzzard,Palis-conj-2005,PY-2009,MY-2010,MMP13}.
In his pioneering works \cite{Newhouse70,Newhouse79},
Newhouse introduced the concept of thickness for Cantor sets in the real line and demonstrated the stable intersection of thick regular Cantor sets, as the milestone of the creation of Newhouse phenomenon
for surface diffeomorphisms.

Recall that a Cantor set is \emph{regular} (or dynamically defined) if it can be generated as the unique attractor of the iterated function systems of a finite family of smooth contracting maps (see \S \ref{sec: pre}). Also, in brief words, two regular Cantor sets $K,K'$ have $\cC^r$ \emph{stable intersection} if $\tilde K$ intersects $\tilde K'$ for any pair of Cantor sets $\tilde K,\tilde K'$ whose generating maps are sufficiently $\cC^r$-close to those of $K,K'$, respectively (see \cite{Moreira-IHP1996}). 
A central question in this direction is the following. 

\begin{question*}\label{problem1}
    Under what conditions does a pair of regular Cantor sets exhibit stable intersection?
\end{question*}

Moreira and Yoccoz \cite{MY} have developed a remarkable criterion for the stable intersection of Cantor sets in the real line. A family of renormalization operators on a suitable function space was associated with a given pair of regular Cantor sets ($K$, $K'$), for which it was shown that the existence of a \emph{recurrent compact set} implies $\cC^{1+\alpha}$ stable intersection of $K$ and $K'$. This allowed them to address Problem \ref{problem1} for all typical Cantor sets in the real line, proving a strong form of Palis Conjecture \cite{Palis-conj-cantors} (see also \cite{PY97,takahashi2019}).

In the complex setting, Buzzard \cite{Buzzard} has extended Newhouse's results by creating a pair of stably intersecting Cantor sets in $\C$. Buzzard's approach, 
 however, relied on different mechanisms, specifically leveraging the good isotropy properties of the conformal maps that generate these Cantor sets.
The Newhouse's criterion (gap lemma) has been recently extended to regular Cantor sets in $\C$ by Biebler \cite{BieblerGapLemma} under some extra assumptions. 

Recently, Araujo, Moreira and Zamudio \cite{AM-2023,Ara-Mor-Zam} extended the method introduced in \cite{MY} to the complex setting, proving a criterion for the stable intersection of regular (conformal) Cantor sets in $\mathbb{C}$ and extend the Moreira-Yoccoz theorem to conformal Cantor sets in $\C$ with large Hausdorff dimension. Their results establish that stability holds within the space of holomorphic (conformal) maps. A Cantor set $K$ is called conformal (resp. homothetical) if the derivative of its generating maps are orthonormal (resp. homothety) on $K$. In dimension one, real or complex, regular Cantor sets are all homothetical (and so conformal) by definition.

In higher dimensions very little is known. The only known examples of stable intersection of a pair of Cantor sets in $\R^d$  ($d>1$),  are given by Asaoka \cite{Asaoka}. Those examples are based on a purely higher dimensional phenomenon that appears in \emph{blenders} \cite{Bonatti-Diaz} (see also \cite{Moreira-daSilva,Nassiri-Pujals}). This yields pairs of Cantor sets $(K,K')$ in $\R^d$ ($d\geq 2$) exhibiting $\cC^1$ stable intersections. 
Such pairs of Cantor sets exhibit a distinctive geometry and are far from being conformal, as their generators satisfy a domination property essential for applying the blender mechanism. 
Furthermore, this phenomenon implies that both Cantor sets must have Hausdorff dimensions greater than one, thereby restricting it to higher-dimensional settings. Specifically, it requires that $\lfloor\dimHD(K)\rfloor  + \lfloor\dimHD(K')\rfloor \geq d$, where $\lfloor t\rfloor$ is the integer part of $t$. Notably, the $\cC^1$ stability of intersections does not occur for any pair of Cantor sets on the real line \cite{Moreira-acta2011}.
This indicates that 
a universal scenario may not be conceivable for all dimensions and raises numerous questions. 

The aim of this work is to address the challenge of identifying a general criterion for the stable intersection of Cantor sets in all dimensions, both in the real and complex settings.
A key step towards this goal is to obtain a fine analysis of the geometry of regular Cantor sets at arbitrarily small scales. It is achieved in $\cC^{1+\alpha}$ regularity and under a mild bunching condition (Theorem \ref{thm: Bounded geometry of regular Cantor sets}). 
This bounded geometry result may be of independent interest and shed some light on the poorly understood geometry of Cantor sets in higher dimensions.
The bunching condition here is always satisfied in dimension one and also applies to perturbations of conformal Cantor sets in higher dimensions.  
Building on this, we establish a criterion for $\cC^{1+\alpha}$ stable intersections of regular Cantor sets in arbitrary dimensions (real or complex).

Relying on our criterion, we present a method for constructing explicit examples of $\cC^{1+\alpha}$ stably intersecting Cantor sets in $\R^d$ (see \S\ref{sec: Affine QC-Cantor sets}). This construction operates in the most fragile and critical regimes, where the Hausdorff dimension of one of the Cantor sets is arbitrarily small and both Cantor sets are nearly homothetical.
In particular, we have the following.

\begin{theoremain}\label{thm: Main: Real example}
For every $d\in\N$, $\alpha\in(0,1)$ and $\epsilon>0$, there exists a pair of Cantor sets $(K, K')$ in $\R^d$ with $\cC^{1+\alpha}$ stable intersection such that $\dimHD(K) <\epsilon$. Moreover, the Cantor sets $K,K'$ can be chosen to be affine and arbitrarily close to the space of conformal (or homothetical) Cantor sets. 
\end{theoremain}
Here, as in the one-dimensional setting, $\cC^{1+\alpha}$ regularity plays a crucial role. In fact, in a forthcoming paper \cite{NZ3}, which extends the work of Moreira \cite{Moreira-acta2011} to higher dimensions, we demonstrate that for the examples constructed in Theorem \ref{thm: Main: Real example} (where $\dimHD(K) < 1$) the intersection of $K$ and $K'$ cannot exhibit $\cC^1$ stability. Thus, despite the results in \cite{Asaoka}, we propose the following problem, for which we expect a positive answer.

\begin{problem}
    Let $(K,K')$ be a pair of regular Cantor sets in $\R^d$ satisfying 
    \begin{equation*}
        \lfloor\dimHD(K)\rfloor  + \lfloor\dimHD(K')\rfloor < d.
    \end{equation*}  
    Is it always possible to remove their intersection by a $\cC^1$ small perturbation?
\end{problem}

Observe that the sum of the Hausdorff dimensions of the Cantor sets appearing in Theorem \ref{thm: Main: Real example} is less than  $d+\epsilon$. By contrast, 
it is well known that two compact sets in $\R^d$ can be separated by  small perturbations (e.g.,  translations) if the sum of their Hausdorff dimensions is less than $d$. 
This raises the question of whether this dimension constraint is optimal for regular Cantor sets in $\R^d$. In the sequel paper \cite{NZ2}, we address this question in the affirmative. 

Our method applies also in the complex setting. So, we have a similar statement for holomorphic Cantor sets in $\C^d$. 

\begin{theoremain}\label{thm: Main: Complex example}
For every $d\in\N$ and $\epsilon>0$, there exists a pair of holomorphic Cantor sets $(K, K')$ in $\C^d$ with stable intersection in the holomorphic topology such that $\dimHD(K)<\epsilon$. Moreover, the Cantor sets $K,K'$ may be affine and arbitrarily close to the  conformal (or homothetical) Cantor sets.
\end{theoremain}

 \subsection{The criterion}

One can induce a family of renormalization operators corresponding to a pair of regular Cantor sets $(K,K')$, which acts on some infinite dimensional space $\cQ$ representing relative $\cC^{1+\alpha}$ embeddings of $K$ and $K'$ in the ambient space $\R^d$.
Indeed, any pair $[h,h']\in \cQ$ represents the class of all $\cC^{1+\alpha}$ embeddings $(A\circ h)(K),(A\circ h')(K')$ in $\R^d$ of $K,K'$, where $A\in \Aff(d,\R)$ and $h,h'$ are $\cC^{1+\alpha}$ maps on $\R^d$.
Here, $\Aff(d,\R)$  is the space of invertible affine maps of $\R^d$.
Roughly speaking, the renormalization operators zoom in to the smaller parts of regular Cantor sets. The key property here is that small parts of regular Cantor sets are geometrically similar to the entire Cantor set.
Moreover, they keep some finite dimensional space $\cQ_{\Aff} \subset \cQ$ invariant.
It turns out that $\cQ_{\Aff}$ represents the relative affine infinitesimal structure of the pair of Cantor sets and $\cQ_{\Aff} \cong \Aff(d,\R)\times \Sigma^- \times \Sigma'^-$, where 
$\Sigma^-,\Sigma'^-$ are one sided shift spaces corresponding to the dynamical definitions of $K,K'$. Analogously, in the holomorphic case we study the action of renormalization operators on certain infinite dimensional space $\cQ^{\rm Hol}$ and its finite dimensional invariant subspace $\cQ_{\Aff}^{\rm Hol}\cong \Aff(d,\C)\times \Sigma^- \times \Sigma'^-$.
 For precise definitions we refer to Section
  \ref{sec: pre} and Section \ref{sec: config and renormalization}.

\begin{figure}[ht]\label{figure: affine cantor}
    \centering
    \includegraphics[width=0.6\textwidth]
    {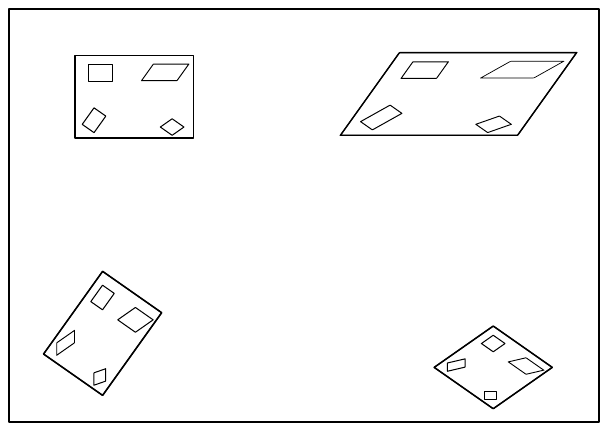}
    \caption{A Cantor set generated by four affine maps: initial approximations.}
    \label{fig:Covering regions}
\end{figure}

In the case that the generating maps of the Cantor sets $K$ and $K'$ are affine maps,
the action of renormalization operators can be simplified as follows. 
In this case renormalization operators are the pairs $(X,X')$ that act on $\Aff(d,\R)$ by  
\begin{equation*}
   H \mapsto X^{-1}\circ H \circ X',
\end{equation*}
where $(X,X')$ is either $(\id, R')$ 
or $(R,\id)$, and $R,R'$ are affine generators of $K,K'$, respectively.

\begin{theoremain}[Covering criterion for stable intersection] \label{thm:A}
Let $d\in\N$, and $(K,K')$ be a pair of bunched Cantor sets in $\R^d$ with the corresponding family $\cR$ of renormalization operators. Assume that there exists a bounded open set $\cU \subset \cQ_{\Aff}$ satisfying the covering condition
\begin{equation}\label{eq:covering-intro}
    \overline{\cU} \subset \bigcup_{\Psi\in \cR} \Psi^{-1}( \cU).
\end{equation}
Then the Cantor sets $h(K)$ and $h'(K')$ have $\cC^{1+\alpha}$ stable intersection for all $[h,h']\in \cU$. In particular, if $[\id,\id]\in \cU$ then $K$ and $K'$ have $\cC^{1+\alpha}$ stable intersection.
\end{theoremain}

Note that the condition \eqref{eq:covering-intro} is finite dimensional, while the conclusion holds $\cC^{1+\alpha}$-stably in the infinite dimensional space of regular Cantor sets.

Theorem \ref{thm:A} extends  to the complex setting and for holomorphic Cantor sets.

\begin{theoremain}
\label{thm:A-complex}
Let $d\in\N$, and $(K,K')$ be a pair of bunched holomorphic Cantor sets in $\C^d$ with the corresponding family $\cR$ of renormalization operators. Assume that there exists a bounded open set $\cU \subset \cQ_{\Aff}^{\rm Hol}$ satisfying the covering condition \eqref{eq:covering-intro}. 
Then $h(K)$ and $h'(K')$ have stable intersection for all $[h,h']\in \cU$. In particular, if 
$[\id,\id]\in \cU$ then $K$ and $K'$ have stable intersection. The stability holds in the space of holomorphic maps.
\end{theoremain}

\begin{remark}
    For Cantor sets on the real line, the covering condition \eqref{eq:covering-intro} is equivalent to the recurrent compact set criterion introduced by Moreira and Yoccoz \cite{MY}. In particular, Theorems~\ref{thm:A} and~\ref{thm:A-complex} generalize this criterion and its complex counterpart developed in \cite{AM-2023}, respectively.
\end{remark}

\subsection{About the proofs}
The general strategy of the proof of Theorem \ref{thm:A} is similar to the work of Moreira-Yoccoz \cite{MY}. However, to generalize their renormalization method to higher dimensions, one needs to deal with the non-conformal behavior of typical maps. 
First, we need to analyze the infinitesimal geometry of regular Cantor sets. Thus, a main part of the proof is devoted to the convergence of limit geometries for the sequences of contracting  maps in arbitrary dimensions. Such convergence results crucially require  
 a bunching condition and also the
 $\cC^{1+\alpha}$ 
 regularity.
Fine geometric control under iterated contracting maps, as studied in \cite{FNR-AiM}, addresses the behavior of balls under iterations of sequences of contracting $\cC^{1+\alpha}$ maps with a quasi-conformality condition. 
Here, we broaden this geometric control to encompass any sequence of contracting maps satisfying a bunching condition (see Definition \ref{def bunched}). This condition yields the bounded geometry (Theorem \ref{thm: Bounded geometry of regular Cantor sets}) and holds in various settings, including sequences of contracting maps in the vicinity of conformal contracting maps. 
 This type of assumptions appears in many other contexts in dynamical systems (e.g. in stable ergodicity, regularity of holonomy maps, etc.) as well as in other areas of mathematics. On the other hand, as we discussed in Remark \ref{remark:bounded geomtery is true for affine Cantors}, the affine case automatically verifies the bounded geometry. 
Thus, it remains a big challenge to understand infinitesimal geometry and the intersection of Cantor sets violating the bunching condition.

To prove Theorems \ref{thm: Main: Real example} and \ref{thm: Main: Complex example}, presenting examples of stably intersecting Cantor sets in arbitrary dimensions, we apply the covering criterion in Theorems \ref{thm:A}, \ref{thm:A-complex}, respectively.
In order to verify the covering criterion in higher dimensions, we require a covering result in the space of linear maps, motivated by the approach of \cite{FNR-AiM}.
A part of the construction is inspired by the examples of affine stably intersecting Cantor sets on the real line developed in \cite{Pourbarat2015,Pourbarat2019}.

\subsection*{Organization} 
In Section \ref{sec: pre}, basic notations and definitions are established, including the bunching condition for regular Cantor sets. In Section \ref{sec: covering condition}, we discuss the abstract covering condition on any metric space.  In Section \ref{sec: convergence geometries}, we prove the control of shape property for a sequence of bunched contracting local diffeomorphisms in $\R^d$. Then we focus on the main object of study: Cantor sets at infinitesimal scales, also known as limit geometries. Section \ref{sec: config and renormalization} defines renormalization operators and studies them on affine Cantor sets as special cases. The covering criterion, which is the main result of this paper, is proved in Section \ref{sec: covering criterion}. In Section \ref{sec: Affine QC-Cantor sets}, we develop a method to verify the covering criterion, specifically proving Theorems
\ref{thm: Main: Real example} and \ref{thm: Main: Complex example}.

\subsection*{Acknowledgement}
The authors express their gratitude to Carlos Gustavo Moreira, Enrique Pujals, Federico Rodriguez Hertz, and Reza Seyyedali for their fruitful comments and conversations. Special thanks are extended to Mehdi Pourbarat for his valuable critiques and discussions during the preparation of this work. The authors also thank ICTP for its support and hospitality, where a part of this work was written.

\section{Preliminaries}\label{sec: pre}

\subsection{Basic notations} Here, we fix some basic notations to be used frequently in the subsequent sections.

\begin{itemize}
\item We may use the notation $\bF$ as either $\R$ or $\C$ when both cases can be described in a same way.
\item We denote the Hausdorff dimension of the set $A\subset \R^d$ as $\dimHD(A)$.
\item We denote $\mathrm{O}(d)$ as the group of orthonormal linear transformations over $\R^d$.
\item 
Given a metric space $(X,d)$, for any $x\in X$ we denote the $\delta$-neighborhood of $x$ in $X$ as $B_{\delta}(x)$. For $V\subset X$ we define 
$B_{\delta}(V):=\bigcup_{v\in V}B_{\delta}(v)$ and
$V_{(\delta)}:=\{v\in V:\; B_{\delta}(v)\subset V \}$.
\item 
Given compact sets $A,B\subset \R^d$, denote $\dist(A,B):=\min_{x\in A,\;y\in B}{|x-y|}$.

\item 
Denote the space of invertible matrices over $\bF^d$ by $\GL(d,\bF)$ and those with determinant equal to $1_{\bF}$ as $\SL(d,\bF)$. For $A\in \GL(d,\bF)$ we denote its norm by $\|A\|_{\it op}:= \sup_{|v|=1} |Av|$ and its co-norm by $m(A):=\inf_{|v|=1} |Av|$.

\item 
We denote  by $\Aff(d,\bF)$ the space of invertible affine transformations $[x\mapsto Ax+a]$ of $\bF^d$ equipped with $\cC^{1}$ topology where $A\in \GL(d,\bF)$ and $a\in \bF^d$.
\item We denote the identity matrix in $\GL(d,\bF)$ as $\id$. We also use the notation $\id$ for the identity function. Each case will be clear in the context.
\item $\GL(d,\C)$ can be interpreted as a subgroup of $\GL(2d,\R)$ since any  
$T= A+iB\in \GL(d,\C)$ with $A,B\in M(d,\R)$ maps a vector $x+iy\in \C^d$ where $x,y\in \R^d$ to the vector $(Ax-By)+i(Bx+Ay)$. Therefore, $T$ acts as a linear map on $\R^{2d}$ which sends the vector $(x,y)\in \R^{2d}$ to $(Ax-By,Bx+Ay)\in \R^{2d}$ and $A,B$ are such that this mapping is not singular. This implies that $\Aff(d,\C)$ is a subgroup of $\Aff(2d,\R)$.
 
\item (H\"older regularity) The $\cC^\alpha$ semi-norm, $\cC^1$ norm and $\cC^{1+\alpha}$ norm of an $\cC^{1+\alpha}$ map $f:X \to \R^d$ on the domain $U\subset X$ are defined, respectively by
\begin{align*}
    \|f\|_{\cC^{0}}&:= \sup_{u\in U} |f(u)|,\\
    |f|_{\alpha}&:=\sup_{x,y\in U} \frac{|f(x)-f(y)|}{d(x,y)^\alpha},\\
    \|f\|_{\cC^1} &:= \|f\|_{\cC^0}+ \sup_{u\in U} \|Df(u)\|_{\it op},\\
     \|f\|_{\cC^{1+\alpha}}&:= \|f\|_{\cC^1}+|Df|_{\alpha}.
\end{align*} 
\item We denote by $\diffloc{r}(M)$, the space of 
all $\cC^r$ diffeomorphisms $f : \mathrm{Dom}(f) \to \mathrm{Im}(f)$ such that $\mathrm{Dom}(f)$,
$\mathrm{Im}(f)$
are open subsets of $M$.  $f,g\in \diffloc{r}(M)$ are $\cC^{r}$-close if there exists a diffeomorphism $h:\mathrm{Dom}(f)\to  \mathrm{Dom}(g)$ $\cC^{r}-$close to $\id|_{\mathrm{Dom}(f)}$ such that the map $g\circ h$ is $\cC^{r}$-close to $f$. In other words, two elements of $\diff_{\rm loc}^{r}(M)$ are $\cC^{r}$-close if their graphs  are $\cC^{r}$-close embedded submanifolds of $M\times M$. 
\end{itemize}

\subsection{Regular Cantor sets}

In this subsection we define regular Cantor sets in $\R^d$ (or $\R^{2d}\cong \C^d$). We keep the notations similar to the ones used in  \cite{AM-2023}.

\begin{definition} \label{def:Cantor}
 A regular (or dynamically defined) Cantor set in $\R^d$ is defined by the following data:

\begin{itemize}

\item A finite set $\mathfrak{A}$ of letters and a set $\mathfrak B \subset \mathfrak{A} \times \mathfrak{A} $ of admissible pairs.

\item For each $a\in \mathfrak{A}$ a compact connected set $G(a)\subset \R^d$. In the complex case $G(a)$ is a compact subset of $\R^{2d}\cong \C^d$.
 
\item A $\cC^{1+\alpha}$ map $g: V \to \R^d$ defined in an open neighborhood $V$ of $\bigsqcup_{a\in \mathfrak{A}}G(a)$. In the complex case $g$ maps $V$ into $\R^{2d}\cong \C^d$.
\end{itemize}

This data must verify the following assumptions:

\begin{itemize}
    
\item The sets $G(a)$, $a \in \mathfrak{A}$ are pairwise disjoint

\item $(a,b)\in \mathfrak B $ implies $G(b) \subset g(G(a))$, otherwise $G(b) \cap g(G(a)) = \emptyset$.

\item For each $a \in \mathfrak{A}$ the restriction $g|_{G(a)}$ can be extended to a $\cC^{1+\alpha}$ embedding (with  $\cC^{1+\alpha}$ inverse) from an open neighborhood of $G(a)$ onto its image such that $m(Dg):= \inf\{m(Dg(x)): x\in V\}>\mu^{-1}$ for some positive constant $\mu<1$. In the complex case $Dg(x)$ is a member of $\GL(2d,\R)$ for all $x\in V$.

\item The subshift $(\Sigma_{\mathfrak{B}}, \sigma)$ is called the symbolic type of the Cantor set 
\begin{equation}\label{eq: definition of symbolic type Sigma}
    \Sigma_{\mathfrak{B}}:=\{ \ua= (a_0, a_1, a_2, \dots  ) \in \mathfrak{A}^{\N}:(a_i,a_{i+1}) \in \mathfrak B,\; \text{for all} \;i \geq 0\}
\end{equation} 
with the topologically mixing shift map
$\sigma (a_0,a_1,a_2, \dots) := (a_1,a_2,a_3, \dots)$.
\end{itemize}

Once we have all these data we can define a Cantor set (i.e. a nonempty, perfect, compact and totally disconnected set) on $\R^d$ (or $\R^{2d}\cong \C^d$),
\[ K=\bigcap_{n \geq 0} g^{-n}\left(  \bigsqcup_{a \in \mathfrak{A}} G(a) \right). \]   
\end{definition}

We say that $K$ is a Cantor set described by the data $(\Sigma_{\mathfrak{B}},g)$. Whenever the set $\mathfrak{B}$ is fixed and clear from the context we simply say that $K$ is constructed via $(\Sigma,g)$.  

\begin{remark}
    Note that here we defined regular Cantor sets with an expanding generator. We also could define them with contracting generators which are inverse of the expanding generator $g$ restricted to domains $G(a)$. Sometimes we will use this definition of regular Cantor sets. More precisely, given contracting maps $f_i : \overline U\to G(a_i)$ for $i=1,\dots,m$ where $m:=\#\mathfrak{A}$, $a_i\in \mathfrak{A}$ and $\overline U\subset \R^d$ is a compact set, the maximal invariant set of the Iterated Function System (IFS) generated by $f_1,\dots,f_m$ with the symbolic type $\Sigma$ is the regular Cantor set $K$. In that case we say that the regular Cantor set $K$ is generated by the Iterated Function Systems $\{f_1,\dots,f_k\}$ and the subshift $(\Sigma,\sigma)$.
    In the case that $(\Sigma,\sigma)$ is the full shift dynamics we say that the Cantor set $K$ is constructed via the IFS $\{f_1,\dots,f_{m}\}$.
\end{remark}
A regular Cantor set $K$ may be constructed in many ways as described above. For instance, the standard middle third Cantor set may be constructed via both $(\Sigma_{\mathfrak{B}_i},g_i)$, 
where $\mathfrak{A}_1 :=\{1,2\}$, $\mathfrak{A}_2:=\{1,2,3,4\}$, $\mathfrak{B}_i= \mathfrak{A}_i\times \mathfrak{A}_i$ for $i=1,2$ and
$$g_1 : \bigcup_{j=0,2} \left[\frac{j}{3},\frac{j+1}{3}\right] \to [0,1],\qquad
g_2 :\bigcup_{j=0,1,3,4} \left[\frac{2j}{9},\frac{2j+1}{9}\right] \to [0,1],$$
such that 
$g_1(x) = 3x-j$ for $x\in[\frac{j}{3},\frac{j+1}{3}]$ and
$g_2(x)=9x-2j$ for  $x\in[\frac{2j}{9},\frac{2j+1}{9}]$.
Therefore, when we say that $K$ is a regular Cantor set we assume a set of data 
$(\Sigma_{\mathfrak{B}},g)$ is fixed. For convenience, we will just mention this set of data most of the times by the expanding map $g$ since all the data can be inferred if we know $g$ or by an IFS of contracting generators as we explained above.

\subsection{Bunched Cantor Sets}\label{sec: space of bunched Cantor sets}

\begin{definition} \label{def bunched}
Given a regular Cantor set $K$ in $\R^d$
 described by $(\Sigma_{\mathfrak{B}},g)$ where $g$ is a $\cC^{1+\alpha}$ map with a uniform expansion rate bigger than $\mu^{-1}>1$, we say that $K$ is
\begin{itemize}
\item \emph{affine} if for any $a\in \mathfrak{A}$, $g|_{G(a)}: G(a) \to \R^d$ is an expanding affine map.
    \item \emph{homothetical} if $Dg(x)/\|Dg(x)\|_{\it op}=\id$ for all $x\in K$,
    \item \emph{conformal} if $Dg(x)/\|Dg(x)\|_{\it op}\in \operatorname{O}(d)$ for all $x\in K$,
    \item  \emph{bunched} whenever $g$ is bunched at the Cantor set $K$, that is, there is $N_g\in \N$ such that for all $x \in K$,
    $Dg(x) \in \GL(d,\R)$ and
    
\begin{equation}\label{eq: definition of bunched cantor sets}
    \|Dg^{N_g}(x)\|_{\textit{op}}\cdot 
   m\left(Dg^{N_g}(x)\right)^{-1} \cdot \mu^{\alpha \cdot N_g } <1. 
\end{equation}
\item \emph{holomorphic} if $d$ is an even integer and for any $x\in K$, $Dg(x)$ lies in $\GL(d/2,\C)$ which is a subgroup of $\GL(d,\R)$.

\end{itemize}
\end{definition}

\begin{remark}
    For convenience and in order to make the presentation as transparent as possible we assume that $N_g=1$ throughout the paper.
    This is not an essential restriction because by using an adapted metric one can assume that $N_g=1$ in \eqref{eq: definition of bunched cantor sets}.
\end{remark}

\begin{remark} 
Conformal Cantor sets are special cases of bunched (or bunched holomorphic) Cantor sets. Conformal Cantor sets in $\R^2$ are defined in \cite{AM-2023} to analyze regular Cantor sets that appear in $\C$. Note that any holomorphic Cantor set in $\R^2$ is bunched and coincides with conformal Cantor sets. Moreover, the bunching condition always holds in a neighborhood of conformal Cantor sets. Therefore, small perturbations of conformal Cantor sets are bunched (holomorphic) Cantor sets.
\end{remark}

Given a regular Cantor set $K$ generated by $g$, $K$ is the maximal invariant compact subset of $\R^d$. Indeed, the pair $(K,g|_K)$ is a dynamical system which is conjugate to the subshift $\Sigma_{\mathfrak B}$. The conjugacy map $h: K \to  \Sigma$ is a homeomorphism that maps each $x\in K$ to the sequence of letters $\{a_n\}_{n \geq 0}$ gained by the orbit of $x$. In particular, $g^n(x) \in G(a_n)$ for all $n$.
Associated to the Cantor set $K$ we define the sets
\begin{align*}
    \Sigma^{\mathrm{fin}} &:=\{(a_0, \dots ,a_n): n\in \N,\;(a_i,a_{i+1}) \in \mathfrak B \; \text{for all} \;i ,\; 0 \leq i < n \},\\
\Sigma^-&:=\{(\dots, a_{-n}, a_{-n+1},\dots,a_{-1},a_0): (a_{i-1},a_i) \in \mathfrak B \; \text{for all} \; i \leq 0\}.
\end{align*}
Given $\ua=(a_0, \dots, a_n)\in \Sigma^{\rm fin}$ we say that the length of $\ua$ is $n$. For $\underline{\theta}=(\dots,\theta_{-2},\theta_{-1},\theta_{0})\in \Sigma^-$ where $\theta_0=a_0$ and $\underline{b}=(b_0,b_1,\dots)$ we denote
\begin{align*}
   \underline{\theta a} &:= (\dots, \theta_{-2},\theta_{-1}, a_0, \dots, a_n ),\\
 \underline{b}_k &:= (b_0,b_1,\dots,b_k),\\\underline{\theta}_k &:= (\theta_{-k},\dots,\theta_{-1},\theta_{0})  
\end{align*}
for all $k$. Furthermore, we define the set 
$$G(\ua):= \{x \in \bigsqcup_{a \in \mathfrak{A}} G(a) ,  g^j(x) \in G(a_j), j=0,1,\dots, n \}$$
and the map $f_{\ua}: G(a_n) \to G(\ua)$ by
$$ f_{\ua} :=  \left(g|^{-1}_{G(a_0)} \circ g|^{-1}_{G(a_1)} \circ \dots \circ g|^{-1}_{G(a_{n-1})}\right)\Big|_{G(a_n)} . $$
Notice that $f_{(a_i, a_{i+1})} = g|^{-1}_{G(a_i)}$, so we have
$$f_{\ua} = f_{(a_0,a_1)}\circ f_{(a_1,a_2)}\circ \cdots \circ f_{(a_{n-1},a_n)}.$$
This implies that the finite collection of maps $\{f_{(a,b)}\;|\; (a,b)\in \mathfrak B\}$ generates the family $\{f_{\ua}\}_{\ua\in \Sigmafin}$
by their composition through the set $\Sigma$. Therefore, $K$ will be the attractor of the IFS generated by this family through to the symbolic type $\Sigma$.

Notice that in the definition of regular Cantor set $K$, pieces $G(a)$ may have an empty interior. However, we can replace pieces $G(a)$ with open relatively compact connected pieces $G^*(a):= B_\delta(G(a))$ for sufficiently small $\delta>0$ and extend the map $g|_{(G(a)}$ to the neighborhood of $G^*(a)$ for all $a\in \mathfrak{A}$ such that the pieces $\overline{G^*(a)}$ and the map $g$ satisfy the properties enumerated in Definition \ref{def:Cantor}. Let also 
$$G^*(\ua):= \{x \in \bigsqcup_{a \in \mathfrak{A}} G^*(a) ,  g^j(x) \in G^*(a_j), j=0,1,\dots, n \}.$$
With this notation, we have the following lemma from \cite[Lemma 2.1]{AM-2023}.
\begin{lemma} \label{lemma introduction diameter of images}
    Let $K$ be a regular Cantor set and $G^*(a)$ the sets defined above. There exists a constant $C > 0$ such that
    $$\diam(G^*(\ua)) < C\mu^{n},$$
    where $\mu < 1$ is such that $m(Dg) > \mu^{-1}$ in $\bigsqcup_{a\in \mathfrak{A}} G^*(a)$ and 
     $\ua = (a_0, \dots, a_n)$.
\end{lemma}
 Lemma \ref{lemma introduction diameter of images} implies that the substitution of the pieces 
 $G(a)$ by $G^*(a)$ will not change the Cantor set $K$ since  $G(\ua)\subset G^*(\ua)$ and diam$(G^*(\ua))\rightarrow 0$. So, we may consider all of the pieces $G(a)$ to be compact sets, which are the closures of their interiors.

\begin{definition}\label{def: space of cantor sets} Let $\beta>1$ an $\Sigma:=\Sigma_{\mathfrak{B}}$ be a symbolic type defined by a given $\mathfrak{A},\mathfrak{B}$ as \eqref{eq: definition of symbolic type Sigma}. We denote $\Omega_{\Sigma,d}^{\beta}$ as the set of all bunched Cantor sets $K$ in $\R^d$ described by the data $(\Sigma,g)$ where $g$ is a 
$\cC^{\beta}$ map. We equip $\Omega_{\Sigma,d}^{\beta}$ with the topology generated by $\delta$-neighborhoods $U_{K,\delta}$ of $K\in \Omega_{\Sigma,d}^{\beta}$. Each $U_{K,\delta}$ consists of bunched Cantor sets $K'$ with symbolic type $\Sigma$ generated by some $\cC^{\beta}$ map $g'$ with pieces $G'(a)$ with the property that $g'|_{G'^*(a)}$ is $\delta$-close to $g|_{G^*(a)}$  in $\diffloc{\beta}(\R^d)$, for all $a\in \mathfrak{A}$.
In the holomorphic case, we denote the set of all bunched holomorphic Cantor sets in $\R^{2d}$ with a topology on it in the same manner and denote it by $\Omega^{\Hol,\beta}_{\Sigma,2d}$.
\end{definition} 

\begin{remark}
    Notice that the bunching condition \eqref{eq: definition of bunched cantor sets} for a Cantor set is stable in the $\cC^{1+\alpha}$ topology for $\alpha>0$.
\end{remark}

\section{The covering conditions}\label{sec: covering  condition}
Given a family of maps $\cF$, we denote $\gen{\cF}^+$ as a semigroup generated by $\cF$.

\begin{definition}[Covering condition]\label{def: Covering}
    Let $\cF$ be a family of continuous maps on a metric space $(X,d)$. We say that a set $V\subset X$ satisfies the covering condition with respect to $\cF$ if
    \begin{equation}\label{eq covering condition}
        \overline{V} \subset \bigcup_{f\in \cF} f^{-1}(V).
    \end{equation}
    \end{definition}
This condition implies that for any $v\in \overline{V}$ there is $f \in \cF$ such that $f(v)\in V$. In other words, one can map any $v\in \overline{V}$ into $V$ by one iteration of the elements of $\cF$.

    \begin{definition}[Strong covering condition]\label{def: Strong covering}
    Let $\cF$ be a finite family of continuous maps on a metric space $(X,d)$. We say that a set $V\subset X$ satisfies the strong covering condition with respect to $\cF$ if there exists $\delta>0$ such that 
    \begin{equation}\label{eq: covering condition delta interior}
         B_{\delta}(V) \subset \bigcup_{f\in \cF} f^{-1}(V_{(\delta)}).
    \end{equation}
    \end{definition}
The following two lemmas are direct consequences of continuity argument and the proofs are omitted. 
    \begin{lemma} \label{lemma covering to strong}
        For a finite family $\cF$ of continuous maps on a locally compact metric space $(X,d)$, the covering condition implies strong covering condition for an open relatively compact set $V\subset X$.
    \end{lemma}
     
Given a finite family $\mathcal{F}=\{f_1,\ldots,  f_k\}$ of continuous maps, we say that the family $\tilde{\cF}$ is $\epsilon$-close to $\cF$ in the $\cC^0$ topology, if $\tilde{\cF}=\{\tilde{f}_1,\ldots,\tilde{f}_k\}$ is such that  $f_i,\tilde{f}_i$ are $\epsilon$-close in the $\cC^0$ topology for $i=1,\ldots,k$.
    
    \begin{lemma}[Stability of strong covering]\label{lemma stability of strong covering} Let $\cF$ be a finite family of continuous maps on a metric space $(X,d)$ and $V\subset X$ be a subset that satisfies strong covering condition with respect to $\cF$. Then for any family $\tilde{\cF}$ sufficiently close to $\cF$ in $\cC^0$ topology, $V$ satisfies strong covering condition with respect to $\Tilde{\cF}$.
    \end{lemma}

In the following lemma we make a bridge from  (strong) covering condition with respect to a family $\cF$ to the (strong) covering condition with respect to its generated semigroup  $\gen{\cF}^+$.

\begin{lemma}\label{lemma: immediately covering proof}
    Let $\cF$ be a finite family of uniformly continuous maps on the metric space $(X,d)$ such that an open set $V\subset X$ satisfies strong covering with respect to a finite family $\cF' \subset \gen{\cF}^+$. Then, there exists $W\subset X$ containing $V$ and satisfying strong covering condition with respect to $\cF$.
\end{lemma}

\begin{proof} Let $\cF:=\{f_1,\cdots,f_k\}$.
By the assumptions there are $\epsilon>0$, $N\in\N$,  such that 
for any $v\in B_{\epsilon}(V)$ there exists a word $\underline{w}=(w_n,\dots, w_2, w_1)\in \{1, \dots, k\}^n$ with $n<N$ such that 
    $B_\epsilon(f_{\underline{w}}(v)) \subset V$  where, $f_{\underline{w}}:= f_{w_n}\circ \cdots \circ f_{w_2} \circ f_{w_1}$. Since elements of the finite family $\cF$ are uniformly continuous then there exist positive real numbers $\delta_0
,\delta_1,\dots,\delta_N<\epsilon$ such that for every $i\in \{1,\dots,N\}$, $f\in \cF$ and $x,y\in X$, if $d(x,y)<\delta_{i-1}$ then $d(f(x),f(y))<\delta_i/4$. We also denote 
\begin{equation}\label{eq: delta def, covering}
    \delta: = \dfrac{1}{2}\min\{\delta_0,\delta_1,\dots, \delta_n\}.
\end{equation}
Fixing $v\in B_{\epsilon}(V)$, then there is $\underline w = (w_n,\dots,w_1)$ such that $B_{\epsilon}(f_{\underline w}(v))\subset V$. 
Now let $v_0 :=v$, $W_0: =B_{\delta_0}(v) $ and $v_i: = f_{w_i}(v_{i-1})$, $W_{i}:= B_{\delta_{i}}(v_i)$ for $i\in {1,2,\dots,N}$. It follows from the properties of $\delta_i$'s that $B_{\epsilon/2}(W_n)\subset V$. In particular, by \eqref{eq: delta def, covering} we have $B_{\delta}(f_{w_i}(B_{\delta}(W_{i-1})))\subset W_{i}$ for $i\in \{1,2,\dots,n\}$. So if we define $U_v: = \bigcup_{i=0}^N W_i$ then $W := \bigcup_{v\in V} U_v$ has strong covering condition. More precisely, 
$$B_{\delta}(W) \subset \bigcup_{f\in \cF} f^{-1}(W_{(\delta)}).$$
\end{proof}
\begin{remark}\label{remark: from covering to immediate covering open compacts}
This proof shows that if $V$ is open relatively compact then one can assume continuity of elements of the finite family $\cF$ instead of uniform continuity.
\end{remark}

\subsection{Covering in linear groups} In the rest of this section we will study group operation of Lie groups $\SL(d,\bF)$, for $\bF$ being either $\bR$ or $\bC$. Here, we recall covering condition for left action of elements of such a group. Assuming $G$ is a topological group, let $\cF\subset G$ be a set of elements of this group. We say that $\cU\subset G$ satisfies covering condition with respect to $\cF$ if
\begin{equation}\label{eq covering condition for topological groups}
        \overline{\cU} \subset \bigcup_{f\in \cF} f^{-1}\cU,
    \end{equation} 
where $f^{-1}\cU = \{f^{-1}u:u\in \cU \}$. By Lemma \ref{lemma covering to strong}, $\cU$ being open relatively compact implies that the covering condition is equivalent to the strong covering. 
We use both real and complex versions of \cite[Lemma 3.8]{FNR-AiM}. Proof of its complex version is analogous to the real case with the same computations in scalar field $\C$ instead of $\mathbb{R}$.

\begin{lemma}\label{lem: generating a convering in SL with d^2 maps} Let $\cU_0$, $\cU_1$ be neighborhoods of $\id$ in $\SL(d,\bF)$. Then, there exists an open relatively compact set $\cU\subseteq \cU_0$ and a finite set $\cF \subset \cU_1$ such that $\cU$ satisfies strong covering condition with respect to $\cF$. Moreover, $\cF$ has at most 
$1+\dim(\SL(d,\bF))$ elements.
\end{lemma}

\section{Convergent geometries}\label{sec: convergence geometries}  

An important tool for study the infinitesimal parts of a regular Cantor set is the convergence of geometries. Roughly speaking, it says that given a sequence of smooth contracting maps satisfying a bunching condition (which always holds in dimension 1, either real or complex), the tail of their composition behaves like the composition of contracting affine maps. 

\subsection{Existence and continuity of limits}

For a sequence $\{h_n\}_{n\in \N}$ we denote $h^n: = h_n\circ \cdots \circ h_1$.
Given $\rho,C>0$ and $\mu,\mu',\alpha \in (0,1)$ and $\kappa \geq 1$, we consider the following hypotheses for a sequence of maps $\{h_n\}_{n\in \N}$ in $\diffloc{1+\alpha}(\mathbb{R}^d)$.
\begin{itemize}
    \item[(H0)] $h_n:B_{\rho}(0)\rightarrow h_n(B_{\rho}(0))$ is a $\cC^{1+\alpha}$ diffeomorphism fixing the origin,
    \item[(H1)] $\|h_n\|_{\cC^{1+\alpha}}<C$,
    \item[(H2)] (\emph{uniform contraction}) for any $y\in B_{\rho}(0)$, $\mu'<\|Dh_n(y)\|_{\textit{op}} <\mu $,
    \item[(H3)] (\emph{bunching}) there is $N>0$ such that 
    $\|Dh^n(0)\|_{\textit{op}}\cdot \|(Dh^n(0))^{-1}\|_{\textit{op}}<\mu^{-n\alpha} $ for all  $n\geq N$.
\end{itemize}
For the convenience and to make the presentation more straightforward, we may replace the hypothesis (H3) above with the following slightly stronger one throughout the paper.
\begin{itemize}
 \item[(H3$'$)] $\|Dh_n(0)\|_{\textit{op}}\cdot \|(Dh_n(0))^{-1}\|_{\textit{op}} \leq \kappa<\mu^{-\alpha}$.
\end{itemize}

\begin{theorem}[Convergence of geometries] \label{thm: convergence of geometries} 
    Let $C,\rho>0$ and $\mu',\mu,\alpha\in (0,1)$ and $\kappa\geq 1$ be real numbers. Let   $\{f_i\}_{i\in \N}$ be a sequence in $\diffloc{1+\alpha}(\mathbb{R}^d)$ 
    satisfying 
    (H0)-(H3).
    Then, the sequence $\{(Df^n (0))^{-1}\circ f^n\}_{n\in \N} $ converges 
    in $\diffloc{1+\alpha}(\mathbb{R}^d)$ to some $F:B_{\rho}(0)\to \R^d$. Moreover, the convergence is uniformly exponential for every sequence 
    satisfying (H0)-(H3) and
    the limit depends uniformly $\cC^1$-continuously on the sequence $\{f_n\}_{n\in \N}$. 
        \end{theorem}
        
Here, by the uniform continuity of the limit we mean that for any $\varepsilon>0$ there is $\delta>0$ such that given two  sequences $\{f_n\}_{n\in \N}$ and $\{\tilde{f}_n\}_{n\in \N}$ satisfying (H0)-(H3) with $\max \{d_{\cC^{1+\alpha}} (f_n,\Tilde{f}_n): n=1,2,\dots\} <\delta$, then their limits are $\varepsilon$-close  in $\diffloc{1+\alpha}(\mathbb{R}^d)$.

In further sections we will use this theorem to obtain infinitesimal properties of regular Cantor sets. Using $\cC^0$-convergence in this theorem one can deduce the following control of shape property.
\begin{corollary}[Control of shape] \label{cor: control-of-shape}
    Let $C,\rho>0$ and $\mu',\mu,\alpha\in (0,1)$ and $\kappa\geq 1$ be real numbers. Then there exist $\xi_0 > 0$ and $\eta_1 > 1>\eta_2$ such that for every sequence $\{f_i\}_{i\in \N}$ in $\diffloc{1+\alpha}(\R^d)$ satisfying 
    (H0)-(H3),  $\xi \in (0, \xi_0]$ and $n \in \mathbb{N}$,
    \begin{equation}
        \eta_2\cdot L_n(B_{\xi}(0)) \subseteq f^n(B_{\xi}(0)) \subseteq \eta_1 \cdot L_n(B_{\xi}(0)),
    \end{equation}
    where  $L_{n}:= Df^n(0)$ is the derivative of $f^n$ at $0$.
\end{corollary}
\begin{figure}[ht]
    \includegraphics[width=.75\textwidth]{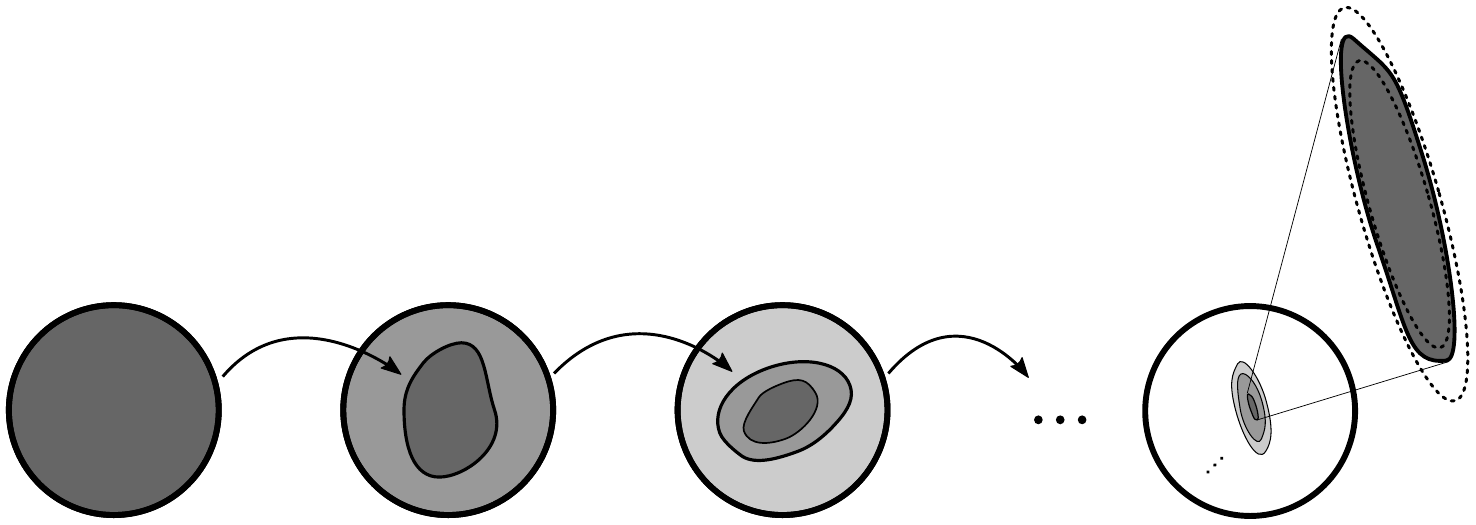}
    \caption{Sequence of maps satisfying (H0)-(H3) in Corollary \ref{cor: control-of-shape}, where the iterations of small balls remain almost ellipsoids. 
    }
     \label{fig:iteration-balls}
     \end{figure}

\begin{remark}
This corollary extends the ``control of shape'' in \cite[Theorem 4.5]{FNR-AiM}, where the authors
assume quasi-conformality of the sequence $\{f^n\}_{n\in \N}$ at the origin. Indeed, their quasi-conformality assumption implies our assumptions (H3) and (H3$'$) for the sequence $\{g_n\}_{n\in \N}$, where 
$g_n : = f_{nk}\circ \cdots \circ f_{nk-k+1}$ for a sufficiently large integer $k$. 
\end{remark}

\begin{remark} The condition (H2) is always satisfied in the study of the geometry of regular Cantor sets. However, within the proof of Theorem \ref{thm: convergence of geometries} it turns out that we could consider the following more general condition (H2$'$) instead of (H2):
    \begin{itemize}
        \item[(H2$'$)] there is a sequence of numbers $\{\mu'_n\}_{n=1}^{\infty}$ in $(0,1)$ with $\lim_{n\to \infty}\frac{1}{n}\ln \mu'_n=0$ such that for all $n\in \N$ and any $y\in B_{\rho}(0)$, $\mu'_n<\|Dh_n(y)\|_{\textit{op}} <\mu $. 
    \end{itemize}
\end{remark}
\begin{proof}[Proof of Theorem \ref{thm: convergence of geometries}]
    Let $r_n$ be the diameter of $U_n$, where $U_n: =f^n\left(B_{\rho}(0)\right)$ for $n\in \N$ and $U_0:= B_{\rho}(0)$. Since $f_n$ is a contracting map with uniform contraction rate $\mu<1$ on $U_{n-1}$, we have 
    $\diam(U_n)\leq \mu \cdot \diam(U_{n-1})$. Therefore, $r_n \leq \rho\cdot \mu^n$. Let $A_k = Df_k(0)$, then
$Df^n(0) = A_n \circ A_{n-1}\circ \cdots \circ A_1.$
We can define
\begin{equation}\label{Control of shape 3 Fn deffinition}
    F_n := (Df^n(0))^{-1}\circ f^n = A_1^{-1}\circ \cdots \circ A_n^{-1} \circ f_n \circ \cdots \circ f_1.
\end{equation}
We have to prove that $\{F_n\}_{n\geq 1}$ is convergent in $\cC^{1+\alpha}$ topology on $B_{\rho}(0)$. 
 For any $z, v \in B_{\rho}(0)$ such that  $z+v \in B_{\rho}(0)$ we have 
\begin{equation}\label{Control of Shape 1}
    \left| f_n(z + v) - \left( f_n(z) + Df_n(z) \cdot v \right) \right| < C\cdot |v|^{1+\alpha}.
\end{equation}
Consequently, $f_n$ is $ C\cdot r_{n-1}^{1+\alpha}$ close to the map $A_n$
in $\cC^0$ norm on $U_{n-1}$. Thus, if $n$ is large enough by (\ref{Control of Shape 1}) we have
\[ r_n \leq \|A_n\|_{\textit{op}} \cdot r_{n-1} + C \cdot r^{1+\alpha}_{n-1} \leq r_{n-1} \cdot \big(\|A_n\|_{\textit{op}} + C_1 \cdot \mu^{(n-1)\alpha}\big), \]
where $C_1 := C\cdot \rho^{\alpha}$. Arguing by induction, using the inequality $\log (x + y) \leq \log x + \frac{y}{x}$ (for positive real numbers $x,y$), we obtain
\[ \log r_n \leq \log \prod_{k=1}^{n}\|A_k\|_{\textit{op}} + C_{0}\cdot \sum_{j=0}^{n-1} (C_1 \cdot \mu^{j\alpha}) \leq
\log \prod_{k=1}^{n}\|A_k\|_{\textit{op}} + C_2, \]
where $C_0 := \sup_{i\geq 1} \|Df_i(0)\|_{\textit{op}}^{-1}<1/{\mu'}$ and $C_2 = \dfrac{C_0\cdot C_1}{1-\mu^{\alpha}}$. So
$r_n \leq C_3 \cdot \prod_{k=1}^{n}\|A_k\|_{\textit{op}},$
where $C_3:= \exp(C_2)$. Analogously using the inequality 
$\log (x-y) \geq \log x- \frac{y}{x-y}$ (for $0<y<x$), we can show, possibly by enlarging $C_3$ , that
\[C_3^{- 1} \cdot \prod_{k=1}^{n}m(A_k) \leq r_n \leq C_3 \cdot \prod_{k=1}^{n}\|A_k\|_{\textit{op}},\]
which means that the size of $f^n(B_{\rho}(0))$ is controlled. This with (\ref{Control of Shape 1}) implies that
\begin{equation} \label{Control of shape 2 gn}
\|g_n\|_{\cC^0} = \left\| f_n - A_n \right\|_{\cC^0} \leq 
C\cdot r_{n-1}^{1+\alpha}\leq C_4 \prod_{k=1}^{n-1}\|A_k\|_{\textit{op}}^{1+\alpha},
\end{equation}
where $g_k$ is $f_k - A_k$ and $C_4 := C\cdot C_3^{1+\alpha}$ and the $\cC^0$ norm is calculated inside the domain $f^n(B_{\rho}(0))$. By (\ref{Control of shape 3 Fn deffinition}) and (\ref{Control of shape 2 gn}) we have 
\begin{equation*}
F_n = \bigg(\id + A_1^{-1}\circ \cdots \circ A_n^{-1} \circ g_n \circ 
A_{n-1} \circ \cdots \circ A_1 \bigg)\circ F_{n-1}.    
\end{equation*}
    Therefore, $F_{n}$ is a sequence of composition of functions 
$$F_{n}=(\id+S_{n})\circ (\id+S_{n-1})\circ \cdots \circ (\id+S_{1}),$$
where 
$S_{n} := A_1^{-1}\circ \cdots \circ A_n^{-1} \circ g_n \circ 
A_{n-1} \circ \cdots \circ A_1.$
We will observe that $S_n$ has exponentially small $\cC^{1+\alpha}$ norm. We have
\begin{align}\label{Control of shape C0 norm}
    \|S_{n}\|_{\cC^0}
    &\leq \prod_{k=1}^{n} \|A_k^{-1}\|_{\textit{op}} \cdot \|g_n\|_{\cC^0} \notag\\ 
    &  \leq  C_4\cdot \prod_{k=1}^{n} \|A_k^{-1}\|_{\textit{op}} \cdot \prod_{k=1}^{n-1}\|A_k\|_{\textit{op}}^{1+\alpha} \notag\\
    &\leq C_4\cdot C_0^{1+\alpha}\cdot \kappa^{n}\cdot \prod_{k=1}^{n}\|A_k\|_{\textit{op}}^{\alpha} \notag\\
    &\leq C_{5} \cdot (\kappa \cdot \mu^{\alpha})^n,
\end{align}
where $C_5 := C_4\cdot C_0^{1+\alpha}$.
We also have the following estimate for $\|DS_n\|_{\it op}$ on $B_{\rho}(0)$.
\begin{align}\label{Control of shape C1 norm}
    \|DS_{n}\|_{\textit{op}}
    & = \|A_1^{-1}\circ \cdots \circ A_n^{-1} \circ Dg_n \circ A_{n-1} \circ \cdots \circ A_1\|_{\textit{op}}\notag\\
    & \leq \prod_{i= 1}^{n} \|A_i^{-1}\|_{\textit{op}} \cdot \|Dg_n\|_{\textit{op}} \cdot \prod_{i= 1}^{n-1} \|A_i\|_{\textit{op}} \notag\\
    &\leq C_0\cdot \kappa^n \cdot \|Dg_n\|_{\textit{op}} \notag\\
    &\leq C_6\cdot (\kappa \cdot \mu^{\alpha})^n, 
\end{align}
where $C_6 := C_0\cdot C\cdot (\rho\mu^{-1})^{\alpha}$ and the last inequality comes from $\alpha$-H\"older regularity of $Df_n$, that for any $z\in U_{n-1}$ we have
$$\|Dg_n(z)\|_{\textit{op}} = \|Df_n(z)-Df_n(0)\|_{\textit{op}} \leq C\cdot r_{n-1}^{\alpha} \leq C\cdot(\rho\mu^{-1})^{\alpha}\cdot\mu^{n\alpha}.$$
We can also bound the $\alpha$ - H\"older semi norm $|DS_n|_{\alpha}$ as below
\begin{align*}
|DS_n|_{\alpha} &= \sup_{a_0,b_0 \in U_0} \dfrac{\|DS_n(a_0)-DS_n(b_0)\|_{\textit{op}}}{|a_0-b_0|^\alpha}&
\\
&= \sup_{a_0,b_0 \in U_0}  \dfrac{1}{|a_0-b_0|^\alpha}\cdot \|A_1^{-1}\circ \cdots \circ A_n^{-1}\\
&\qquad \qquad \qquad \qquad\qquad
\circ (Dg_n(a_{n-1})-Dg_n(b_{n-1})) \circ A_{n-1} \circ \cdots \circ A_1\|_{\textit{op}}
\\
&\leq \sup_{a_0,b_0 \in U_0}  \dfrac{1}{|a_0-b_0|^\alpha}\cdot \prod_{i= 1}^{n} \|A_i^{-1}\|_{\textit{op}} \cdot \|Dg_n(a_{n-1})-Dg_n(b_{n-1})\|_{\textit{op}} \cdot \prod_{i= 1}^{n-1} \|A_i\|_{\textit{op}}
\\
&\leq \sup_{a_0,b_0 \in U_0}  \dfrac{C_0\cdot \kappa^n}{|a_0-b_0|^\alpha}\cdot \|Dg_n(a_{n-1})-Dg_n(b_{n-1})\|_{\textit{op}} 
\\
&= \sup_{a_0,b_0 \in U_0}  \dfrac{C_0\cdot\kappa^n}{|a_0-b_0|^\alpha}\cdot  \|Df_n(a_{n-1})-Df_n(b_{n-1})\|_{\textit{op}}
\\
& \leq \sup_{a_0,b_0 \in U_0}  \dfrac{C_0\cdot\kappa^n}{|a_0-b_0|^\alpha}\cdot C\cdot |a_{n-1}-b_{n-1}|^\alpha 
\\
&= C_0\cdot C\cdot \sup_{a_0,b_0 \in U_0} \big(\dfrac{|a_{n-1}-b_{n-1}|}{|a_0-b_0|}\big)^\alpha\cdot \kappa^n,
\end{align*}
where $a_{n-1} = A_{n-1}\circ \cdots \circ A_1(a_0)$ and $b_{n-1} = A_{n-1}\circ \cdots \circ A_1(b_0)$. So 
$$\dfrac{|a_{n-1}-b_{n-1}|}{|a_0-b_0|}\leq \dfrac{1}{|a_0-b_0|}\cdot \prod_{i\geq 1}^{n-1} \|A_i\|_{\textit{op}} \cdot |a_0-b_0| \leq \mu^{n-1}$$
which implies that
\begin{equation}\label{Control of shape alpha norm}
    |DS_n|_{\alpha}\leq C_7\cdot (\kappa\cdot \mu^\alpha)^n,
\end{equation}
where $C_7:=C_0\cdot C\cdot \mu^{-1}$. Therefore, by (\ref{Control of shape C0 norm}), (\ref{Control of shape C1 norm}) and (\ref{Control of shape alpha norm}) we conclude that $\|S_n\|_{\cC^{1+\alpha}} \leq C_8\cdot (\kappa\cdot \mu^\alpha)^n$ where $C_8 := C_7+C_6+C_5$.
Bunching assumption implies that $\kappa\cdot \mu^{\alpha}<1$, so  the series $\sum_{n\geq 1} \|S_{n}\|_{\cC^{1+\alpha}}$ converges.  Thus, by Lemma \ref{lem: appendix- convergence lemma} the sequence $\{F_n\}_{n\in\N}$ of compositions 
converges in $\diffloc{1+\alpha}(\R^d)$ to some $F:B_{\rho}(0)\to \R^d$.

All the constants appearing above depend continuously on $C,C_0,\mu,\kappa$. Indeed, $C$ is the $\cC^{1+\alpha}$ norm bound for the sequence $\{f_n\}_{n\in \N}$ and $C_0 = \sup_{i\geq 1} \|(Df_i(0)\|_{\textit{op}}^{-1}<1/{\mu'}$ so they depend continuously on the sequence $\{f_n\}_{n\in \N}$.  Thus, for any other sequence $\{\tilde{f}_n\}_{n\in \N}$ sufficiently close to $\{f_n\}_{n\in \N}$ in $\diffloc{1+\alpha}(\R^d)$, all of those estimates would be the same
except with a minor pre-fixed error. This implies that the convergence is uniform for any sequence of functions $\{\tilde{f}_n\}_{n\in \N}$ satisfying (H0)-(H4) which admits same constants $C,C_0$. Finally these observation together with Lemma \ref{lem: appendix continuity second version} implies that $F$ depends continuously on the sequence $\{f_n\}_{n\in \N}$.
\end{proof}

\subsection{Strong bounded distortion} 

\begin{theorem}\label{thm: Bounded distortion limit geom improved}
Given the assumptions of the Theorem \ref{thm: convergence of geometries}, there is a positive constant $C>1$ such that for any $x,y\in B_{\rho}(0)$ and all $n\geq 1$,
$$C^{-1}<
\dfrac{\|Df^n(x)\|_{\textit{op}}}{\|Df^n(y)\|_{\textit{op}}},\;
\dfrac{m\big(Df^n(x)\big)}{m\big(Df^n(y)\big)},\;
\dfrac{\operatorname{det}Df^n(x)}{\operatorname{det} Df^n(y)}
< C.$$
The closer $x$ and $y$ are to each other, the closer the above quantities will be to 1.
\end{theorem}
\begin{proof} 
We prove that the sequence of matrices 
    $$P_m(y,x) = Df^m(y)^{-1} \circ Df^m(x) $$ is convergent and the limit matrix $P(y,x)$ is non-singular and continuous in $y,x$. We can write
    $$
    Df^m(x)= \tilde A_{m}\circ \cdots \circ \tilde A_{1},\;\;
    Df^m(y)= \tilde B_{m}\circ \cdots \circ \tilde B_{1},    
    $$
    where $\tilde A_{k},\tilde B_k$ are linear maps equal to derivatives of $f_m$ at the points $x_{m-1}:=f^{m-1}(x),y_{m-1}:=f^{m-1}(y)$, respectively. We have
    \begin{equation}\label{eq distortion corol naming Pn(y,x)}
        P_m(y,x) = Df^{m}(y)^{-1} \circ Df^{m}(x) 
    = \tilde B_{1}^{-1}\circ \cdots \circ \tilde B_{m}^{-1} \circ \tilde A_{m}\circ \cdots \circ \tilde A_{1}.
    \end{equation}
Openness of bunching assumption which is defined in the origin implies that it is extendable to a small neighborhood around the origin. Indeed, because $x_k,y_k\to 0$ as $k\to \infty$ we may assume the bunching conditions (H3) or (H3$'$) for the sequence of linear maps 
$\{B_m\}_{m\geq N}$ for some $N$.
Hence, one can assume that there is $\Tilde{\kappa}$ such that $\Tilde{\kappa}\cdot \mu^{\alpha}<1$ and for any $i>N$,
$ \|\tilde B_{i}\|_{\textit{op}}\cdot \|\tilde B_{i}^{-1}\|_{\textit{op}}\leq \Tilde{\kappa}.$
Let $r_n$ be the diameter of $f^n\left(B_{\rho}(0)\right)$. Analogous with the proof of Theorem \ref{thm: convergence of geometries} we infer that the sequence of matrices 
$X_k:= \tilde A_k-\tilde B_k$ has small norm $\|X_{k}\|_{\textit{op}}\leq C_1\cdot |x-y|^\alpha \cdot \mu^{k\alpha}$ 
 where $C_1=\rho \cdot C$. Also, we can write $P_m(y,x)$ as a sequence of products of matrices
$$P_m(y,x)=(\id+S_{m})\circ (\id+S_{m-1})\circ \cdots \circ (\id+S_{1}),$$
where 
$S_m$ has small norm 
$\|S_{m}\|_{\textit{op}}\leq  C_{2}\cdot(\Tilde{\kappa} \cdot \mu^{\alpha})^{m}\cdot |x-y|^\alpha$ for some constant $C_2$.
Since $\Tilde{\kappa}\cdot \mu^{\alpha}<1$ we have $\sum_{m>N} \|S_{m}\|_{\textit{op}}< C_1\cdot |x-y|^\alpha\cdot (1-\Tilde{\kappa}\cdot \mu^{\alpha})^{-1}$.
Thus the series $\sum_{m>N} \|S_{m}\|_{\textit{op}}$ converges and so via Lemma \ref{lem: appendix- convergence lemma} the infinite product 
$\prod_{m>N} (\id+S_m)$
converges to a matrix $Q(y,x)$ which give us the convergence of $P_m(y,x)$.
Continuity part of the Lemma \ref{lem: appendix- convergence lemma} with the upper bound of $\sum_{m>N} \|S_{m}\|_{\textit{op}}$ implies that this infinite product converges to $\id$ matrix as $|x-y|\rightarrow 0$. In addition, $\cC^1$ continuity of $g$ implies that the finite product $\prod_{i=1}^{N}(\id+S_i)$ converges to $\id$ as $|x-y|\rightarrow 0$. 
This with the fact that 
$P(y,x)\circ P(x,z) = P(y,z)$ for all $x,y,z\in G(a)$, gives the continuity of $P(y,x)$ in $(x,y)$.
Therefore, there is a constant $\Bar{C}>1$ such that for all $m\in \mathbb{N}$ and any $x,y\in B_{\rho}(0)$
$$\Bar{C}^{-1} < \|P_m(x,y)\|_{\textit{op}},\; m\big(P_m(x,y)\big),\; |P_m(x,y)|<\Bar{C}.$$
\end{proof}

\subsection{Sequence of geometries}\label{sec: sequence of geometries}
In this section we study the geometry of infinitesimal pieces of regular Cantor sets generalizing the results in \cite{MY}.
We use the results of previous subsections to show that the sequences of normalized geometries associated to a Cantor set  
have a well-defined limit, provided that the bunching condition holds. The limit geometries (or the limit charts) are the main tools to study the infinitesimal geometry of regular Cantor sets.
Given a bunched Cantor set $K$, denote $K(a)= K\cap G(a)$ and fix a base point $c(a)\in K(a)$  for any $a\in \mathfrak{A}$. Additionally, given 
$\underline{\theta}=(\dots, \theta_{-2},\theta_{-1},\theta_{0})\in \Sigma^{-}$, we write 
$\underline{\theta}_n=(\theta_{-n},,\dots ,\theta_{0})$ and $r_{\underline{\theta}_n}:= \diam(G^*(\underline{\theta}_n))$.
Given $\underline{\theta} \in \Sigma^-$ and $n \geq 1$,
If $\theta_{0}=a$ we can define $c_{\underline{\theta}_0} := c(a) $ and
$$c_{\underline{\theta}_n} := f_{\underline{\theta}_n}(c_{\underline{\theta}_0})=f_{(\theta_{-n},\theta_{-n+1})}\circ f_{(\theta_{-n+1},\theta_{-n+2})} \circ \cdots \circ f_{(\theta_{-1},\theta_{0})}(c_{\underline{\theta}_0}),$$
$$k_n^{\underline{\theta}}:= A_{\underline{\theta}_{n}}^{-1}\circ f_{\underline{\theta}_{n}}:G^*(\theta_0)\rightarrow \R^d$$
called normalized geometries where  \(A_{\underline{\theta}_{n}}\) is an affine map such that 
$$DA_{\underline{\theta}_{n}} = Df_{\underline{\theta}_{n}}(c_{\underline{\theta}_{0}}) \;\text{and} \;A_{\underline{\theta}_{n}}(c_{\underline{\theta}_{0}}) = c_{\underline{\theta}_{n}}.$$
Observe that, $k_n^{\underline{\theta}}(\theta_0) = \theta_0$, $Dk_n^{\underline{\theta}}(\theta_0) = \id$ and 
\begin{equation}\label{eq Affine definition limit geometry}
A_{\underline{\theta}_{n}} = A_{n}^{\underline{\theta}}\circ A_{n-1}^{\underline{\theta}}\circ \cdots \circ A_{1}^{\underline{\theta}}
\end{equation}
where for any $j$,  $A_{j}^{\underline{\theta}}$ is
the affine estimate of $f_{(\theta_{-j},\theta_{-j+1})}$ at the point $c_{\underline{\theta}_{j-1}}$. 
Therefore, $A_{\underline{\theta}_{n}}\in \Aff(d,\R)$. 
Analogously, in the case that $K$ is a bunched holomorphic
Cantor set in $\R^{2d}$ then $A_{\underline{\theta}_{n}}\in \Aff(d,\C)\subset \Aff(2d,\R)$.
This is because derivatives of the map $g$ at the points $c_{\underline{\theta}_{j-1}}$ are  matrices in $\GL(d,\C) \subset \GL(2d,\R)$.

Inspired by \cite{Sullivan1988}, the existence of limit geometries (or limit charts) is a key ingredient in establishing the covering criterion. In the following lemma, we show that limit geometries exist under the bunching assumption; this result is a direct corollary of Theorem~\ref{thm: convergence of geometries}. Indeed, given $\underline{\theta}\in \Sigma$, under appropriate changes of coordinates the sequence of maps $h_n:=f_{(\theta_{-n},\theta_{-n+1})}$ satisfies hypotheses (H0)-(H3).
 
 \begin{lemma}[Limit geometries]\label{lemma normalized limit geometries convergence} 
     For any $\underline{\theta}\in \Sigma^{-}$ the sequence of maps $\{k_n^{\underline{\theta}}\}_{n\in \N}$ converges in $\diffloc{1+\alpha}(\R^d)$ to a map $k^{\underline{\theta}}:G^*(\theta_0)\mapsto \R^d$. Moreover, the convergence is uniform over $\underline{\theta}\in \Sigma^{-}$.
     and in a small neighborhood of $g$ in $\Omega^{1+\alpha}_{\Sigma,d}$ (or $\Omega_{\Sigma,2d}^{\Hol,1+\alpha}$ for the complex case). 
     The map $ k^{\underline{\theta}}: G^*(\theta_0) \to \R^d $ defined for any $\underline{\theta} \in \Sigma^-$ is called the limit geometries of $K$. 
 \end{lemma}

 \begin{remark}
     Compared with the definitions in \cite{AM-2023}, here it holds that $k^{\underline{\theta}}(c_{\underline{\theta}_0})=c_{\underline{\theta}_0}$, while in that paper one has $k^{\underline{\theta}}(c_{\underline{\theta}_0})=0$. Note that the maps $A_{\underline{\theta}_{n}}$ and $k_n^{\underline{\theta}}$ depend on the base point $c({\theta}_0)$; consequently, the limit geometry $k^{\underline{\theta}}$ also depends on it. However, if we vary the base point $c({\theta}_0)$ to $c'({\theta}_0)$ the resulting limit geometry changes up to a composition with a bounded affine map. This is a consequence of the convergence of $\{P_m(x,y)\}_{m\in \N}$ in the proof of Theorem \ref{thm: Bounded distortion limit geom improved}.
 \end{remark}

Next, we deduce 
that the limit geometry $k^{\underline{\theta}}$ is $\cC^1$-continuous with respect to the Cantor set $K$ and H\"older continuous with respect to $\underline{\theta}$. Indeed, one can define a metric on the set $\Sigma^-$ for a given Cantor set $K$ that  
\begin{align}\label{eq: metric shift Sigma}
    d(\underline{\theta}',\underline{\theta}) := a^{k},
\end{align}
where 
$a:=\mu\cdot \kappa^{1/\alpha}<1$
and $k\geq 0$ is the the smallest integer that $\theta_{-k} \neq \theta'_{-k}$.
We use the notations from Definition \ref{def: space of cantor sets} in the following lemma.

 \begin{lemma}(Continuity of limit geometries)\label{lem: continuity of limit geoms under perturbation of Cantors} the map $(K,\underline{\theta}) \mapsto k^{\underline{\theta}}$ from $\Omega_{\Sigma,d}^{1+\alpha}\times \Sigma^-$ to ${\rm {Emb}}^{1+\alpha}(G^*(\theta_0),\R^d)$ is continuous. Moreover, 
 
 \begin{itemize}
    \item[(i)]  Given $\underline{\theta}\in\Sigma^-$, for any $\varepsilon>0$ there is $\delta>0$ such that for any $\tilde K\in U_{K,\delta}$, $\tilde k^{\underline{\theta}}$ is $\varepsilon$-close to $k^{\underline{\theta}}$ on $G^*(\theta_0)\cap G'^*(\theta_0)$ in $\cC^{1}$ topology.
     \item[(ii)]  Given $K$, the map $\underline{\theta}\mapsto k^{\underline{\theta}}$  from $\Sigma^-$ to $\bigcup_{a\in \mathfrak{A}} {\rm {Emb}}^{1+\alpha}(G^*(a),\R^d)$ is $\alpha$-H\"older continuous with respect to $\underline{\theta}$ in the metric defined in \eqref{eq: metric shift Sigma}.
 \end{itemize}
 \end{lemma}
 
 \begin{proof} First item follows directly from
      the continuity part of the Theorem \ref{thm: convergence of geometries}. To prove item (ii), we use Remark \ref{remark: appendix speed of convergence}. It implies that for all $n\in \N$ we have  that 
     $$\|k^{\underline{\theta}}-k_n^{\underline{\theta}}\|_{\cC^{1+\alpha}}\leq \hat{C_8}\cdot (\kappa\cdot \mu^{\alpha})^n$$
     for some $\hat{C_8}>0$ depended on $K$. Therefore, for $\underline{\theta}',\underline{\theta}\in \Sigma^-$ that agree on their first $n$ letters we have $k_n^{\underline{\theta}}=k_n^{\underline{\theta}'}$ so
     $$\|k^{\underline{\theta}}-k^{\underline{\theta}'}\|_{\cC^{1+\alpha}}\leq
     2\hat{C_8}\cdot (\kappa\cdot \mu^{\alpha})^n =2\hat{C_8}\cdot d(\underline{\theta}',\underline{\theta})^{\alpha}.$$
     The above estimate shows that the map 
     $(K,\underline{\theta}) \mapsto k^{\underline{\theta}}$ is continuous.
 \end{proof}

\subsection{Infinitesimal geometry of  Cantor sets}
\label{sec: infinitesimal geometry of cantor} 
Let $K$ be a bunched regular  Cantor set in $\R^d$. With the notations fixed in \S \ref{sec: sequence of geometries} we have the following shape control  property, which provides a fine description of the geometry of bunched Cantor sets at any scale.
It is a direct consequence of Corollary \ref{cor: control-of-shape}.

\begin{theorem}[Bounded geometry of bunched Cantor sets]\label{thm: Bounded geometry of regular Cantor sets}
There are positive constants $\gamma_1<1<\gamma_2$ and $\xi_0>0$ such that for any $n\geq 1$ and any $x\in K$ and $\xi \in [0,\xi_0]$ we have
\begin{equation}\label{eq: Bounded geometry of bunched Cantors}
A_{\underline{\theta}_{n}}\big(B_{\gamma_1\cdot \xi}(x)\big) \subseteq f_{\underline{\theta}_{n}}\big(B_{\xi}(x)\big)\subseteq A_{\underline{\theta}_{n}}\big(B_{\gamma_2\cdot \xi}(x)\big).
\end{equation}
\end{theorem}
In particular, the diameter of $G^*(\underline{\theta}_{n})$, denoted by $r_{\underline{\theta}_{n}}$, is of order $\|Df_{\underline{\theta}_{n}}(c_{\underline{\theta}_{0}})\|_{\textit{op}}$. More precisely, there are positive constants $\gamma_1'< 1< \gamma_2'$ such that for any $n\in \N$ and any base point $c(\theta_0)\in K(\theta_0)$
\begin{equation}\label{eq control of diameter of nth image by derivatives}
    \gamma_1' \cdot \|Df_{\underline{\theta}_{n}}(c(\theta_0))\|_{\textit{op}}\leq r_{\underline{\theta}_{n}} \leq \gamma_2' \cdot \|Df_{\underline{\theta}_{n}}(c(\theta_0))\|_{\textit{op}}.
\end{equation}
Similar statement also holds for the inner radius of $G^*(\underline{\theta}_{n})$.

\section{Configurations and renormalizations}\label{sec: config and renormalization}
The renormalization operators are fundamental tools to study intersection of regular Cantor sets, introduced by Moreira-Yoccoz \cite{MY} for Cantor sets in dimension one. 
In this section, we see that the results of the previous section allow one to extend the Moreira-Yoccoz's method for the study of bunched Cantor sets in arbitrary dimensions both in real and complex spaces.  
\subsection{Configurations of regular Cantor sets}
The bunching condition \eqref{eq: definition of bunched cantor sets} is invariant under (smooth) conjugation, possibly with bigger $N_g$. Therefore, if \(h : \bigsqcup_{a \in \mathfrak{A}} G(a) \to U \subset \R^d\) is a \(\cC^{1+\alpha}\) diffeomorphism then \(h(K)\) is a bunched Cantor set with the generator $\Tilde{g} = h \circ g \circ h^{-1}$ and pieces $\Tilde{G}(a) := h(G(a))$.
We name such a re-embedding of $K$ in $\R^d$ as a \emph{configuration}. 
Given a piece $G(a)$ for $a\in \mathfrak{A}$, such a $\cC^{1+\alpha}$ diffeomorphism $h: G(a)\to h(G(a))\subset \R^d$ is called \emph{configuration of the piece $G(a)$} of Cantor set $K$.
  We write \(\mathcal{P}(a)\) for the space of all \(\cC^{1+\alpha}\) configurations of the piece \(G(a)\) equipped with the \(\cC^{1+\alpha}\) topology and we denote
  \begin{equation}
   \cP: =  \bigcup_{a\in \mathfrak{A}} \cP(a).  
  \end{equation}
  In the holomorphic case we denote the space of configurations as $\cP^{\rm Hol}$ containing configurations $h$ that $Dh(x)\in \GL(d,\C)\subset \GL(2d,\R)$ for any $x\in K$.
  
    If $h\in \cP(a)$ (or $\cP^{\rm Hol}$) is an affine map in
  $\Aff(d,\R)$ (or $\Aff(d,\C)$) on its domain $G(a)$, then we call it an \emph{affine configuration}. If $h = A\circ k^{\underline{\theta}}$ for some $A\in \Aff(d,\R)$ (or $\Aff(d,\C)$) and $\underline{\theta}\in \Sigma^-$, we call it an \emph{affine configuration of limit geometry}.

The space of all configurations is a function space equipped with $\cC^{1+\alpha}$ topology which allows us to analyze linear and non-linear re-positionings of $K$ in $\R^d$.  
For each $\ua=(a_0,a_1,\dots,a_n)\in \Sigma^{\mathrm{fin}}$ we define the map 
$T_{\ua}:\cP(a_0)\to \cP(a_n)$
with 
$$T_{\ua}(h):=h\circ f_{\ua}.$$
This definition implies that 
$$T_{\ua}= T_{(a_0,a_1)}\circ T_{(a_1,a_2)}\circ \cdots \circ T_{(a_{n-1},a_n)}.$$
Hence, 
this semi-pseudogroup of operators $\cT:=\{T_{\ua}:\ua\in \Sigmafin\}$ is generated by the finite family of operators 
$
\cT_1: = \{T_{(a,b)}:\; (a,b)\in \mathfrak B\}.    
$

Recall that given a bunched Cantor set $K$
 it follows from Theorem \ref{thm: Bounded geometry of regular Cantor sets} that the infinitesimal geometries of $K$ are approximately affine. This allows us to implement the strategy of Moreira-Yoccoz. To do so, we first define the space of \emph{representation of affine configurations of limit geometries} \begin{equation}
     \cA := \Aff(d,\R) \times \Sigma^-
 \end{equation}
 with the continuous map 
 \begin{align}
 \Phi: \cA &\rightarrow \cP    \notag\\
 (A,\underline{\theta})&\mapsto A\circ k^{\underline{\theta}}.
 \end{align}  
 For $K$ being a  bunched holomorphic Cantor set in $\R^{2d}$  we will consider the space 
 \begin{equation}
     \cA^{\Hol}: = \Aff(d,\C) \times \Sigma^-
 \end{equation}
 The following lemma shows that the family $\cT$ induces an action on the space $\cA$ (or $\cA^{\rm{Hol}}$).
 It is indeed a natural extension of \cite[Lemma 3.6]{AM-2023} to higher dimensions with a similar proof, thanks to the convergence results in \S \ref{sec: sequence of geometries}-\ref{sec: infinitesimal geometry of cantor}.

\begin{lemma} \label{lem: affineness of action of renormaliztions on limit geoms} Given $\ua= (a_0,\dots,a_n)\in \Sigma^{\rm{fin}}$ and $\underline{\theta}\in \Sigma^-$ with $a_0=\theta_0$, then 
\[F^{\underline{\theta}\ua}:= k^{\underline{\theta}}\circ f_{\ua}\circ (k^{\underline{\theta}\ua})^{-1}\in \Aff(d,\R),\]
and the following diagram is commutative.
\begin{equation}\label{diagram}
     \begin{tikzcd}
\cA \arrow{r}{\Phi} \arrow[swap]{d}{T^{*}_{\ua}} & \cP \arrow{d}{T_{\ua}} \\%
\cA \arrow{r}{\Phi}& \cP
\end{tikzcd}
\end{equation}
where $T^{*}_{\ua}:\cA\to \cA$ is a (H\"older) continuous map defined by
\begin{equation}\label{eq: hat T renormalization}
    T^{*}_{\ua}(A,\underline{\theta}) := (A \circ F^{\underline{\theta}\ua},\underline{\theta}\ua).
\end{equation}
In the holomorphic case, \(F^{\underline{\theta}\ua}\) is in $\Aff(d,\C)\subset \Aff(2d,\R)$.
\end{lemma}
\begin{proof}
We have
\begin{align}
    k^{\underline{\theta}}\circ f_{\ua}\circ (k^{\underline{\theta}\ua})^{-1}&= 
\lim_{m\to \infty} k_m^{\underline{\theta}}\circ f_{\ua}\circ (k_{m+n}^{\underline{\theta a}})^{-1}\notag\\
& = \lim_{m\to \infty} A^{-1}_{\underline{\theta}_{m}}\circ f_{\underline{\theta}_m}\circ f_{\ua} \circ (A^{-1}_{\underline{\theta a}_{m+n}}\circ f_{\underline{\theta a}_{m+n}})^{-1}\notag\\
& = \lim_{m\to \infty} A^{-1}_{\underline{\theta}_m}\circ A_{\underline{\theta a}_{m+n}}\notag
\end{align}
which is a convergent sequence in the closed space of affine maps. Therefore, $F^{\underline{\theta}\ua}$ is an affine map which is invertible by definition.
For $(A,\underline{\theta})\in \cA$, we have
\begin{align}
T_{\ua}\circ \Phi\left((A, \underline{\theta})\right)
&= T_{\ua}(A\circ k^{\underline{\theta}})
\notag\\&= A\circ k^{\underline{\theta}}\circ f_{\ua} \notag\\
&= A\circ F^{\underline{\theta}\ua} \circ k^{\underline{\theta a}}\notag\\
& = \Phi \circ T^{*}_{\ua} \left((A, \underline{\theta})\right).\notag
\end{align}

H\"older continuity of $T^{*}_{\ua}$ follows directly from Lemma \ref{lem: continuity of limit geoms under perturbation of Cantors}.
\end{proof}
\begin{lemma}\label{lem: continuity of renormalization operators} Given $\ua\in \Sigmafin$,
   the operator $T^{*}_{\ua}$ varies continuously with respect to the Cantor set $K\in \Omega_{\Sigma,d}^{1+\alpha}$.  
\end{lemma}
\begin{proof}  
According to Lemma \ref{lem: continuity of limit geoms under perturbation of Cantors} the relation $F^{\underline{\theta a}}= k^{\underline{\theta}}\circ f_{\ua}\circ (k^{\underline{\theta a}})^{-1}$ implies that
given $\underline{\theta}\in \Sigma^-$ and $\ua\in \Sigmafin$, the association 
$K \mapsto F^{\underline{\theta a}}$ is continuous with respect to $K$. In particular, for any $\varepsilon>0$ there is $\delta>0$ such that for any $\tilde K\in U_{k,\delta}$ we have that 
$\|F^{\underline{\theta a}}-\tilde F^{\underline{\theta a}}\|_{\it op}<\varepsilon$. This implies the continuity.
\end{proof}
Consequently, the space of affine configurations of limit geometries $\Phi(\cA)\subset \cP$ is invariant under the action of the family $\cT$ and the action is continuous with respect to $K$.
The commutative diagram \eqref{diagram} shows that the family $\cT^{*}:=\{T^{*}_{\ua}:\ua \in \Sigmafin\}$ of operators acting on $\cA$ (or $\cA^{\rm Hol}$) is generated by the finite family
  $  \cT^{*}_1:=\{T^{*}_{(a,b)}:(a,b)\in \mathfrak{B}\}$
which implies that
$$F^{\underline{\theta}\ua}=
F^{\underline{\theta} a_{1}} \circ F^{\left(\underline{\theta} a_{1}\right) a_{2}} \circ \cdots \circ F^{\left(\underline{\theta} a_{1} \cdots a_{n-1}\right) a_{n}}.$$
Assuming $A=\id$ in \eqref{eq: hat T renormalization}, we deduce that
\begin{equation}\label{eq affine action of renormalization operator with word an}
    F^{\underline{\theta a}} \circ k^{\underline{\theta a}} = k^{\underline{\theta}} \circ f_{\ua}.
\end{equation}
Therefore, $\diam(k^{\underline{\theta}} \circ f_{\ua}(G(a_n)))$ is of order $\|DF^{\underline{\theta a}}\|_{\textit{op}}$.
On the other hand, since $k^{\underline{\theta}}$ is uniformly bounded for all $\underline{\theta}\in \Sigma^-$ and $G(\ua)=f_{\ua}(G(a_n))$, we get that $\diam(k^{\underline{\theta}} \circ f_{\ua}(G(a_n)))$ is of order $\diam(G(\ua))$. In addition, by $F^{{\underline{\theta a}}}= k^{\underline{\theta}}\circ f_{\ua}\circ (k^{\underline{\theta}\ua})^{-1}$ we deduce that $DF^{{\underline{\theta a}}}$ and $Df_{\ua}$ have uniformly comparable norms.  These observations give us the following corollary.

\begin{corollary} \label{cor:operator bound of renormalization linear part}
     Given a bunched (or bunched holomorphic) Cantor set $K$, there is a constant $C>1$ such that for any $\underline{\theta} \in \Sigma^{-}$ and 
     $\ua = (a_0,\dots,a_n) \in \Sigmafin$ with $a_0=\theta_0$,
    \begin{align}
    C^{-1}\leq\; \dfrac{\|DF^{\underline{\theta a}_{n}}\|_{\textit{op}}}{\|Df_{\ua_{n}}\|_{\textit{op}}},\;
    \dfrac{\|(DF^{\underline{\theta a}_{n}})^{-1}\|_{\textit{op}}}{\|Df_{\ua_{n}}^{-1}\|_{\textit{op}}},\; \dfrac{\|DF^{\underline{\theta a}_{n}}\|_{\textit{op}}}{\diam(G(\ua))}\; \leq C.\notag
    \end{align}    
\end{corollary}

\subsection{Normalized configurations}
  
Since the maps $f_{\ua}$ are all contracting, in order to understand the action of the maps $T_{\ua}$ on the space $\cP$ (or $\cP^{\rm Hol}$) we consider the normalization of configurations.
Given  $h\in \cP(a)$ (or $\cP^{\rm Hol}(a)$) and $p\in K(a)$, we define $A_{h,p}\in \Aff(d,\R)$ (or $\Aff(d,\C)$) as the affine estimate of $h$ at $p$, 
so that $A_{h,p}(p):=h(p)$ and $DA_{h,p}:=Dh(p)$.
We define the \emph{normalization
of $h$ at $p$} by 
\begin{align}\label{eq: definition of normalized configurations}
    \widehat h_p:=A_{h,p}^{-1}\circ h.
\end{align}
 Notice that $\widehat h_p$ fixes the point $p$ and its derivative at $p$ is $\id$. In the following lemma we show that limit geometries are attracting normalization of configurations under the action of family $\cT$.

\begin{lemma} \label{lemma limit geometries are attractor}
 Let $K$ be a bunched (or bunched holomorphic) Cantor set and $h \in \cP\left(a_{0}\right)$ a configuration of a piece in $K$.
 Then, there exists a constant $C_{h,K}>0$ depended continuously on $K$ and $h$ such that for any finite word
 $\ua=(a_0,a_1,\dots,a_n)\in \Sigmafin$,
 $\underline{\theta} \in \Sigma^{-}$with $\theta_{0}=a_{0}$ and a base point $p\in K(a_n)$ we have
$$\| \widehat{T_{\ua}(h)}_{p}- k^{\underline{\theta a}}\|_{\cC^{1+\alpha}}<C_{h,K}\cdot (\kappa\cdot \mu^{\alpha})^n.$$ \end{lemma}
\begin{proof}
By \eqref{eq: definition of normalized configurations} we have to prove that
$$\|A_{f_{\ua},p}^{-1} \circ A_{h,f_{\ua}(p)}^{-1} \circ h \circ f_{\ua}- k^{\underline{\theta a}}\|_{\cC^{1+\alpha}}<C_{h,K}\cdot (\kappa\cdot \mu^{\alpha})^n.$$
We know from Remark \ref{remark: appendix speed of convergence} and uniformity of the convergence of limit geometries $\lim_{m\rightarrow \infty} k_m^{\underline{\theta}}= k^{\underline{\theta}}$ with respect to $\underline{\theta}$ 
that $\|k_n^{\underline{\theta a}}-k^{\underline{\theta a}}\|_{\cC^{1+\alpha}}<\tilde C\cdot (\kappa\cdot \mu^{\alpha})^n$ for some $\tilde C$ depended only on $K$. Thus it is sufficient to prove that there is a constant $\tilde C_1>0$ depended on $K,h$ such that
\begin{equation}\label{eq lemma attractor statement for analysis}
\|A_{f_{\ua},p}^{-1} \circ A_{h,f_{\ua}(p)}^{-1} \circ h \circ f_{\ua}- k_n^{\underline{\theta a}}\|_{\cC^{1+\alpha}}
    <\tilde C_1\cdot (\kappa\cdot \mu^{\alpha})^n.
\end{equation}
Denote   
\begin{align}
F_{\ua}:  & = A_{f_{\ua},p}^{-1} \circ A_{h,f_{\ua}(p)}^{-1} \circ h \circ f_{\ua}- k_n^{\underline{\theta a}}  \notag\\
& = \left(A_{f_{\ua},p}^{-1} \circ A_{h,f_{\ua}(p)}^{-1} \circ h \circ A_{h,f_{\ua}(p)}- \id \right) \circ k_n^{\underline{\theta a}}
\end{align}
where the second equation comes from $k_n^{\underline{\theta a}} = A_{f_{\ua},p}^{-1} \circ f_{\ua}$. By $\cC^{1+\alpha}$- H\"older regularity of $h$, there are constants $C', C''>0$ such on the domain  $G(\ua) = 
    f_{\ua}(G(a_n))$ 
    \begin{align}
        \label{eq holder estimate attractorness of limit geoms 1}
        \|h -A_{h,f_{\ua}(p)}\|_{\cC^{0}}\leq & ~ C'\cdot r_{\ua}^{1+\alpha}\leq C''\cdot r_{\ua}^{\alpha}\cdot \|Df_{\ua}\|_{{\textit{op}}},\\
        \label{eq holder estimate attractorness of limit geoms 2}
        \|Dh - DA_{h,f_{\ua}(p)}\|_{{\textit{op}}} < & ~ C'\cdot r_{\ua}^{\alpha},\\
        \label{eq holder estimate attractorness of limit geoms 3}
    \|DA_{h,f_{\ua}(p)}^{-1}\|_{{\textit{op}}} \leq & ~ \sup_{p\in G(\ua)} \{\|Dh^{-1}(p)\|_{\it op}\}<C'_h,
    \end{align}    
    where the first inequality is obtained by \eqref{eq control of diameter of nth image by derivatives} and $r_{\ua} := \mathrm{diam}(G(\ua))$.
       Similar to the previous estimates in proof of Theorem \ref{thm: Bounded distortion limit geom improved} for $\kappa>1$ with $\kappa\cdot \mu^{\alpha}<1$ we have 
    $$ \|Df_{\ua}^{-1}\|_{{\textit{op}}} \cdot  \|Df_{\ua}\|_{{\textit{op}}} \leq \kappa^n ,\;\; \|Df_{\ua}\|_{{\textit{op}}} \leq C_1\cdot \mu^{n\alpha}.$$
    Therefore, according to Lemma \ref{lemma estimation of chosing linear instead of non linear part} for 
    $$H_{\ua}:= A_{f_{\ua},p}^{-1} \circ A_{h,f_{\ua}(p)}^{-1} \circ h \circ A_{h,f_{\ua}(c_0)}- \id$$
    one has that for some constant $C_{h,k}$ depended continuously on $h,K$,
    \begin{equation}\label{eq closenes of linear estimate in attractorness proof}
        \|H_{\ua}\|_{\cC^{1+\alpha}}\leq C_{h,K}\cdot (\kappa\cdot \mu^{\alpha})^n.
    \end{equation}
    Hence, since the set $\{k_m^{\underline{\theta a}}:m\geq 1,\; \ua\in \Sigma^{\rm fin} \}$ is bounded in $\cC^{1+\alpha}$ topology by \eqref{eq closenes of linear estimate in attractorness proof}  and the H\"older estimate  
    \eqref{Hormander} we obtain \eqref{eq lemma attractor statement for analysis} in the following way
    \begin{align}
\|F_{\ua}\|_{\cC^{1+\alpha}}
    &\leq 
    C_{1+\alpha}\cdot \big(\|H_{\ua}\|_{\cC^{1+\alpha}}\cdot \|k_n^{\underline{\theta a}}\|_{\cC^{1}}^{\alpha} + \|H_{\ua}\|_{\cC^{1}}\cdot \|k_n^{\underline{\theta a}}\|_{\cC^{1+\alpha}}+\|H_{\ua}\|_{\cC^{0}}\big)
    \notag\\
   & \leq C_{h,K} \cdot (\kappa\cdot \mu^{\alpha})^n.
    \end{align}
\end{proof}

\subsection{Relative configurations and renormalization operators}
Given a pair of bunched Cantor sets $(K,K')$
the group $\Aff(d,\R)$ 
naturally acts on the space 
$\cP\times \cP'$ by the map $(h,h')\mapsto (A\circ h,A\circ h')$ for $A\in \Aff(d,\R)$. Let $\cQ$ be the quotient of $\cP\times \cP'$ with the quotient map $Q:\cP\times \cP'\to \cQ$ under this action equipped with the quotient topology.
We call $\cQ$ as the space of \emph{relative configurations} of the pair of Cantor sets. We denote the representative of $(h,h')\in \cP\times \cP'$ in $\cQ$ as $[h,h']$. In the holomorphic case we define $\cQ^{\rm Hol}$ for the quotient of the action of $\Aff(d,\C)$ on the space $\cP^{\rm Hol}\times \cP'^{\rm Hol}$.
When we are analyzing the topology of $\cQ$ (or $\cQ^{\rm Hol}$)
we write
\begin{align}\label{eq: normalization in Q}
    (\widehat{h}_p,A_{h,p}^{-1}\circ A_{h',p'}\circ \widehat{h'}_{p'})=(A_{h,p}^{-1}\circ h,A_{h,p}^{-1}\circ h')
\end{align}
as a representative of $[h,h']\in \cQ$
(or $\cQ^{\rm Hol}$) 
for some points $p\in K$ and $p'\in K'$.

Given a pair $(K,K')$, we define \emph{renormalization operators} on the space $\cP\times \cP'$.
For $(\ua,\ua')\in \Sigma^{\rm{fin}}\times \Sigma'^{\rm{fin}},(h,h')\in \cP\times \cP'$, we have the renormalization operator
\begin{align}
T_{\ua, \ua'}:  \cP\times \cP'&\to \cP\times \cP',\notag\\
(h,h')&\mapsto \left(T_{\ua}(h), T'_{\ua'}(h')\right).
\end{align}
We can also allow one of the words $\ua$ or $\ua'$ be void and define 
$$T_{\ua,\emptyset}(h,h'):=\left(T_{\ua}(a),h'\right),\quad T_{\emptyset,\ua'}(h,h'):=\left(h,T_{\ua'}(h')\right).$$
 Let $\Sigmafin_{*}: = \Sigmafin\cup \emptyset$. The family of renormalization operators 
 \begin{align}\label{eq: definition of family of renormalization operators}
 \cR:=\left\{T_{\ua,\ua'}:(\ua,\ua')\in \Sigmafin_{*} \times \Sigma'^{\rm fin}_{*} \backslash \{(\emptyset,\emptyset)\}\right\}
 \end{align}
 is acting form right on $\cP\times \cP'$ and invertible affine maps $\Aff(d,\R)$ acting from left on the space of relative configurations. So they commute with each other and thus renormalization operators naturally act on the space $\cQ$ (or $\cQ^{\rm Hol}$) with the map
 $$T_{\ua,\ua'}:[h,h']\mapsto [T_{\ua}(h),T_{\ua'}(h')].$$
The family of operators $\cT=\{T_{\ua}:\ua\in \Sigmafin\}$ is generated by the finite family $\cT_1$. Thus, $\cR$ is generated by the finite family 
$\cR_1:=\{\id\}\times \cT'_1 ~\cup~ \cT_1\times \{\id\}$
  consisting of operators 
  $T_{\ua,\ua'}$,
  where the sum of the lengths of $\ua,\ua'$ equals one.
  
We define the space of \emph{relative affine configurations of limit geometries} as the quotient of the space $\Phi(\cA)\times \Phi'(\cA')\subset \cP\times \cP'$ under the left action of the group $\Aff(d,\R)$ containing relative configurations with representatives $[A\circ k^{\underline{\theta}},A'\circ k'^{\underline{\theta}'}]\in \cQ$. 
Like the space $\cP\times \cP'$, the group $\Aff(d,\R)$ acts on the space $\cA\times \cA'$ by $$(A,A',\underline{\theta},\underline{\theta}')\mapsto (BA,BA',\underline{\theta},\underline{\theta}').$$

Given a pair of bunched Cantor sets $(K,K')$, we define the space of 
\emph{representation of relative affine configuration of limit geometries} as the quotient of the space $\cA\times \cA'$ under the action of the group $\Aff(d,\R)$ and denote it by $\cQ_{\Aff}$ with the quotient map 
\begin{align}\label{eq: quotient map definition of Q_aff}
Q_{\cA}:\cA\times \cA'&\to \cQ_{\Aff},\notag\\
(A,A',\underline{\theta},\underline{\theta}')&\to (A^{-1}\circ A,\underline{\theta},\underline{\theta}').
\end{align}
Therefore, we have that  $\cQ_{\Aff}\cong\Aff(d,\R)\times \Sigma^-\times \Sigma'^-$. The map
\begin{align}
    \Phi_q: \cQ_{\Aff}&\to \cQ,\notag\\
    (B,\underline{\theta},\underline{\theta}')&\mapsto [k^{\underline{\theta}},B\circ k^{\underline{\theta}'}]
\end{align}
is the quotient of the map 
$\Phi\times \Phi' : \cA\times\cA' \to \cP\times \cP'$ 
such that the following diagram commutes.
\begin{equation}\label{diagram quotient}
     \begin{tikzcd}
\cA\times \cA' \arrow{r}{\Phi\times \Phi'} \arrow[swap]{d}{Q_{\cA}} & \cP\times \cP' \arrow{d}{Q} \\%
\cQ_{\Aff}\arrow{r}{\Phi_q}& \cQ
\end{tikzcd}
\end{equation}
Similar to Lemma \ref{lem: affineness of action of renormaliztions on limit geoms}, the action of family $\cR$ on the space $\cP\times \cP'$ has a pull back via the map $\Phi_{q}$. There is a family $\cR^*$ of renormalization operators acting on the space $\cA\times \cA'$ . More precisely, for $(\ua,\ua')\in \Sigmafin\times \Sigma'^{\rm fin}$ there is an operator
$T^{*}_{\ua,\ua'}:\cA\times \cA'\to \cA\times \cA'$ defined by
$$T^{*}_{\ua,\ua'}(A,A',\underline{\theta},\underline{\theta}'):=(A\circ F^{\underline{\theta a}},A'\circ F^{\underline{\theta}' \ua'},\underline{\theta a},\underline{\theta}'\ua').$$
Similar to above, we can allow one of $\ua$ or $\ua'$ be void and define
$$T^{*}_{\emptyset,\ua'}(A,A',\underline{\theta},\underline{\theta}'):=(A,A'\circ F^{\underline{\theta}' \ua'},\underline{\theta},\underline{\theta}'\ua'),\;
T^{*}_{\ua,\emptyset}(A,A',\underline{\theta},\underline{\theta}'):=(A\circ F^{\underline{\theta a}},A',\underline{\theta},\underline{\theta}').$$
Therefore, 
 $\cR^*:=\left\{T^*_{\ua,\ua'}:(\ua,\ua')\in \Sigmafin_{*} \times \Sigma'^{\rm fin}_{*} \backslash \{(\emptyset,\emptyset)\}\right\}$. Similar to above, we have that $\cR^*$ is generated with the finite family $\cR^*_1$ of operators consisting of maps $T^*_{\ua,\ua'}$ with sum of lengths of $\ua$ and $\ua'$ equal to one.
 
 \noindent Since $\cQ_{\Aff}$ is the quotient of $\cA\times \cA'$, the family $\cR^*$ naturally acts on $\cQ_{\Aff}$.
 In particular, the action of $T^*_{\ua,\ua'}\in \cR^*$ on $\cQ_{\Aff}$ is defined by
 \begin{align}\label{eq: main action of renormalization operators}
     T^*_{\ua,\ua'}(B,\underline{\theta},\underline{\theta}'):=
 \left((F^{\underline{\theta a}})^{-1}\circ B\circ  F^{\underline{\theta}' \ua'},\underline{\theta}\ua,\underline{\theta}'\ua'\right),
 \end{align}
where whenever $\ub$ is void we set $F^{\underline{\theta b}}:= \id$. Similar to Lemma \ref{lem: affineness of action of renormaliztions on limit geoms}, operators $T^*_{\ua,\ua'}$ are H\"older continuous maps over $\cQ_{\Aff}$ such that the following diagram commutes.
 \begin{equation}\label{diagram renormalization of pair of cantors}
     \begin{tikzcd}
\cQ_{\Aff} \arrow{r}{\Phi_q} \arrow[swap]{d}{T^{*}_{\ua,\ua'}} & \cQ \arrow{d}{T_{\ua,\ua'}} \\%
\cQ_{\Aff} \arrow{r}{\Phi_q}& \cQ
\end{tikzcd}
\end{equation}
 The family $\cR_1^*$ of renormalization operators with total length of 1 can be partitioned in to two disjoint sets of operators. In particular, since
 $F^{\underline{\theta'} \ua'}$
 and $(F^{\underline{\theta a}})^{-1}$ are contracting and 
 expanding affine maps, respectively, $\cR_1^*=\cD_1^*\sqcup \cE_1^*$ such that $\cD_1^*, \cE_1^*$ are families of contracting and expanding operators of $\cR^*$, respectively where
 \begin{equation}\label{eq: contracting and expanding renormalization operators definition}
     \cD_1^*:=\{T^*_{\ua,\ua'}\in \cR_1^*: \ua\; \text{is void}\},\;\;
 \cE_1^*:=\{T^*_{\ua,\ua'}\in \cR_1^*: \ua'\; \text{is void}\}.
 \end{equation}
 Similarly, there is the partition $\cR_1 = \cD_1\sqcup \cE_1$. In the holomorphic case $\cQ_{\Aff}^{\rm {Hol}}$ and renormalization operators acting on this space are defined in a similar manner.

\subsection{Affine Cantor sets} \label{sec: Affine Cantor sets}
Recall that an affine Cantor set $K$ is generated by an expanding generator $g$ which is an affine map on a neighborhood of $G(a)$ for any $a\in \mathfrak{A}$. When the symbolic type of $K$ is the full shift it is more convenient to define $K$ by the family of contracting maps  $\{f_1,f_2,\dots,f_{m}\}\subset \Aff(d,\R)$
     where $m:=\#\mathfrak{A}$ and $f_j : = (g|_{G(a)})^{-1}$ for $j=1,\dots,m$.

Limit geometries of affine Cantor sets are trivial identity maps, since for all $n\in \N$ we have $k_n^{\underline{\theta}}= A_{\underline{\theta}_{n}}^{-1}\circ f_{\underline{\theta}_{n}} = \id$.
Hence, for all $\underline{\theta}\in \Sigma^-$ we have that $k^{\underline{\theta}}=\lim_{n\to \infty} k_n^{\underline{\theta}}=\id$.
Consequently, the family of operators $\cT^*$  acting on the space $\cA$ of representations of affine configuration of limit geometries consists of the composition of the generators of $K$. More precisely, for any 
\((A ,\underline{\theta})\in \cA\) (or $\cA^{\rm Hol}$) and $\ua=(a_0,\cdots,a_n)\in \Sigma^{\rm {fin}}$ with $a_0=\theta_0$ we have
$
F^{\underline{\theta a}}= k^{\underline{\theta}}\circ f_{\ua}\circ (k^{\underline{\theta}\ua})^{-1}=f_{\ua}
$. So $F^{\underline{\theta a}}$ is independent of $\underline{\theta}$.
Thus, we can write
the action of $ T^*_{\ua}$ over \(\cA\) (or $\cA^{\rm Hol}$) 
as
\begin{equation}\label{eq: renormalization operators on affines}
  T^*_{\ua}(A,\underline{\theta}) = (A \circ f_{\ua},\underline{\theta a}).     
\end{equation}
Furthermore, the action of renormalization operators $T^*_{\ua,\ua'}\in \cR^*$ on the space $\cQ_{\Aff}$ defined by \eqref{eq: main action of renormalization operators} could be simplified as \begin{align}\label{eq: Affines, main action of renormalization operators}
     T^*_{\ua,\ua'}(B,\underline{\theta},\underline{\theta}'):=
 \left((f_{\ua})^{-1}\circ B\circ  f_{\ua'},\underline{\theta a},\underline{\theta}' \ua'\right).
 \end{align}
 Another consequence of the relation $F^{\underline{\theta a}}=f_{\ua}$ is that in the case of affine bunched Cantor sets the family of renormalization operators $\cR^*$ has an action on the group $\Aff(d,\R)$. Indeed, any $T^*_{\ua,\ua'}\in \cR^*$ defines a map on $\Aff(d,\R)$ defined by
 \begin{align}\label{eq: full shift operator act on affine space}
 T^*_{\ua,\ua'}(B):=(f_{\ua})^{-1}\circ B\circ f_{\ua'}.
 \end{align}
 The partitioning
 $\cR_1^* = \cD_1^*\sqcup \cE_1^*$ to the contracting and expanding operators in the affine case is special since $f_{\ua},f_{\ua'}$ are affine generators of $K,K'$ respectively when lengths of $\ua,\ua'$ are equal to 1. Thus, $\cD_1^*,\cE_1^* \subset \Aff(d,\R)$ and
 \begin{equation}\label{eq: partitioning to conracting and expanding operators Affine}
     \cD_1^*:=\{f_{(a'_0,a'_1)}: (a'_0,a'_1)\in \mathfrak{B}'\},\;\;
 \cE_1^*:=\{f_{(a_0,a_1)}^{-1}: (a_0,a_1)\in \mathfrak{B}\}.
 \end{equation}
 The action of the sets $\cD_1^*\subset \Aff(d,\R)$ on the group $\Aff(d,\R)$ is the right action while the action of the set $\cE_1^*\subset \Aff(d,\R)$ is the left action. Indeed, action of $\cD_1^*$ and $\cE_1^*$ on $B\in \Aff(d,\R)$ can be described as maps
 \begin{equation}\label{eq: left and right action of contracting and expanding operators}
     f_{(a_0',a_1')}: B\mapsto B\circ f_{(a_0',a_1')},\;\; f_{(a_0,a_1)}^{-1}: B\mapsto f_{(a_0,a_1)}^{-1}\circ B.
 \end{equation}
 We have explored these actions in more details in Appendix \ref{appendix affine actions}. Indeed, 
 by 
 Lemma \ref{lem: affine left and right action details} 
we can compute the action of $\cE_1^*$ and $\cD_1^*$ on $\Aff(d,\R)$. 
If $\psi:=[x\mapsto Px+w] ,\phi:=[x\mapsto P'x+w']\in\Aff(d,\R)$ are two affine generators of $K,K'$,
respectively,
then $\psi^{-1}\in \cE^*_1$ and $\phi\in\cD^*_1$ are
 operators of the pair $(K,K')$ acting on $\Aff(d,\R)$ by
\begin{align}
  \psi^{-1}:\;   (v,t,A)&\mapsto (P^{-1}(v-w),t/ s_{P}, \hat P^{-1}  A), \label{eq: expanding action computation}\\ 
 \phi:\;
 (v,t,A)&\mapsto (v+t\cdot Aw',t \cdot s_{P'}, A\hat P'),\label{eq: contracting action computation}
\end{align} 
where $s_A :=e^{i\operatorname{arg}(\operatorname{det}(A))}  \sqrt[d]{|\operatorname{det}(A)|}$, $\hat A : = s_{A}^{-1}\cdot A$, for any $A\in \GL(d,\mathbb{F})$ and $(v,t,A)\in \R^d \rtimes  (\R^* \times\SL(d,\R)) \cong \Aff(d,\R)$. 
\begin{corollary}\label{corol: action of renormalization is product and semi direct product}
    Given a pair of affine Cantor sets $(K,K')$, expanding renormalization operators $\cE_1^*$ are the generators of $K$ acting on $\Aff(d,\bF)$ as the product of their action on $\bF^*\times \bF^d$ and $\SL^\star(d,\bF)$ defined by \eqref{eq: expanding action computation}, while the contracting renormalization operators $\cD_1^*$ are the generators of $K'$ acting on $\Aff(d,\bF)$ as the semidirect product of their action on $\bF^d$ and $\bF^*\times\SL^\star(d,\bF)$ defined by \eqref{eq: contracting action computation}. 
\end{corollary}

\section{The covering criterion}\label{sec: covering criterion}

In this section we prove a key result that shows that how the the action of renormalization operators of  a pair of bunched Cantor sets $(K,K')$ is relevant to the problem of stable intersection of $K$ and $K'$. 
\subsection{Intersecting configurations}\label{sec: intersecting configurations} 
Given a pair of configurations $\left(h_{a}, h'_{a'}\right) \in \mathcal{P}(a) \times \mathcal{P}'(a')$  (or $\cP^{\rm Hol}(a)\times \cP'^{\rm Hol}(a')$) we say that it is

\begin{itemize}
  \item \emph{linked} whenever $h_{a}(\overline{G(a)}) \cap h'_{a'}(\overline{G'\left(a'\right)}) \neq \emptyset$;

  \item \emph{intersecting} whenever $h_{a}(K(a)) \cap h'_{a'}\left(K'\left(a'\right)\right) \neq \emptyset$;

  \item  \emph{stably intersecting} whenever $\tilde{h}_{a}(\tilde{K}(a)) \cap \tilde{h}'_{a'}(\tilde{K}'(a')) \neq \emptyset$ for any pairs of Cantor sets $(\tilde{K}, \tilde{K}') \in \Omega_{\Sigma,d}^{1+\alpha} \times \Omega_{\Sigma',d}^{1+\alpha}$ in a small neighborhood of $(K, K')$ and any configuration pair $(\tilde{h}_{a}, \tilde{h}'_{a'})$ that is sufficiently close to $\left(h_{a}, h'_{a'}\right)$ in the $\cC^{1+\alpha}$ topology at $G(a) \cap \tilde{G}'(a)$ and $G\left(a'\right) \cap \tilde{G}'\left(a'\right)$ for some $\varepsilon>0$.

\end{itemize}
The action of $\Aff(d,\R)$ on 
$\cP_a\times \cP'_{a'}$ (or of $\Aff(d,\C)$ on $\cP_a^{\rm Hol}\times \cP'^{\rm{Hol}}_{a'}$) preserves above properties.
Therefore, if any of these properties does holds, we say that the relative configuration $[h_a,h'_{a'}]\in \cQ$ is linked, intersecting or stably intersecting (respectively).

\begin{lemma}\label{lemma intersecting vs linked}
    A pair of configurations $(h,h')\in \cP\times \cP'$ (or $\cP^{\rm Hol}\times \cP'^{\rm Hol}$) is intersecting if and only if there exists
    $(\ua,\ua')\in \Sigma  \times \Sigma'$ such that the pairs $T_{\ua_n,\ua'_n}(h,h')$ are linked for all $n\in \N$.
\end{lemma}

\begin{proof}   Assume that  $p\in h(K(a_0))\cap h'(K'(a'_0))$ for some $a_0\in \mathfrak{A}$ and $a'_0\in \mathfrak{A}'$. Then, there exists letters $a_1\in \mathfrak{A}$ and $a'_1\in \mathfrak{A}'$ such that $p\in T_{(a_0,a_1)}(h)(G(a_1)) \cap T'_{(a'_0,a'_1)}(h)(G'(a'_1))$. By repeating this argument infinitely many times we may build 
$\ua:=(a_0,a_1,\dots)$ and $\ua':= (a'_0,a'_1,\dots)$ such that $p\in T_{\ua_n}(h)(G(a_n))\cap T_{\ua'_n}(h')(G'(a'_n))$.

    For the other direction, let $h_n : = T_{\ua_n}(h)$, $h'_n : = T'_{\ua'_n}(h')$, $p_{n} \in h_{n}(K(a_n)) \subset h(K)$ and $p_n'\in h_n'(K'(a'_n)) \subset h'(K')$ be arbitrary points. We know that $\diam(h_{n}(K))$ and $\diam(h'_{n}(K'))$ converge to 0 as $n \rightarrow \infty$ since they are controlled by $\diam(G(\ua_{n}))$ and $\diam(G'(\ua'_{n}))$ respectively. Therefore, knowing that $h_n(G(a_n))\cap h'_n(G'(a'_n))\neq \emptyset$, the diameter of the compact set $H_n : = h_n(G(a_n))\cup h_n'(G'(a'_n))$ converges exponentially to zero as $n\to \infty$. Thus, $\{H_k\}_{k\in \N}$ is a chain of non-empty compact sets with the property $H_{k+1}\subset H_{k}$ for all $k\in \N$ which implies that  $\bigcap_{k=1}^{\infty} H_k$ is a non-empty singleton $\{p\}$. Thus, $p = \lim _{n \rightarrow \infty} p_{n}=\lim _{n \rightarrow \infty} p_{n}'\in h(K(a_0))\cap h(K'(a'_0))$.
\end{proof}

The following lemma unveils  the main idea of the stable intersection criterion. 
 Later we will use stability of strong covering in this lemma to deduce the stable intersection. 
\begin{lemma}\label{lemma intersection criteria}
 A relative configuration $\left[h, h'\right]\in \cQ$ (or $\cQ^{\rm Hol}$) is intersecting if it belongs to a bounded set $V\subset \cQ$ (or $\cQ^{\rm Hol}$) satisfying the covering condition with respect to the action of the family $\cR$ (defined in \eqref{eq: definition of family of renormalization operators}) on the space $\cQ$ (or $\cQ^{Hol}$).
 \end{lemma}

\begin{proof} 
Since $V$ satisfies the covering condition, there is a sequence of operators $\{\Psi_n\}_{n\in \N}$ in $\cR$ such that $[h_n,h_n']:= \Psi_n([h_{n-1},h'_{n-1}])\in V$ for all $n\in \N$ where $[h_0,h_0']:=[h,h']$. $V$ being bounded implies that the sequence $[h_n,h_n']$ is bounded. To obtain that $[h,h']$ is intersecting it is enough to show that relative configurations $\left[h_{n}, h_{n}'\right]$ are all linked (See Lemma \ref{lemma intersecting vs linked}).
Assume the contrary that $(h_k,h_k')\in \cP(a_k)\times \cP'(a'_k)$ is not linked for some $k$, where  $a_k\in \mathfrak{A}$ and $a'_k\in \mathfrak{A}'$ are corresponding letters determined by $\Psi_k$. Then, the sets $h_k(G(a_k))$ and $h'_k(G'(a'_k))$ has some distance $\gamma>0$.
Thus, $\dist(h_{n}(G(a_n)),h'_{n}(G'(a'_n)))\geq \gamma$ for all $n\geq k$,  since $h_n(G(a_n))\subset h_k(G(a_k))$ and $h'_n(G'(a'_n))\subset h'_k(G'(a'_k))$.
  Let $p$ and $p'$ be some base points in $K,K'$. Then by \eqref{eq control of diameter of nth image by derivatives} for some constant $C$ depended only on $h$ and $K$,
 $\dist(A_{h_n,p}^{-1}\circ h_n, A_{h_n,p}^{-1}\circ h'_n)\geq C\cdot \gamma \cdot \mu^{-n} $
which contradicts with $[h_n,h'_n]$ being bounded in $\cQ$.
\end{proof}

As a direct consequence of Lemma \ref{lemma intersection criteria} we have the following corollary.

\begin{corollary}\label{corol:  intersecting pairs in affine configs inside a covering region}
 An affine configurations $[A\circ k^{\underline{\theta}}, A'\circ k^{\underline{\theta}'}]\in \cQ$ (or $\cQ^{\rm Hol}$) is intersecting if its representative
   $( A^{-1} \circ A',\underline{\theta}, \underline{\theta}')\in \cQ_{\Aff}$ (or $\cQ_{\Aff}^{\rm Hol}$) belongs to the closure of some
  open relatively compact set $\cW\subset \cQ_{\Aff}$ (or $\cQ_{\Aff}^{\rm{Hol}}$) that satisfies covering condition \eqref{eq covering condition} with respect to the family $\cR^*$. 
\end{corollary}
\begin{proof} 
$\Phi_q(( A^{-1} \circ A',\underline{\theta}, \underline{\theta}'))\in \cQ$ is a relative configuration which belongs to the bounded set $\Phi_q(\cW)$. In addition, the commutative  diagram \eqref{diagram renormalization of pair of cantors} implies that $\Phi_q(\cW)$ satisfies covering condition \eqref{eq covering condition}
with respect to the family $\cR$. So Lemma \ref{lemma intersection criteria}
is applicable in here and we're done. The holomorphic case follows similarly.
\end{proof}

Here, we recall from \S\ref{sec: Affine Cantor sets} that when $K$ and $K'$ are affine bunched (or bunched holomorphic) Cantor sets then the family $\cR^*$ acts on the group $\Aff(d,\R)$ (or $\Aff(d,\C)$) via the map defined in \eqref{eq: full shift operator act on affine space}.
\begin{theorem}\label{thm: covering criterion in affine cantor sets with full shift}
Let $(K,K')$ be a pair of affine bunched (or bunched holomorphic) Cantor sets such that both of their symbolic types 
$\Sigma,\Sigma'$
are full shift. Let $\cL\subset \Aff(d,\R)$ (or $\Aff(d,\C)$) be an open relatively compact set satisfying the covering condition \eqref{eq covering condition} with respect to the action of $\cR_1^*$ on $\Aff(d,\R)$ (or $\Aff(d,\C)$). Then, for any $A\in \overline{\cL}$, $\underline{\theta}\in \Sigma^-$ and $\underline{\theta}'\in \Sigma'^-$, $[k^{\underline{\theta}},A\circ k^{\underline{\theta}'}]\in \cQ$ (or $\cQ^{\rm Hol}$) is intersecting.
\end{theorem}
\begin{proof}
   Let $\cW:= \cL\times \Sigma^-\times \Sigma'^- \subset \cQ_{\Aff}$. Since the symbolic types of $K$ and $K'$ are full shift, \eqref{eq: Affines, main action of renormalization operators} implies that $\cW$ satisfies the covering condition \eqref{eq covering condition} with respect to the action of $\cR_1^*$ on the space $\cQ_{\Aff}$. We know that $\cR_1^*$ generates $\cR^*$. So by Corollary \ref{corol:  intersecting pairs in affine configs inside a covering region} the relative configuration $[k^{\underline{\theta}},A\circ k^{\underline{\theta}'}]\in \cQ$ with representative $(A,{\underline{\theta}}, {\underline{\theta}'})\in \overline{\cW}\subset \cQ_{\Aff}$ is intersecting since $\cW$ is an open relatively compact subset of $\cQ_{\Aff}$.  
\end{proof}

\subsection{Covering criterion for stable intersection}

The following theorem is one of the main results of this paper. Indeed, Theorems \ref{thm:A} and \ref{thm:A-complex} are its immediate consequences.

\begin{theorem}[Covering criterion for stable intersection] \label{thm: covering criterion for intersection}
Let $(K,K')$ be a pair of bunched (or bunched holomorphic) Cantor sets in $\R^d$ (or $\C^d$).  Assume that an open relatively compact set $\cW\subset \cQ_{\Aff}$ (or $\cW\subset \cQ_{\Aff}^{\rm Hol}$) satisfies strong covering condition \eqref{eq: covering condition delta interior} with respect to the finite family $\cR_1^*$ of renormalization operators of the pair $(K,K')$. Then,
\begin{enumerate}
    \item  for every $(\tilde K,\tilde K')$ in an open neighborhood of $(K,K')$, $\cW$ also satisfies strong covering condition with respect to the family $\tilde{\cR}_1^*$ of renormalization operators of the pair $(\tilde K,\tilde K')$;
    \item each affine relative configuration contained in $\Phi_q(\overline{\cW})$ is stably intersecting.
\end{enumerate}
\end{theorem}
\begin{proof}
We only present the proof for the real case. The holomorphic case follows from a same argument.

{\it Proof of} (1).
The family $\cR^*_1$ acts on the locally compact space $\cQ_{\Aff}\cong \Aff(d,\R)\times \Sigma^-\times \Sigma'^-$. Moreover, it consists of continuous maps varying continuously with respect to the pairs of Cantor sets $(K,K')$ due to Lemmas \ref{lem: affineness of action of renormaliztions on limit geoms}, \ref{lem: continuity of renormalization operators}. So we can apply Lemma \ref{lemma stability of strong covering} which concludes this item.

{\it Proof of} (2). According to Corollary \ref{corol:  intersecting pairs in affine configs inside a covering region} 
any relative configuration $[ k^{\underline{\theta}}, B\circ k^{\underline{\theta}'}]\in \cQ$ with 
 $(B,\underline{\theta}, \underline{\theta}')\in \overline{\cW}$ is intersecting. Moreover, by item (1), for $(\tilde K,\tilde K')$ sufficiently close to $(K,K')$ the relative configuration  
 $[\tilde k^{\underline{\theta}}, B\circ \tilde k^{\underline{\theta}'}]\in \tilde{\cQ}$ is also intersecting.
 So, in order to prove that $[ k^{\underline{\theta}},B\circ k^{\underline{\theta}'}]\in \cQ$ is stably intersecting it suffices to show that all relative configurations $\left[h \circ k^{\underline{\theta}},h' \circ B \circ k^{\underline{\theta}'}\right]\in \cQ$ in a neighborhood of $[ k^{\underline{\theta}},B\circ k^{\underline{\theta}'}]= \Phi_q\left((B,\underline{\theta},\underline{\theta}')\right)$ in $\cQ$ are
 intersecting.
 
To do so, we show that for small enough $\delta'>0$ the neighborhood $B_{\delta'}\left((\Phi_q(\overline{\cW})\right)\subset \cQ$ of $\Phi_q(\overline{\cW})$, consisting of relative configurations $\left[h \circ k^{\underline{\theta}},h' \circ B \circ k^{\underline{\theta}'}\right]\in \cQ$ with $h,h'$ being $\delta'$-close to $\id$, satisfies the covering condition with respect to the family
$\cR$. This implies that any relative configuration in this neighborhood is intersecting (see  Lemma \ref{lemma intersection criteria}).
More precisely, we will prove that there is $\lambda \in (0,1)$ and an integer $l\in \N$ such that for any $[\hat h,\hat h']\in B_{\delta'}\left(\Phi_q(\overline{\cW})\right)$ there exist 
$\ua\in \Sigma^{\rm{fin}},\ua'\in\Sigma'^{\rm{fin}}$ with lengths less than or equal to $l$ such that $T_{\ua,\ua'}([\hat h,\hat h'])\in B_{\lambda \cdot \delta'}\left(\Phi_q(\overline{\cW})\right)$. This gives the (strong) covering condition for $B_{\delta'}(\Phi_q(\overline{\cW}))$ with respect to the finite family $\cR_{l}\subset \cR$ consisting of operators $T_{\ua,\ua'}\in \cR$ with lengths of $\ua,\ua'$ both less than or equal to $l$.

Within the proof, to 
estimate the distance of a relative configuration  $[\phi,\phi']\in \cQ$ from $\Phi_q(\overline{\cW})$ in the quotient topology we first write  
\begin{align}\label{eq: normalized form for relative configs in proof of main thm}
\left[\phi, \phi'\right]=\left[\eta_{\phi} \circ k^{\underline{\theta}(\phi)}, \eta'_{\phi'} \circ B_{\phi,\phi'} \circ k^{\underline{\theta}'(\phi')}\right]    
\end{align}
such that the affine estimates of $\eta_{\phi}$ and $\eta'_{\phi'}$ at points $c_{{\theta(\phi)}_0}$ and $B_{\phi,\phi'}(c'_{{\theta'(\phi')}_0})$ (respectively) are $\id\in \Aff(d,\R)$ and $\left(B_{\phi,\phi'},{\underline{\theta}(\phi)}, {\underline{\theta}'(\phi')}\right)\in \cQ_{\Aff}$ for some $\underline{\theta}(\phi)$ and $\underline{\theta}'(\phi')$ in $\Sigma^-$ and $\Sigma'^-$, respectively. Then, we will analyze the distance of $\eta_{\phi},\eta'_{\phi'}$ from $\id$, which gives the required estimate if $\left(B_{\phi,\phi'},{\underline{\theta}(\phi)}, {\underline{\theta}'(\phi')}\right)\in \overline{\cW}$.

Let $[\hat h,\hat h']\in B_{\delta'}(\Phi_q(\overline{\cW}))$ 
 be a relative configuration near some $\left[ k^{\underline{\theta}}, B \circ k^{\underline{\theta}'}\right]\in \Phi_q(\overline{\cW})$ where $\delta'>0$ is a constant which will be defined in \eqref{eq: definition of delta' proof of the main theorem}.
 Thus, there are $h,h'$ $\delta'$-close to $\id$ such that
 \begin{align}
    [\hat h,\hat h']=& \left[h \circ k^{\underline{\theta}},h' \circ B \circ k^{\underline{\theta}'}\right]\notag\\
=& \left[ A_{h,p}^{-1}\circ h\circ k^{\underline{\theta}}, A_{h,p}^{-1}\circ h' \circ B\circ k^{\underline{\theta}'}\right]\notag
 \end{align}
where $A_{h,p}$ is the affine estimate of $h$ at $p:=c_{\theta_0}$.
Now, we choose
\begin{align*}
    \eta_{\hat h}&:= A_{h,p}^{-1} \circ h,\\
    \eta'_{\hat h'} &:=  A_{h,p}^{-1}\circ h'\circ A'^{-1}_{h',p'}\circ A_{h,p},\\
    B_{\hat h, \hat h'} &:=  A_{h,p}^{-1}\circ A'_{h',p'} \circ  B,
\end{align*}
with $A_{h',p'}$ be the affine estimate of $h'$ at $p':=B(c'_{\theta'_0})$. Notice that affine estimates of $\eta_{\hat h},\eta'_{\hat h'}$ at the points $c_{\theta_0},B_{\hat h,\hat h'}(c'_{\theta'_0})$ are $\id$. Hence, we have 
$$[\hat h,\hat h']= [\eta_{\hat h} \circ k^{\underline{\theta}(\hat h)},
\eta'_{\hat h'} \circ B_{\hat h, \hat h'}\circ k^{\underline{\theta'}(\hat h')}],$$
where $\underline{\theta}(\hat h) := \underline{\theta},\underline{\theta'}(\hat h'):=\underline{\theta'}$.
Continuity of $A_{h,p}^{-1}\circ A'_{h',p'}$ with respect to $h,h'$ implies that
for any $\delta>0$ (which will be determined in few lines later) there exists $\delta_1: = \delta_1(\cW,\delta)>0$ such that if $h,h'$ are $\delta_1$-close to $\id$ then 
$(B_{\hat h, \hat h'},\underline{\theta},\underline{\theta}')\in \cQ_{\Aff}$
is $\delta$-close to $(B,\underline{\theta},\underline{\theta}')\in \overline{\cW}$. Thus, $(B_{\hat h,\hat h'},\underline{\theta},\underline{\theta}') \in B_{\delta}(\cW)$
. On the other hand, since 
$\cW\subset \cQ_{\Aff}$ satisfies strong covering \eqref{eq: covering condition delta interior} with respect to $\cR_1^*$, 
there exists $\delta>0$ such that 
$$
B_{\delta}(\cW)\subset \bigcup_{\Psi \in \cR^*_1} \Psi^{-1}(\cW_{(\delta)}).
$$
Hence, 
for any $l_0\in \N$ there are $\ua\in \Sigmafin$, $\ua'\in
\Sigma'^{\rm fin}$
with lengths at least $l_0$ such that 
\begin{equation}\label{eq delta interior of renormalization process}
   T_{\ua, \ua'}\left(\left[k^{\underline{\theta}}, B_{\hat h,\hat h'} \circ k^{\underline{\theta}'}\right]\right) \in \Phi_q(\cW_{(\delta)}).
\end{equation}
Now we define
$
[\hat h_1, \hat h'_{1}]:=T_{\ua, \ua'}\left([\hat h, \hat h']\right).
$
According to \eqref{eq affine action of renormalization operator with word an} we have
$$
\begin{aligned}
    {[\hat h_1, \hat h'_{1}]} & = [\eta_{\hat h} \circ k^{\underline{\theta}(\hat h)}\circ f_{\ua},\eta'_{\hat h'} \circ B_{\hat h,\hat h'} \circ k^{\underline{\theta'}(\hat h')}\circ 
f_{\ua'}]\\
& = [\eta_{\hat h} \circ F^{\underline{\theta}(\hat h)\ua}\circ
k^{\underline{\theta}(\hat h) \ua},
\eta'_{\hat h'} \circ B_{\hat h,\hat h'}\circ F^{\underline{\theta'}(\hat{h}') \ua'}\circ k^{\underline{\theta'}(\hat{h}') \ua'}].
\end{aligned}
$$
Next, we write $[\hat h_1,\hat{h}'_1]$ in the form of \eqref{eq: normalized form for relative configs in proof of main thm}:
$$
 [\hat h_1,\hat{h}'_1]= [\eta_{\hat h_1}\circ k^{\underline{\theta}(\hat h_1)},\eta'_{\hat h'_1}\circ B_{\hat h_1,\hat h'_1}\circ  k^{\underline{\theta'}(\hat h'_1)}],
$$
where 
$\underline{\theta}(\hat h_1):=\underline{\theta}(\hat h)\ua$, $\underline{\theta'}(\hat h'_1):=\underline{\theta}'(\hat h')\ua'$ and
\begin{align*}\label{eq renormalization estimate for first coordinate}
    \eta_{\hat h_1}&:= \left(F^{\underline{\theta}(\hat{h})\ua}\right)^{-1} \circ S_{\hat h}^{-1}\circ\eta_{\hat h} \circ F^{\underline{\theta}(\hat{h})\ua},\\
     \eta'_{\hat h'_1}&:= \left(S'^{-1}_{\hat h'}\circ S_{\hat h} \circ F^{\underline{\theta}(\hat{h})\ua}\right)^{-1} \circ S'^{-1}_{\hat h'}\circ \eta'_{\hat h'} \circ  \left(S'^{-1}_{\hat h'}\circ S_{\hat h} \circ F^{\underline{\theta}(\hat{h})\ua}\right ),\\
    B_{\hat h_1,\hat h'_1}&:=\left(F^{\underline{\theta}(\hat{h})\ua}\right)^{-1} \circ S_{\hat h}^{-1}\circ S'_{\hat h'}\circ B_{\hat h,\hat h'} \circ F^{\underline{\theta'}(\hat{h}') \ua'}.
\end{align*}
Here, $S_{\hat h}: = A_{\eta_{\hat h},p_{{\hat h}}}$ and $S'_{\hat h'}:= A_{\eta'_{\hat h'},p'_{{\hat h'}}}$ are affine estimates of 
$\eta_{\hat h},\eta'_{\hat h'}$ at 
$p_{{\hat h}}:= F^{\underline{\theta}(\hat{h}) \ua}(c_a),p'_{{\hat h'}}: =B_{\hat h,\hat h'}\left(F^{\underline{\theta'}(\hat{h}') \ua'}(c'_{a'})\right)$ where $a,a'$ are last letters of $\ua,\ua'$ respectively.

In order to estimate
the distances of $\eta_{\hat h_1},\eta'_{\hat h'_1}$ from $\id$ in 
$\cC^{1+\alpha}$ topology we use Lemma \ref{lemma estimation of chosing linear instead of non linear part} twice.
Note that the key observation here is that we can take the length of the words $\underline a$ and $\underline a'$ as long as we want such that the contraction rate of maps 
$ F^{\underline{\theta}(\hat{h}) \ua}$ and 
$F^{\underline{\theta'}(\hat{h}') \ua'}$ become exponentially small.
Applying Lemma \ref{lemma estimation of chosing linear instead of non linear part} with 
$B := F^{\underline{\theta}(\hat{h}) \ua}$, $\phi =\eta_{\hat h_1}$, $p := F^{\underline{\theta}(\hat{h}) \ua}(c_a)$ and 
$X: =F^{\underline{\theta}(\hat{h}) \ua} \left( k^{\underline \theta (\hat h_1)}(G(a))\right )$ in addition to inequalities derived in Corollary \ref{cor:operator bound of renormalization linear part} and  Lemma \ref{lemma introduction diameter of images} gives us that there exists bounded constants 
$C_{K,\eta_{\hat h}}, C_K$ depended on $\eta_{\hat h}$ and  $K$  such that 
\begin{align*}
\|\eta_{\hat h_1} - \id \|_{\cC^{1+\alpha}}
&= \|\left(F^{\underline{\theta}(\hat{h})\ua}\right)^{-1} \circ A_{\eta_{\hat h},p_{{\hat h}}}^{-1}\circ\eta_{\hat h} \circ F^{\underline{\theta}(\hat{h})\ua}-\id\|_{\cC^{1+\alpha}}\\ 
&\leq  C_{K,\eta_{\hat h}} \cdot \|\left(DF^{\underline{\theta}(\hat{h})\ua}\right)^{-1}\|_{\textit{op}} \cdot \|DF^{\underline{\theta}(\hat{h})\ua}\|_{\textit{op}}\cdot \diam (X)^{\alpha} 
\\
 &\leq C_{K,\eta_{\hat h}}\cdot C_K\cdot  \|Df_{\ua}^{-1}\|_{\it op} \cdot 
 \|Df_{\ua}\|_{\it op}
\cdot \diam (X)^{\alpha} \\
&
\leq  
C_{K,\eta_{\hat h}}  \cdot (C_K\cdot \kappa^{l_0}) \cdot  (C_K \cdot\mu^{l_0\cdot \alpha}),
\end{align*}
where $l_0$ is the length of $\underline a$. Given $\lambda<\frac{1}{4},\delta'>0$, the bunching condition $\kappa\cdot \mu^{\alpha}<1$ and uniform bound for $C_K,C_{K,\eta_{\hat h}}$ imply that there exist $l_1\in \N$ such that 
$C_K^2  \cdot 
C_{K,\eta_{\hat h}} \cdot (\kappa \cdot\mu^{\alpha})^{l_1}\leq \lambda\cdot \delta'$. Thus, length of $\ua$ being bigger than $l_1$ is sufficient to obtain that 
$$\|\eta_{\hat h_1} - \id \|_{\cC^{1+\alpha}}\leq \lambda\cdot \delta'.$$
Same calculations beside that $S_{\hat h},S'_{\hat h'}$ are affine maps with uniform bound concludes that (possibly by increasing $l_1$) if the length of $\ua'$ is bigger than $l_1$ then 
$$
\|\eta'_{\hat h'_1}-\id\|_{\cC^{1+\alpha}}\leq \lambda\cdot \delta'.
$$

\begin{figure}[ht]\label{fig: 3D. strong covering for the mais theorem.}
\def\svgwidth{\textwidth}
\centering
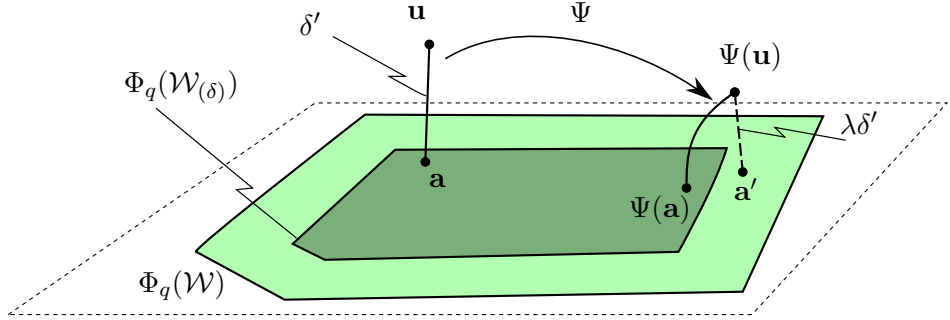
\caption{ $B_{\delta'}\left(\Phi_q(\cW)\right)$ satisfies strong covering condition with respect to the family $\cR$.
In this figure, 
${\bf{u}}:= [\eta_{\hat h} \circ k^{\underline{\theta}},
\eta'_{\hat h'} \circ B_{\hat h, \hat h'}\circ k^{\underline{\theta}'}]$, ${\bf{a}}:=[k^{\underline{\theta}}, B_{\hat h,\hat h'} \circ k^{\underline{\theta}'}]
$
and ${\bf{a}'}: = [k^{\underline{\theta}(\hat h_1)}, B_{\hat h_1,\hat h'_1} \circ k^{\underline{\theta}'(\hat h'_1)}]$. 
}
\end{figure}

Summarizing above, for $\delta'\in (0,\delta]$ and $\lambda<\frac{1}{4}$ there is $l>l_1>0$ such that for 
${\bf{a}}:=[k^{\underline{\theta}}, B_{\hat h,\hat h'} \circ k^{\underline{\theta}'}]\in \Phi_q(B_{\delta}(\cW))$ and 
${\bf{u}}:=[\hat h,\hat h']= [\eta_{\hat h} \circ k^{\underline{\theta}},
\eta'_{\hat h'} \circ B_{\hat h, \hat h'}\circ k^{\underline{\theta}'}]$ 
in the $\delta'$ neighborhood of $\Phi_q(\cW)$, there exists $\Psi:= T_{\ua,\ua'}\in \cR_l$ with lengths of words $\ua,\ua'$ both less than or equal to $l$ such that $\operatorname{dist}_{\cC^{1+\alpha}}(\Psi({\bf{u}}),{\bf{a}'})\leq \lambda \cdot \delta'$, where
$$\Psi({\bf{u}}) := [\eta_{\hat h_1}\circ k^{\underline{\theta}(\hat h_1)},\eta'_{\hat h'_1}\circ B_{\hat h_1,\hat h'_1}\circ  k^{\underline{\theta}'(\hat h'_1)}],\; 
{\bf{a}'}: = [k^{\underline{\theta}(\hat h_1)}, B_{\hat h_1,\hat h'_1} \circ k^{\underline{\theta}'(\hat h'_1)}].
$$
It only remains to choose the constant $\delta'>0$ and prove ${\bf{a}'}\in \Phi_q(\overline{\cW})$. Since $\Phi_q(\overline{\cW})$ is compact and renormalization operators $\Psi\in \cR$ with bounded length of words are finite, there exists 
$\delta_2<\delta/4$ such that for any ${\bf{b}}\in \overline{\cW}$, whenever $\operatorname{dist}_{\cC^{1+\alpha}}({\bf{b}},{\bf{v}})<\delta_2$ for some ${\bf{v}}\in \cQ$ then $\operatorname{dist}_{\cC^{1+\alpha}}(\Psi({\bf{b}}),\Psi({\bf{v}}))\leq \delta/4$ for all $\Psi\in \cR_l$. Denote 
\begin{equation}\label{eq: definition of delta' proof of the main theorem}
    \delta' : = \min \{\delta_1,\delta_2,\delta\}.
\end{equation}
By triangle inequality we have
\begin{align*}
\operatorname{dist}_{\cC^{1+\alpha}}({\bf{a}'},\Psi({\bf{a}})) & \leq \operatorname{dist}_{\cC^{1+\alpha}}({\bf{a}'},\Psi({\bf{u}}))+\operatorname{dist}_{\cC^{1+\alpha}}(\Psi({\bf{u}}),\Psi({\bf{a}}))  \\  
& \leq \lambda\cdot \delta' + \delta/4 < \delta/2.
\end{align*}
This implies that ${\bf{a}'}\in \Phi_q(\cW)$ since $\Psi({\bf{a}}) \in \Phi_q(\cW_{(\delta)})$ by \eqref{eq delta interior of renormalization process}.
\end{proof}

\begin{remark}\label{remark:bounded geomtery is true for affine Cantors} If we restrict ourselves to the setting of affine Cantor sets and their affine configurations, as the proofs throughout this and previous sections show, the bunching condition \eqref{eq: definition of bunched cantor sets} is not necessary. Indeed,
    in this case limit geometries do exist and are identity maps as discussed in \S \ref{sec: Affine Cantor sets} regardless of whether they satisfy the bunching condition or not.
    Therefore, Corollary \ref{corol:  intersecting pairs in affine configs inside a covering region} can be applied in this case; implying the stability of intersection within the space of affine configurations of affine Cantor sets. Thus, Theorem \ref{thm:A} has an analogous version for general affine Cantor sets and their affine perturbations even without the bunching assumption.
\end{remark}

\section{Examples of stably intersecting Cantor sets}\label{sec: Affine QC-Cantor sets} 
In this section, we provide explicit examples of pairs of Cantor sets having stable intersection by showing that their corresponding renormalization operators satisfy the covering criterion in Theorem \ref{thm: covering criterion for intersection}. Moreover, one of the Cantor sets can have arbitrarily small Hausdorff dimension. 

\begin{definition}[Expanding $n$-cover]\label{def: n-cover of expanding operators}
Let $(K,K')$ be a pair of bunched (or bunched holomorphic) Cantor sets, $\cD_1^*,\cE_1^*\subset \cR_1^*$ be the families of contracting and expanding renormalization operators of the pair $(K,K')$ respectively (see \eqref{eq: contracting and expanding renormalization operators definition}) and $\cW\subset \cQ_{\Aff}$ be an open relatively compact set. We say that $\cW$ has an \emph{expanding $n$-cover} for some $n>1$ with respect to $\cR_1^*=\cD_1^*\sqcup \cE_1^*$ if there are disjoint subsets $\cE_1,\cE_2,\dots,\cE_n\subset \cE_1^*$ such that $\cW$ satisfies the strong covering condition with respect to the
family of operators $\cD_1^* \sqcup \cE_i$ for $i=1,\dots,n$.
\end{definition}

Our construction of stably intersecting Cantor sets in $\R^d$ has the following steps. 
\begin{itemize}
    \item[(1)] Introducing a pair of affine Cantor set $(K_1,K_1')$ in $\R$ and a bounded open set $W_1\subset \Aff(1,\bR)$ such that $W_1$ has an expanding $3$-cover with respect to the pair $(K_1,K_1')$. Moreover, $\dimHD(K_1)$ can be arbitrarily close to 0. 
    \item[(2)] Constructing 
    a bounded open set $W_d\subset \Aff_{\id}(d,\R)\cong \bR^*\times \bR^d$ such that $W_d$ has expanding 
    $3^d$-cover with respect to the action of renormalization operators of the pair $(K_1^d,K_1'^d)$ restricted to the subgroup 
    $\Aff_{\id}(d,\R)$ of $\Aff(d,\bR)$ (see Lemma \ref{lem: affine left and right action details} for the definition of $\Aff_{\id}(d,\mathbb{F})$).
    
    \item[(3)]  Perturbing the generators of $K_1^d,K_1'^d$ to obtain Cantor sets $K_d,K'_d$ such that for an open set $\cU \subset \SL(d,\R)$ in the neighborhood of
    $\id$
    the set 
    $\cW:=W_d\times \cU\subset \Aff(d,\bR)$ satisfies strong covering condition with respect renormalization operators of the pair 
    $(K_d,K'_d)$.
\end{itemize}

\subsection{Cantor sets in $\R$ 
with expanding $3$-cover}\label{sec: cantor sets in dimension one}
We begin with the construction of some Cantor sets in $\R$. Let  $N\geq 7$ be an integer, and $\tau>0$ be a small enough constant which will be determined later. Define the following maps on the real line
\begin{align*}
 g_{1,1}&: x\mapsto \frac{x}{N+\tau}, 
& g_{1,2}&:x\mapsto \frac{x-1}{N+\tau}+1,\\  
 g_{2,1}&: x\mapsto \frac{x}{N+\tau}+\dfrac{1}{N}, &
g_{2,2}&:x\mapsto \frac{x-1}{N+\tau}+\dfrac{N-1}{N},\\
 g_{3,1}&:x\mapsto \dfrac{x}{N+\tau}+\dfrac{2}{N}, & g_{3,2}&: x\mapsto \frac{x-1}{N+\tau}+\dfrac{N-2}{N},\\ 
 f_1&:x\mapsto \frac{1}{2+\tau}x,&  f_2&: x\mapsto  \frac{x-1}{2+\tau}+1.
\end{align*}

We denote the Cantor sets $K_1$ and $K_1'$ as the invariants sets of the IFSs generated by the families of maps $\{g_{1,1},g_{2,1},g_{3,1},g_{1,2},g_{2,2},g_{3,2}\}$ and $\{f_1,f_2\}$, respectively. Hence, 
\begin{align}\label{eq: hausdorff dimension of K1,K1'}
    \dimHD(K_1)=\dfrac{\ln 6}{\ln (N+\tau)},\quad
\dimHD(K_1')=\dfrac{\ln 2}{\ln(2+\tau)}.
\end{align}
Therefore, we can take $\dimHD(K_1)$ as small as we want by taking $N$ large enough.

~\medskip
\begin{figure}[ht]
    \centering
    \includegraphics[width=0.86\textwidth]{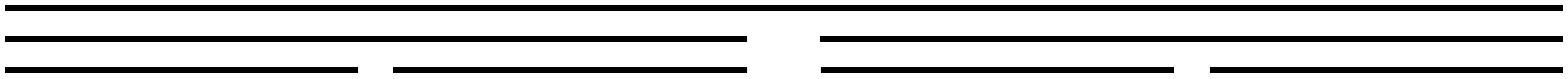}
    \caption{Cantor set $K_1'$, approximation in first and second steps.}\label{figure: 1D cantor K1'}
\end{figure}
~\medskip
\begin{figure}[ht]
    \centering
    \includegraphics[width=0.86\textwidth]{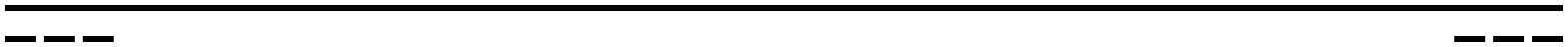}
    \caption{Cantor set $K_1$, approximation in first step.}\label{figure: 1D cantor K1}
\end{figure}

Since $K_1$ and $K_1'$ are affine Cantor sets, the corresponding family of renormalization operators $\cR^*$ has an action on the group $\Aff(1,\R)$ (see \S \ref{sec: Affine Cantor sets}).
We will present an open relatively compact set 
$W_1\subset \Aff(1,\R) \cong \R \times \R^*$  which satisfies the strong covering condition with respect to the action of the family of renormalization operators $\cR_1^*$ of the pair $(K_1,K_1')$ on the space 
$\Aff(1,\R)$. Moreover, we shall show that $W_1$ has an expanding 3-cover with respect to the pair $(K_1,K'_1)$. 
For this aim, we define the following operators acting on 
$\Aff(1,\R)\cong \R\times \R^*$.
\begin{align*}
G_{1,1}&: (s,t) \mapsto (Ns,Nt), &G_{1,2}&:(s,t)\mapsto (Ns,Nt-(N-1)),\\
G_{2,1}&: (s,t)\mapsto (Ns,Nt-1), &G_{2,2}&:(s,t)\mapsto (Ns,Nt-(N-2)),\\
G_{3,1}&: (s,t)\mapsto (Ns,Nt-2), &G_{3,2}&:(s,t)\mapsto (Ns,Nt-(N-3)),\\
      F_1&: (s,t)\mapsto (\dfrac{s}{2},t),&F_2&: (s,t)\mapsto (\dfrac{s}{2},t+\dfrac{s}{2}).
  \end{align*}
  Note that $G_{i,j}$ is the limit of the expanding operator $g_{i,j}^{-1}$ and $F_i$ is the limit of the contracting operator $f_i$ where $\tau \to 0$. So, these operators are in fact the limit of renormalization operators of the pair $(K_1,K'_1)$ when $\tau\to 0$. 

Let $W_1$ be the interior of closed convex hull of points ${\bf a},{\bf b},{\bf c},{\bf d}$ with coordinates
\begin{align*}
&{\bf a} := \left(a,\dfrac{2}{N-1}-a+\delta\right), &&{\bf b}:= \left(a,\dfrac{N-3}{N-1}-\delta\right),\\
&{\bf c}:=\left(b,\dfrac{N-3}{N-1}-2\delta\right),
&&{\bf d}: = \left(b,\dfrac{2}{N-1}-b+2\delta\right),
\end{align*}
where $\delta>0$ is sufficiently small positive number and $a:=1, b := 2N+1$.

\begin{figure}
    \centering
    \includegraphics{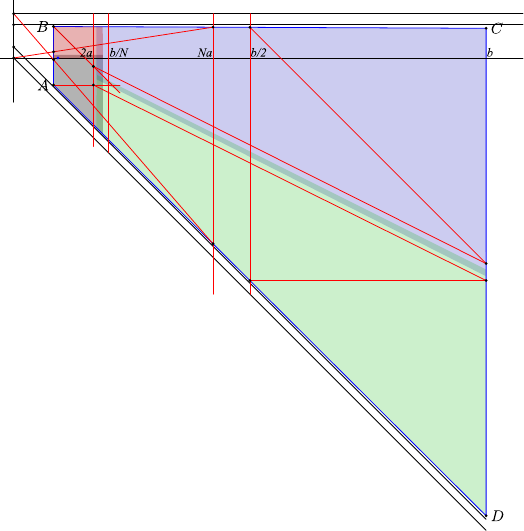}
    \caption{The open set $W_1$.
    The red lines show how to find the colored regions $P'_1,P'_2,P_2,P_1$.}\label{figure: colored regions. idea of proof}
\end{figure}

\begin{proposition}\label{prop: multi covering for K1,K2,K3}
    The open set $W_1\subset \Aff(1,\R)$ satisfies the covering condition with respect to the families $\{F_1,F_2,G_{j,1},G_{j,2}\}$ for $j\in \{1,2,3\}$. More precisely, there are polygons 
    $P_1,P_2,P_1',P_2'\subset \R\times \R^*$ covering $\overline{W_1}$ such that for each $j\in \{1,2,3\}$, 
    \begin{align}\label{eq: example double recurrence W1}
    P_1\cup P_2 \cup P_1'\cup P_2' \subset \bigcup_{i=1,2} F_i^{-1}(W_1)\cup G_{j,i}^{-1}(W_1).
    \end{align}
\end{proposition} 
\begin{proof} We prove that for any $j\in \{1,2,3\}$ and $i\in\{1,2\}$, 
\begin{equation}\label{eq: expanding 3 cover W1 example}
  F_i(P_i') \cup G_{j,i}(P_i)
 \subset W_1.
\end{equation}
Clearly, \eqref{eq: expanding 3 cover W1 example} implies  \eqref{eq: example double recurrence W1}.
    As illustrated in the Figure \ref{figure: colored regions. idea of proof}, $P_i$ and $P_i'$ are all compact convex polygons for $i=1,2$ so we will define them with their vertices. We will use the notation $\overline{\operatorname{ch}}({\bf a},{\bf b},{\bf c},{\bf d})$ to describe the closed convex hull of ${\bf a},{\bf b},{\bf c},{\bf d}\subset \R\times \R^*$. For $i=1,2$ we shall describe coordinates of points ${\bf a}_i,{\bf b}_i,{\bf c}_i,{\bf d}_i$ and ${\bf a}_i',{\bf b}_i',{\bf c}_i',{\bf d}_i'$ to define $P_i := \overline{\operatorname{ch}}({\bf a}_i,{\bf b}_i,{\bf c}_i,{\bf d}_i)$ and $P_i' := \overline{\operatorname{ch}}({\bf a}_i',{\bf b}_i',{\bf c}_i',{\bf d}_i')$. Let $\gamma>0$ be sufficiently smaller than $\delta$ and define
    \begin{align*}
      &{\bf a}_1 :=(a,\dfrac{2}{N-1}-a+\dfrac{\delta}{2}), &&{\bf b}_1 := \left(a, \dfrac{(N-4)}{N(N-1)}\right),\\
      &{\bf c}_1:= \left(\dfrac{b}{N}-\gamma, \dfrac{(N-4)}{N(N-1)}\right),
    &&{\bf d}_1:=(\dfrac{b}{N}-\gamma,\dfrac{2}{N-1}+\dfrac{\delta}{2}-\dfrac{b}{N}+\gamma),\\
    &{\bf a}_2: = \left(a,\dfrac{-(N-4)}{N(N-1)}\right),&&{\bf b}_2 := \left(a,\dfrac{N-3}{N-1}-\dfrac{\delta}{2}\right),\\
        &{\bf c}_2: = \left(\dfrac{b}{N}-\gamma,\dfrac{N-3}{N-1}-\dfrac{\delta}{2}\right),
    &&{\bf d}_2: = \left(\dfrac{b}{N}-\gamma,\dfrac{-(N-4)}{N(N-1)}\right),\\
    &{\bf a}_1':= \left(2a+\gamma,\dfrac{2}{N-1}-a+2\delta-\dfrac{\gamma}{2}\right), &&{\bf b}_1' := \left(2a+\gamma,\dfrac{N-3}{N-1}-\delta-\dfrac{\delta(a+\gamma)}{b-a}\right),\\
    &{\bf c}_1':=\left(b,\dfrac{N-3}{N-1}-2\delta\right),
    &&{\bf d}_1':= \left(b,\dfrac{2}{N-1}-\dfrac{b}{2}+2\delta\right),\\
    &{\bf c}_2':=\left(b,\dfrac{N-3}{N-1}-2\delta-\dfrac{b}{2}\right),
    &&{\bf d}_2': = \left(b, \dfrac{2}{N-1}-b+2\delta\right),
\end{align*}
\begin{align*}
    &{\bf a}_2' := \left(2a+\gamma,\dfrac{2}{N-1}-2a-\gamma+\delta+\dfrac{\delta(a+\gamma)}{b-a}\right),\\
    &{\bf b}_2': = \left(2a+\gamma,\dfrac{N-3}{N-1}-2\delta-a-\dfrac{\gamma}{2}\right).
    \end{align*}
   Given a polygon $P$ and a convex open set $Q$, to show that $\overline{P}\subset Q$ it suffices to prove that the vertices of $P$ lie inside the $Q$. 
    Therefore, observing that $F_i,G_{i,j}$ are all affine maps from $\R\times \R^*$ to $\R^2$, the relation \eqref{eq: expanding 3 cover W1 example} satisfies because by the definitions above for each $j\in \{1,2,3\}$ we have 
    \begin{align*}
       F_i({\bf a}_i'),F_i({\bf b}_i')&,F_i({\bf c}_i'),F_i({\bf d}_i') \subset W_1,\\ 
       G_{j,i}({\bf a}_i),G_{j,i}({\bf b}_i)&,G_{j,i}({\bf c}_i),G_{j,i}({\bf d}_i)\subset W_1,
    \end{align*}
    for $i=1,2$. Note that $W_1$ is surrounded by 4 lines. We can describe $W_1$ by 
    $$W_1= \left\{(s,t): a< s< b,\;\; \dfrac{2}{N-1}-s+\delta+\dfrac{\delta(s-a)}{b-a}< t< \dfrac{N-3}{N-1}-\delta-\dfrac{\delta(s-a)}{b-a}\right\}.$$
    This implies that $\overline{W_1}\subset P_1\cup P_2\cup P_1'\cup P_2'$. 
    \end{proof}
    \begin{remark}\label{remark: vertical lines multi cover}
                        Note that in \eqref{eq: expanding 3 cover W1 example} contracting maps $F_i$ are acting on the regions $P_i'$, while expanding maps $G_{1,i},G_{2,i},G_{3,i}$ are acting on the regions $P_i$ for $i=1,2$. 
    An important property of these regions is that any vertical line intersecting $\overline W_1$ either intersects $P_1\cup P_2$ or $P_1'\cup P_2'$. Thus, \eqref{eq: expanding 3 cover W1 example} implies that we can map all points on a vertical segment inside  $\overline{W_1}$ into $W_1$ via only expanding maps or by only contracting maps. 
    Indeed, all points on such a vertical segment returns back into $W_1$ by either the contracting map $F_i$ or each of the expanding maps $G_{1,i},G_{2,i},G_{3,i}$, for some $i=1,2$. This observation will play a crucial role in our construction of Cantor sets in higher dimensions.
\end{remark}

It follows from Lemma \ref{lemma covering to strong} that for a relatively open compact set and a finite family of operators the covering condition \eqref{eq covering condition} is equivalent to strong covering \eqref{eq: covering condition delta interior}. Moreover, the strong covering \eqref{eq: covering condition delta interior} with respect to a finite family is stable under small perturbations (see Lemma \ref{lemma stability of strong covering}). Hence, by Proposition \ref{prop: multi covering for K1,K2,K3}, $W_1\subset \Aff(1,\R)$ satisfies strong covering condition with respect to the families of affine maps 
$\{f_1,f_2,g_{j,1}^{-1},g_{j,2}^{-1}\}$
for $j\in \{1,2,3\}$ when $\tau>0$ is small enough. Action of renormalizaion operators $\cR_1^*$ on $\Aff(1,\R) \cong \R \times \R^*$ for the pair $(K_1,K_1')$ consists of two contracting maps $f_1,f_2$ which are the generators of $K_1'$ and six expanding maps 
$g_{1,1}^{-1},g_{1,2}^{-1},g_{2,1}^{-1},g_{2,2}^{-1},g_{3,1}^{-1},g_{3,2}^{-1}$ obtaining from $K_1$. Consider the partition of expanding renormalization operators of the pair $(K_1,K_1')$ into two families 
$\cA_i:=\{g_{1,i}^{-1},g_{2,i}^{-1},g_{3,i}^{-1}\}$ for $i=1,2$. Consequently, fro sufficiently small enough $\tau$ we have the following.
\begin{corollary}\label{corol: W1 has expanding 3 cover, example part 1}
     $W_1$ 
     has expanding $3$-cover with respect to the pair $(K_1,K_1')$. Indeed, any ${\bf x}\in \overline{W_1}$ returns back into $W_1$ by the contracting operator $F_i$ when ${\bf x}\in P_i'$ or by each of the expanding operators in the family $\cA_i$ whenever ${\bf x}\in  P_i$, for $i\in\{1,2\}$.
\end{corollary}

\subsection{Cantor sets in $\R^d$ with expanding $3^d$-cover on $\Aff_{\id}(d,\R)$}
It is a simple observation that if $K$ is a regular Cantor set in $\bF$ generated by an IFS with full shift symbolic type consisting of $n$ contracting maps $f_1,f_2,\dots,f_n$, then $K^d: = K\times K\times \cdots \times K$ is a regular Cantor set in $\bF^d$ generated by an IFS consisting of $n^d$ maps 
$f_{i_1}\times  f_{i_2}\times \cdots \times f_{i_d}$.
In the case that $K$ is an affine Cantor set, the map $f_{i_1}\times  f_{i_2}\times \cdots \times f_{i_d}: \bF^d \to \bF^d$ is an affine map of the form $[x\mapsto Ax+v]$ where $A\in \GL(d,\bF),v\in \bF^d$. Hence, $K^d$ is an affine Cantor set generated by an IFS consisting of $n^d$ affine maps with full shift symbolic type. Furthermore, $\dimHD(K^d) = d\cdot \dimHD(K)$ (see \cite{Palis-Takens_book93}).

\begin{figure}[ht]\label{figure: 2D cantors K,K'}
    \centering
    \includegraphics[width=0.26\textwidth]{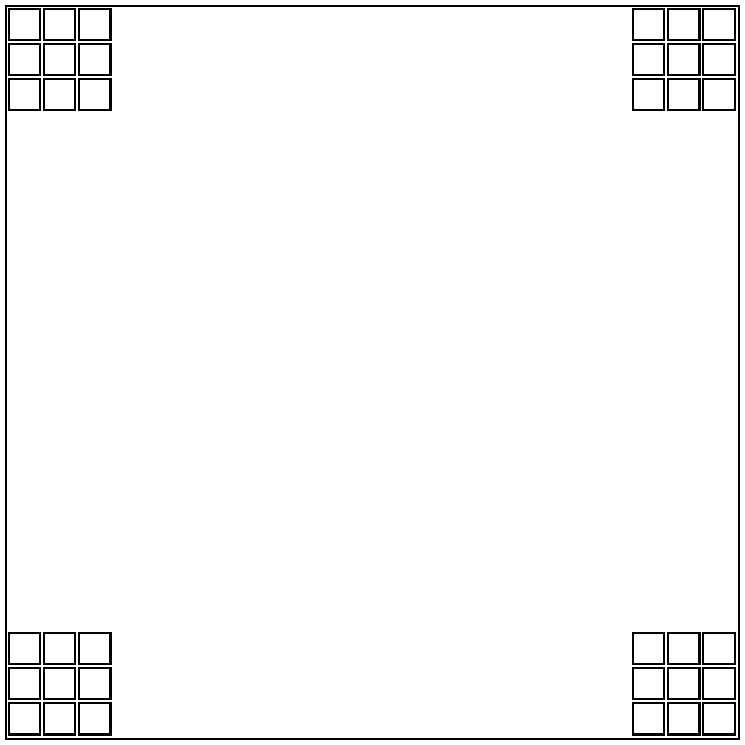}
    \quad\includegraphics[width=0.26\textwidth]{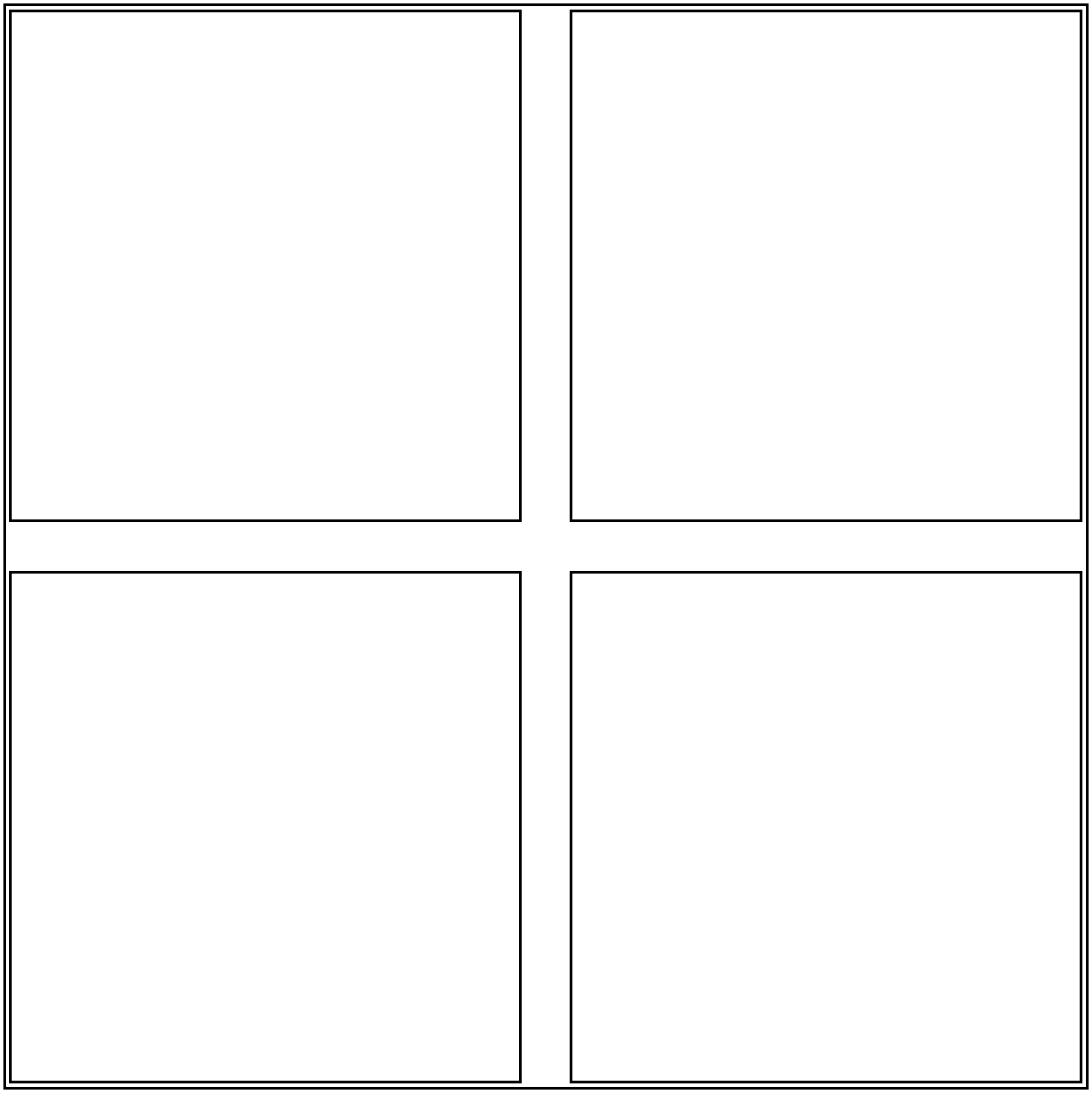}
    \caption{Cantor sets $K_1^d, K_1'^d$ in dimension $d=2$, first step approximation.}
\end{figure}
Each of the affine generators of the affine Cantor sets $K_1^d,{K'_1}^d$ is equal to composition of a homothety and a translation which is like $[x\mapsto s \cdot x+b]\subset \Aff_{\id}(d,\R)$ (see Lemma \ref{appendix affine actions} for the definition of $\Aff_{\id}(d,\R)$).  Because the generators $f_i$ of $K_1'$ have same contraction rate $1/(2+\tau)$ which implies that the generators of ${K'_1}^d$ are members of $\Aff_{\id}(d,\R)$. Similarly, the generators of $K_1^d$ belong to $\Aff_{\id}(d,\R)$.
Thus, the renormalization operators of the pair $(K_1^d,K_1'^d)$ belong to $\Aff_{\id}(d,\R)$. So, we can study their action on this subgroup.

Denote $W_d\subset \Aff_{\id}(d,\R)$ as the fiber product of $d$ copies of $W_1$ over the base $\R^*$ (See Proposition \ref{prop: multi covering for K1,K2,K3}). More precisely, 
$$W_d : = \left\{(s,t_1,t_2,\dots,t_d) : (s,t_i)\in W_1,\;\; \forall i=1,2,\dots,d \right\}\subset \R^d \rtimes \R^* \cong \Aff_{\id}(d,\R) .$$

\begin{lemma}\label{lem: Wd satisfies covering} 
$W_d$ has expanding $3^d$-cover with respect to the family of renormalization operators of the pair $(K_1^d,K_1'^d)$. Moreover, 
there exists $\varepsilon_1,\varepsilon_2>0$ such that
for any $(w,u)\in  B_{\varepsilon_2}(W_d) \times \overline{B_{\varepsilon_1}(\id)}\subset \Aff_{\id}(d,\R)\times \SL(d, \R)$ either there exists a contracting operator or at least $3^d$ expanding operators that maps $(w,u)$ to 
$(w',u)\in( W_d)_{(\varepsilon_2)}\times  \overline{B_{\varepsilon_1}(\id)}$.
\end{lemma}
\noindent Here, by $B_{\varepsilon_1}(\id)$ we mean the $\varepsilon_1$ neighborhood of $\id\in \SL(d,\R)$. Also, $B_{\delta}(W_d)$ and $(W_d)_{(\delta)}$ are $\delta$-neighborhood and $\delta$-interior of $W_d$ in $\Aff_{\id}(d,\R)$ which should not be confused with the ambient space $\Aff(d,\R)\supset W_d$.

\begin{proof} Let $w = (s,t_1,\dots,t_d)\in \overline{W_d}$. For each $j=1,2,\dots, d$ there is a renormalization operator 
 $\Psi_j$ of the pair $(K_1,K_1')$ that maps each $(s,t_j)$ to $(s_j,t_j')\in W_1$. Denote the vertical line $L_s$ inside $\overline{W_1}$ as 
 $$L_s:= \{(s,t): (s,t)\in \overline{W_1}\}\subset \overline{W_1}$$
 By Remark \ref{remark: vertical lines multi cover}, we have either $L_s\subset P'_1\cup P'_2$ or $L_s\subset P_1\cup P_2$. 
 Subsequently, $\Psi_j$'s are either all contracting operators or all expanding operators. When $L_s\subset P'_1\cup P'_2$, 
 $\hat \Psi := \Psi_1 \times \Psi_2 \times \cdots \times \Psi_d\in \Aff_{\id}(d,\R)$ will be a contracting renormalization operator of the pair $(K_1^d,K_1'^d)$ which maps $w$ into $W_d$. If $L_s\subset P_1\cup P_2$, then by Corollary \ref{corol: W1 has expanding 3 cover, example part 1} we conclude that there are  least three options for each of $\Psi_j$, $j=1,\dots,d$.
Therefore, there are at least $3^d$ expanding renormalization operators of the form 
$\hat \Psi =\Psi_1 \times \Psi_2 \times \cdots \times \Psi_d\in \Aff_{\id}(d,\R)$ which map $w$ into $W_d$.
     This implies that $W_d$ has expanding $3^d$-cover with respect to the family of the renormalization operators of the pair $(K_1^d,K_1'^d)$. By Corollary \ref{corol: W1 has expanding 3 cover, example part 1}, we can partition the family of expanding renormalization operators of the pair 
    $(K_1^d,K_1'^d)$ which consists of $6^d$ elements as disjoint union 
    $\bigsqcup_{\alpha \in \{1,2\}^d} \cH_{\underline{\alpha}}$
    of $2^d$ sets such that for $\underline{\alpha}=(\alpha_1,\dots,\alpha_d)$,
    \begin{align}\label{eq: partition Hi}
    \cH_{\underline{\alpha}}:=\{\Psi_1\times \Psi_2 \times\cdots \times\Psi_{d} ~:~ \Psi_j\in \cA_{\alpha_j},\; 1\leq j\leq d\}.
    \end{align}
    This partitioning is such that any $w\in \overline{W_d}$ maps into $W_d$ by either a contracting operator or any expanding operator in the set $\cH_{\underline{\alpha}}$ for some $\underline{\alpha}\in \{1,2\}^d$. Lemma \ref{lemma covering to strong} implies that there is $\delta>0$ such that any $w\in B_{\delta}(W_d)\subset \Aff_{\id}(d,\R)$ maps to $(W_d)_{(\delta)}$ via either a contracting operator or any expanding one from the set $\cH_{\underline{\alpha}}$ for some $\underline{\alpha}\in \{1,2\}^d$.
    Corollary \ref{corol: action of renormalization is product and semi direct product} helps us to study the action of renormalization operators of the pair $(K_1^d,K_1'^d)$ on $\overline{W_d}\times \overline{B_{\varepsilon}(\id)}$. Let $(w,u)\in B_{\delta}({W_d})\times \overline{B_{\varepsilon}(\id)}$.
    In the case that $w$ maps to $(W_d)_{(\delta)}$ by an expanding operator $\Psi\in\cH_{\underline{\alpha}}$ for some $\underline{\alpha}\in \{1,2\}^d$,
    because $\Psi\in \Aff_{\id}(d,\R)$ and the family of expanding operators act on $\Aff(d,\R)$ as the product action on 
    $\R^d\times (\R^*\times \SL^\star(d,\R))$ then the image of $(w,u)$ is $(\Psi(w),u)\in (W_d)_{(\delta)}\times \overline{B_{\varepsilon}(\id)}$.
    In the other case, when $w\in B_{\delta}({W_d})$ is mapped into $(W_d)_{(\delta)}$ by a contracting operator $\Psi$, then the image of $(w,u)$ can be written as $(w'_u,u)$.
    Note that $w'_{\id}=\Psi(w)\in (W_d)_{(\delta)}$. But for $u\neq \id$,  $w'_u$ may differ from $\Psi(w)$, since the action of contracting operators on $\Aff(d,\R)$ is semidirect product of its action on $\R^d$ and $\GL(d,\R)$.
    However, $w'_u$ varies continuously with respect to $u\in \SL(d,\R)$.
    Therefore, there exists $\varepsilon_1>0$ such that for any $u\in \overline{B_{\varepsilon_1}(\id)}\subset \SL(d,\R)$, $w'_u\in (W_d)_{(\delta/2)}$.
It suffices to denote $\varepsilon_2:= \delta/2$
to yield the result.
\end{proof}

\begin{remark}\label{remark: finitely many consequtive contraction}
    Since $W_d\subset \R^d\times \R^*$ is a bounded set, a point cannot remain in $W_d$ solely by the action of contracting operators, nor solely by the action of expanding ones. In particular, there exists $T>0$ such that for any $w\in \overline{W_d}$ and any sequence $\Psi_1,\Psi_2,\dots,\Psi_T$ of the operators of the pair $(K_1^d,{K_1'}^d)$ such that 
    $$\Psi_j\circ \Psi_{j-1}\circ \cdots \circ \Psi_1(w)\in W_d$$
    for all $1\leq j\leq T$, then there are both expanding and contracting operators among $\Psi_1,\dots,\Psi_T$.
\end{remark}

\subsection{Examples of stable intersection in any dimension
}
We are ready to construct a pair of bunched Cantor sets $(K_d,K_d')$ in $\R^d$ satisfying the covering criterion. 

\begin{theorem}\label{thm: there exists covering criterion in Rd}
There exists a pair $(K_d,K'_d)$ of bunched affine Cantor sets with full shift symbolic type in $\R^d$ arbitrarily close to the pair $(K_1^d,K_1'^d)$ and an open relatively compact set $\cU\subset \SL(d,\R)$ such that $\cW: = W_d \times \cU\subset \Aff(d,\R)$ satisfies the strong covering condition with respect to the action of a finite family of renormalization operators of the pair $(K_d,K'_d)$ on the group $\Aff(d,\R)$. 
\end{theorem}

\begin{proof}[Proof of Theorem \ref{thm: there exists covering criterion in Rd}] Let $\varepsilon<\varepsilon_1$ be a positive number where $\varepsilon_1$ is from Lemma \ref{lem: Wd satisfies covering} and $\delta>0$ be a small  constant which will be determined later.
According to Lemma \ref{lem: generating a convering in SL with d^2 maps} there exists an open relatively compact set $\cU \subset B_{\varepsilon}(\id)\subset \SL(d,\R)$ and $d^2$ matrices $M_1,\dots,M_{d^2}\in B_{\delta}(\id) \subset \SL(d,\R)$ such that $\cU$ satisfies the (strong) covering condition with respect to action of the family $\{M_1^{-1},\dots,M_{d^2}^{-1}\}$ on $ \SL(d,\R)$. 
  
  Consider the families of expanding renormalization operators $\{\cH_{\underline{\alpha}}\}_{\underline{\alpha}\in \{1,2\}^d}$ of the pair 
  $(K_1^d,K_1'^d)$ defined in \eqref{eq: partition Hi}. These are inverses of the generators of $K_1^d$ as has been shown in \S \ref{sec: Affine Cantor sets}. Thus, their $\SL$ parts are all $\id\in \SL(d,\R)$. Now we perturb the affine generators of the Cantor set $K_1^d$ to obtain the Cantor set $K_d$.  This perturbation is such that for each $\underline{\alpha}\in \{1,2\}^d$ the $\SL$ part of $d^2$ generators among the generators corresponded to operators in the family $\cH_{\underline{\alpha}}$
  be equal to $\{M_1,M_2,\dots,M_{d^2}\}$. This perturbation can be done because $3^d>d^2$ and implies that 
  for each $\underline{\alpha}\in \{1,2\}^d$ 
  there are $d^2$ operators in $\cH_{\underline{\alpha}}$ with $\SL$ component equal to 
  $M_1^{-1},M_2^{-1},\dots,M_{d^2}^{-1}$.
  It is important to observe that we do not change determinant or the translation component of the affine generators.
  The difference of the affine generators of $K_d$ and $K_1^d$ is just in their $\SL$ component.
  Note that 
  $K_1^d$ and $K_1'^d$ are both conformal Cantor sets. So, their sufficiently small perturbations will be bunched Cantor sets.

\begin{figure}[ht]
    \centering
    \includegraphics[width=0.7\textwidth]{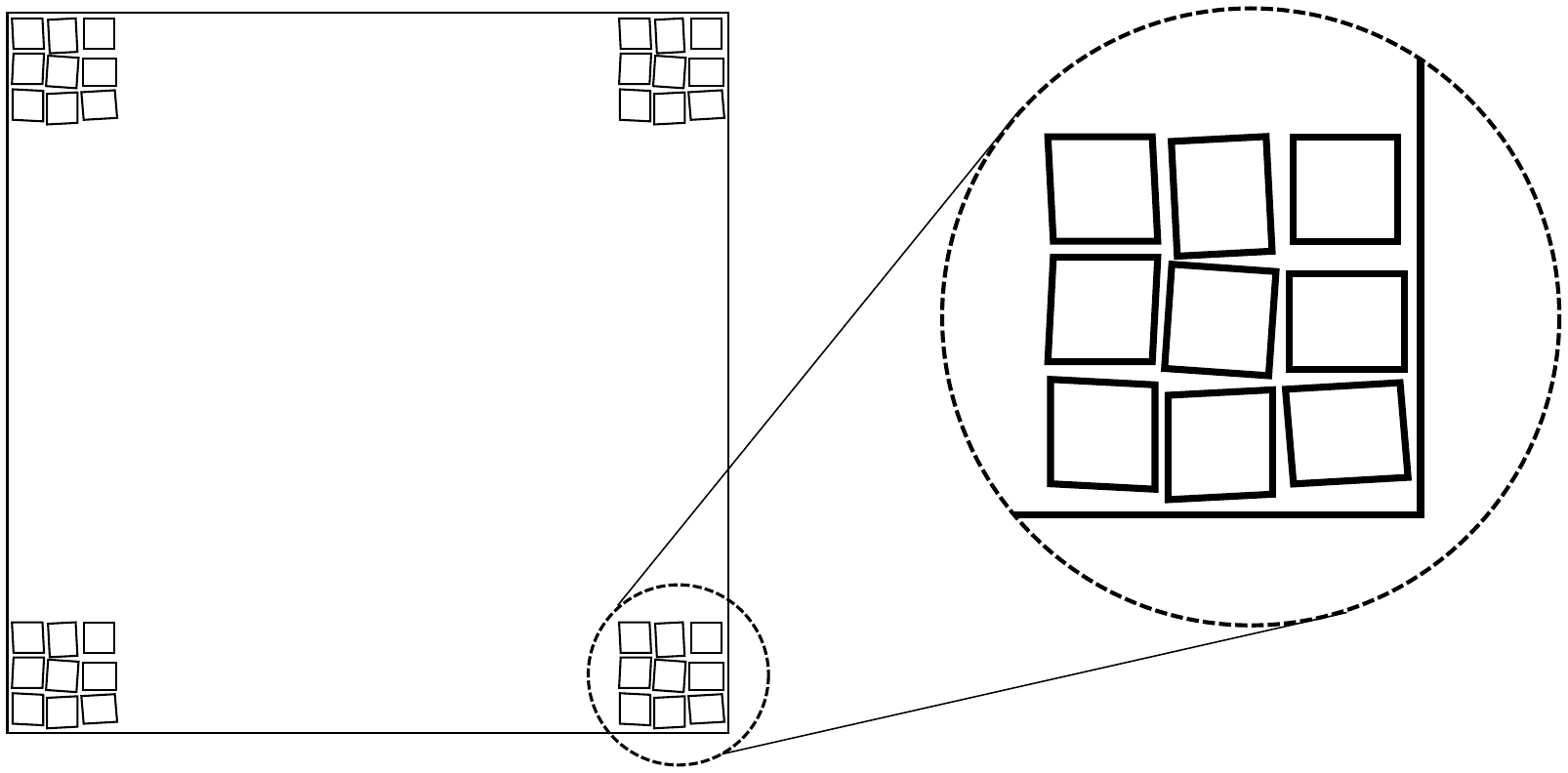}
    \caption{Cantor set $K_d$ for $d=2$ is illustrated in first step approximation.}\label{figure: 2D cantor K'}
\end{figure}

\noindent   Let $K_d':= K_1'^d$. 
   The contracting operators of the pair $(K_1^d,K_1'^d)$ coincides with the ones of  $(K_d,K_d')$, since they are the generators of $K_1'^d=K_d'$.
   The expanding operators of the pairs $(K_d,K_d')$ and $(K_1^d,K_1'^d)$ belong to $(\R^d\times \R^*)\times \SL(d,\R)$ which differ only on their $\SL$ part.
   By \eqref{eq: expanding action computation}, small changes on $\SL$ part of expanding operators affect small changes on their action on the subgroup $\Aff_{\id}(d,\R)\cong \R^d\times \R^*$.
    So, for small enough $\delta>0$ the action of expanding operators of the pairs $(K_d,K'_d)$ and $(K_1^d,{K'_1}^d)$ on $\Aff_{\id}(d,\R)$ are $\varepsilon_3$-close where $\varepsilon_3:=\varepsilon_2/2$ and $\varepsilon_2$ is obtained from Lemma \ref{lem: Wd satisfies covering} for the pair $(K_1^d,{K'_1}^d)$.
   
   We claim that any 
   $(w,u)\in \overline{W_d}\times \overline{\cU}\subset B_{\varepsilon_2}(W_d)\times \overline{B_{\varepsilon}(\id)}$ maps to $(W_d)_{(\varepsilon_3)}\times \cU$ via some renormalization operator of the pair $(K_d,K'_d)$.
  According to Lemma \ref{lem: Wd satisfies covering}, there is either a contracting operator of the pair $(K_1^d,K_1'^d)$ that maps $(w,u)$ to some $(w'_u,u)\in (W_d)_{(\varepsilon_2)} \times \overline{\cU}$ or an index $\underline{\alpha}\in \{1,2\}^d$ such that maps $(w,u)$  to $(w',u)\in (W_d)_{(\varepsilon_2)}\times \overline{\cU}$ via each of expanding operators of the pair $(K_1^d,K_1'^d)$ in $\cH_{\underline{\alpha}}$.
In the first case, since the contracting operators of the pairs $(K_1^d,K_1'^d)$ and $(K_d,K_d')$ coincide, one can map $(w,u)$ to some $(w'_u,u)\in (W_d)_{(\varepsilon_2)}\times \overline{\cU}$ by operators of $(K_d,K'_d)$. By Remark \ref{remark: finitely many consequtive contraction}, 
  after a finite number of iteration of contracting operators we reach to some pair $(w^*_{u},u)\in (W_d)_{(\varepsilon_2)}\times \overline{\cU}$ such that we are able to map $w^*_{u}$ into $(W_d)_{(\varepsilon_3)}\supset (W_d)_{(\varepsilon_2)}$ by an expanding operator of the pair $(K_1^d,{K_1'}^d)$. Thus, it only remains to resolve the second case. 
  Since $\cU$ satisfies covering condition with respect to the set $\{M_1^{-1},\dots,M_{d^2}^{-1}\}$, there is $1\leq j\leq d^2$ such that $M_j^{-1}u=u'\in \cU$.
  Thus, the expanding operator from the set $\cH_{\underline{\alpha}}$ which is perturbed to have the $\SL$ part equal to $M_j^{-1}$ maps $(w,u)$ to $(w'',u')$ where $w''$ is $\varepsilon_3$-close to $w'$. Therefore, $(w'',u')\in (W_d)_{(\varepsilon_3)}\times \cU$.
  This implies that $\cW$ satisfies covering with respect to the finite family of renormalization operators with lengths of words less than $T$.
  The strong covering is concluded via Lemma \ref{lemma covering to strong}.
\end{proof} 

\begin{proof}[Proof of Theorem \ref{thm: Main: Real example}] According to Lemma \ref{lemma: immediately covering proof}, since $\cW$ satisfies strong  covering with respect to a finite subfamily of renormalization operators of the pair $(K_d,K'_d)$, there exists a relatively compact $\cW'\supset \cW$ that satisfies the strong covering condition with respect to the (finite) generating family $\cR_1^*$ of renormalization operators of the pair $(K_d,K_d')$. Then, since $\id\in \Aff(d,\R)$ lies in $\cW$, stable intersection of $K_d$ and $K_d'$ follows from Theorems \ref{thm: covering criterion in affine cantor sets with full shift} and \ref{thm: covering criterion for intersection}.

To estimate the dimension of $K_d$, we use the system of covering $\{G(\ua)\}_{\ua\in X_n}$, where $X_n\subset \Sigmafin$ denotes the set of all elements of $\Sigmafin$ with length $n$. In particular, for each $n\in \N$, $K_d\subset \bigcup_{\ua\in X_n} G(\ua)$. Note that for all $a\in\mathfrak{A}$, $G(a)$ is contained in the unit square $B$ and $G(\ua):= f_{\ua}(G(a_n))$ for each $\ua\in X_n$.
Recall that the construction of $K_1$ in \S \ref{sec: cantor sets in dimension one} depends on two parameters $\tau$ and $N$. 
Now, let $\tau'<\tau$ be positive. One can take $\delta>0$ in the proof of Theorem \ref{thm: there exists covering criterion in Rd} in the construction of $K_d$ such that for any affine generator $f_{(a,b)}$ of $K_d$ we have that $\|Df_{(a,b)}\|_{\it op}\leq (N+\tau')^{-1}$. Thus, $\diam(G(\ua))\leq (N+\tau')^{-n}$ for any $\ua\in X_n$  and $n\in \N$. 
Therefore,
$$\sum_{\ua\in X_n} \diam(G(\ua))^s \leq 6^{nd}\cdot (N+\tau')^{-sn}= (6^d/(N+\tau')^s)^n$$
tends to zero as $n \to \infty$ for any $s>{d\cdot \ln 6}/{\ln (N+\tau')}$.
This implies that 
$$\dimHD(K_d)\leq \dfrac{d\cdot \ln 6}{\ln (N+\tau')}.$$
In particular, given $\varepsilon>0$, if  $N$ is taken large enough, then $\dimHD(K_d)<\varepsilon$.
\end{proof}

\subsection{Affine bunched holomorphic Cantor sets: proof of Theorem \ref{thm: Main: Complex example}} 
To construct a pair of bunched holomorphic Cantor sets in $\C^d\cong \R^{2d}$ with an open relatively compact set satisfying the strong covering condition we follow the same argument of previous sections with a small change in the last steps.  
From Lemma \ref{lem: affine left and right action details} we have that 
$\Aff(d,\C)\cong \R^{2d} \rtimes \left(\R^*\times \T^1\times \SL^*(d,\C)\right)$ where $\T^1$ is the one-dimensional torus. In addition, we have $\R^{2d}\rtimes \R^*$ is isomorphic to a subgroup of $\Aff_{\id}(d,\C)$ under the injection map 
$(v,a)\hookrightarrow (v,a,0,\id)$. 
Therefore, renormalization operators of the pair $(K_1^{2d},K_1'^{2d})$ are in the subgroup 
$\R^{2d}\rtimes \R^*$ of $\Aff(d,\C)$. We will have the following analog of the Lemma \ref{lem: Wd satisfies covering} with the same proof.

\begin{lemma}\label{lem: Wd satisfies covering complex}  $W_{2d}$ has expanding $3^{2d}$-cover with respect to renormalization operators of the pair $(K_1^{2d},K_1'^{2d})$. Moreover, 
there exists $\varepsilon_1.\varepsilon_2>0$ such that if 
$B_{\varepsilon_1}(\id)$ and $I_{\varepsilon_1}$ are the $\varepsilon_1$-neighborhoods of $\id$ and $0$ in $\SL(d,\bC)$ and  $\T^1$ respectively, then for any $(w,\theta,u)\in  B_{\varepsilon_2}(W_{2d}) \times \overline{I_{\varepsilon_1}}\times \overline{B_{\varepsilon_1}(\id)}$ either there exists a contracting operator or at least $3^{2d}$ expanding operators that maps $(w,\theta,u)$ to 
$(w',\theta,u)\in (W_{2d})_{(\varepsilon_2)}\times \overline{I_{\varepsilon_1}} \times \overline{B_{\varepsilon_1}(\id)}$.
\end{lemma}

We have the following theorem, analogous to Theorem \ref{thm: there exists covering criterion in Rd}.

\begin{theorem}\label{thm: there exists covering criterion in Cd}
There exists a pair $(K^{\rm Hol}_d,K'^{\rm Hol}_d)$ of affine holomorphic Cantor sets with full shift symbolic type in $\C^d$ arbitrarily close to the pair $(K_1^{2d},K_1'^{2d})$ and  $I\times \cU\subset \T^1\times  \SL(d,\C)$ which is an open relatively compact set such that 
 $\cW: = W_{2d} \times I\times \cU$ satisfies strong covering condition with respect to a finite family of renormalization operators of the pair $(K^{\rm Hol}_d,K'^{\rm Hol}_d)$.
\end{theorem}
\begin{proof}
    The argument is similar to the proof of Theorem \ref{thm: there exists covering criterion in Rd}. The only difference is that we use the complex version of Lemma \ref{lem: generating a convering in SL with d^2 maps} to obtain the set of matrices $\cM:=\{M_1,\dots,M_{2d^2-1}\}\subset \SL(d,\C)$ together with $\cU\subset \SL(d,\C)$ such that
    $\cU$ satisfies covering condition with respect to $\cM$ (note that the dimension of the Lie group $\SL(d,\C)$ as a subgroup of $\SL(2,\R)$ is $2d^2-2$). Then since
    $3^{2d}>4d^2-2$ we may duplicate  $\cM$ into the set of matrices 
$$\tilde{\cM}:=\{e^{i\delta}M_1,e^{-i\delta}M_1,\dots, e^{i\delta}M_{2d^2-1},e^{-i\delta}M_{2d^2-1}\}\subset \T^1\times \SL(d,\C)$$
to obtain the covering condition on the set $I\times \cU$ with respect to $\tilde{\cM}$ where $I: = I_{\varepsilon}$ is the $\varepsilon$-neighborhood of $0$ in $\T^1$. 
\end{proof}
\begin{proof}[Proof of Theorem \ref{thm: Main: Complex example}]
   The proof is identical to the proof of Theorem \ref{thm: Main: Real example} except that here we use Theorem \ref{thm: there exists covering criterion in Cd} and Lemma \ref{lem: Wd satisfies covering complex} in place of Theorem \ref{thm: there exists covering criterion in Rd}
 and Lemma \ref{lem: Wd satisfies covering}, respectively. The dimension estimate remains unchanged. 
\end{proof}

\begin{remark} Similar arguments can be applied to other Lie groups. In the case of the group $O(d)$, this method yields pairs of stably intersecting Cantor sets within the space of conformal Cantor sets.
\end{remark} 

\appendix

\section{A convergence lemma}\label{sec: appendix convergence lemma}
In this appendix we prove a known convergence lemma
which allowed us to study infinitesimal geometry of bunched Cantor sets.

\begin{lemma}
\label{lem: appendix- convergence lemma}
Let $\alpha > 0$, $V, U\subset \R^d$ be bounded open sets that $\overline{V}\subset U$ and $\{g_m\}_{m\in \N}$ be a sequence 
in $\diffloc{1+\alpha}(\R^d)$ such that $\operatorname{Dom}(g_m)\cap \operatorname{Im}(g_m)\supset U$. If $\sum_{m=1}^{\infty} d_{\cC^{1+\alpha}}(g_m, \id) < \infty$, then the sequence
$\{g_m\circ g_{m-1} \circ \cdots \circ g_1\}_{m\in \N}$ converges in $\diffloc{1+\alpha}(\R^d)$ on $\overline{V}$.
Moreover, the limit varies $\cC^1$-continuously with respect to the sequence $\{g_m\}_{m\in \N}$ in the following sense. For every $\varepsilon>0$ there exists $\delta>0$ such that if $\{\tilde{g}_m\}_{m\in \N}$ be a sequence in $\diffloc{1+\alpha}(\R^d)$ with 
$\sum_{m=1}^{\infty} d_{\cC^{1+\alpha}}(g_m, \tilde g_m) < \delta$ then the resulting limits 
$$G:= \lim_{m\rightarrow \infty} g_m\circ g_{m-1} \circ \cdots \circ g_1,\quad\tilde G:= \lim_{m\rightarrow \infty} \tilde g_m\circ \tilde g_{m-1} \circ \cdots \circ \tilde g_1$$
restricted to $\overline{V}$ are $\varepsilon$-close in $\diffloc{1}(\R^d)$.
\end{lemma}
To prove this lemma we shall use the following H\"older estimates from \cite[Theorem A.8]{Hormander}. For $U,V$ as in Lemma \ref{lem: appendix- convergence lemma}, $\beta\geq 1 $, $\varphi:U\to \R^d$ be a $\cC^\beta$ map and $g,h \in \diffloc{\beta}(\R^d)$ whose domain and image both contain $U$, there is a constant $\lambda_\beta$ such that the following holds on $\overline{V}$
    \begin{align}\label{Hormander}
        \|\varphi\circ g\|_{\cC^\beta} &\leq \lambda_\beta\left(\|\varphi\|_{\cC^1}\|g\|_{\cC^\beta}+\|\varphi\|_{\cC^\beta}\|g\|_{\cC^1}^\beta+\|\varphi\|_{\cC^0}\right)\notag\\
        &\leq \lambda_{\beta}\cdot \|\varphi\|_{\cC^\beta}\left(\|g\|_{\cC^\beta}+ \|g\|_{\cC^1}^\beta+ 1\right).
    \end{align}
    Moreover, if $\beta>1$ on $\overline{V}$ we have
    \begin{align}\label{eq: appendix holder estimate g-h}
        \|\varphi\circ g-\varphi\circ h\|_{\cC^1}&\leq 2\|\varphi\|_{\cC^1}\cdot \|g-h\|_{\cC^1}+\|h\|_{\cC^1}\cdot \|g-h\|_{\cC^0}^{\beta-1}\cdot \|\varphi\|_{\cC^{\beta}}\notag \\
        &\leq \|\varphi\|_{\cC^{\beta}}\cdot \left(2\|g-h\|_{\cC^1}+\|h\|_{\cC^1}\cdot \|g-h\|_{\cC^0}^{\beta-1}\right)\notag \\
        & \leq \|\varphi\|_{\cC^{\beta}}\cdot \left(2+\|h\|_{\cC^1}\right)\cdot \left(\|g-h\|_{\cC^1}+\|g-h\|_{\cC^1}^{\beta-1}\right).
    \end{align}
    
\begin{proof}[Proof of Lemma \ref{lem: appendix- convergence lemma}] Denote 
$$G_m:=g_m\circ g_{m-1} \circ \cdots \circ g_1,\quad
    \phi_m := g_m - \id,\quad\Phi_m := G_m -\id.$$
    Note that $\operatorname{Dom}(\phi_m)=\operatorname{Dom}(g_m)\supset U$ and 
    $\operatorname{Dom}(\Phi_m)=\operatorname{Dom}(G_m)\supset U$. We have
    $$\id+\Phi_m  = G_m = g_m \circ G_{m-1} =  \id+ \Phi_{m-1}+\phi_m \circ (\id+\Phi_{m-1}),$$
    so
    \begin{equation}\label{eq: Appendix-psi_m}
    \Phi_m - \Phi_{m-1} = \phi_m \circ (\id+\Phi_{m-1}).
    \end{equation}
    In addition, the assumption $\sum_{m=1}^{\infty} d_{\cC^{1+\alpha}}(g_m, \id) < \infty$ implies that 
    \begin{equation}\label{eq: Appendic-phi_m norm}
        \sum_{m=1}^{\infty} \|\phi_m\|_{\cC^{1+\alpha}}<\infty.
    \end{equation}
    Since 
    $\|\phi_m \circ (\id+\Phi_{m-1})\|_{\cC^0}\leq \|\phi_m\|_{\cC^0}$, by (\ref{eq: Appendix-psi_m}) and (\ref{eq: Appendic-phi_m norm}) we get that $\{\Phi_n\}_{n\geq 1}$ is a Cauchy sequence in $\cC^0$ topology which yields the $\cC^0$ convergence.

    Let $a_n := \|\phi_n\|_{\cC^{1+\alpha}}$, $b_n := \|\Phi_n\|_{\cC^1}$ and $c_n := \|\Phi_n\|_{\cC^{1+\alpha}}$. Then (\ref{eq: Appendix-psi_m}) and (\ref{Hormander}) for $\beta=1$ implies that
    \begin{align}\label{Hormander2}
    \|\Phi_m - \Phi_{m-1}\|_{\cC^1} &=\|\phi_m \circ  (\id+\Phi_{m-1})\|_{\cC^1}\notag\\ 
    &\leq \lambda_1\cdot \|\phi_m\|_{\cC^1} \cdot\left( 2\|\id+\Phi_{m-1}\|_{\cC^1} + 1\right)\notag\\
    & \leq 2\lambda_1\cdot a_m \cdot \left(b_{m-1}+ \ell_0\right),
    \end{align}
    where $\ell_0 :=2+\diam(U)$. Then, since $b_m - b_{m-1}\leq \|\Phi_m - \Phi_{m-1}\|_{\cC^1}$ we have
        $$b_m +\ell_0\leq (b_{m-1}+ \ell_0)(1+2\lambda_1\cdot a_m).$$
    It follows that for $\ell_1 := (b_1+\ell_0)/(1+2\lambda_1\cdot a_1)$ and $m\in \N$,
    \begin{align*}
    b_m +\ell_0 &\leq  \ell_1\cdot \prod_{j=1}^{m} (1+2\lambda_1\cdot a_j).
    \end{align*}
    Since $\sum_{j=1}^{\infty} a_j<\infty$, the product $\prod_{j=1}^{\infty}(1+2\lambda_1\cdot a_j)$ converges and the sequence $\{b_n\}_{n\in \N}$ is bounded by some $b>0$. 
    In particular, $\|\id+ \Phi_{m-1}\|_{\cC^1}$ is bounded by $b+\ell_0$.
    To prove $\cC^{1+\alpha}$ convergence we proceed with similar calculations. We have
    \begin{align}\label{Hormander3}
    \|\Phi_m - \Phi_{m-1}\|_{\cC^{1+\alpha}} &=\|\phi_m \circ  (\id+\Phi_{m-1})\|_{\cC^{1+\alpha}}\notag\\ 
    &\leq \lambda_{1+\alpha}\cdot \|\phi_m\|_{\cC^{1+\alpha}}\cdot \left(\|\id+\Phi_{m-1}\|_{\cC^{1+\alpha}}+ \|\id+\Phi_{m-1}\|_{\cC^{1}}^{1+\alpha} +1\right)\notag\\
    &\leq 2\lambda_{1+\alpha}\cdot \|\phi_m\|_{\cC^{1+\alpha}}\cdot \left(\|\id+\Phi_{m-1}\|_{\cC^{1+\alpha}}+(b+\ell_0)^{1+\alpha} + 1\right)\notag\\
    & \leq 2\lambda_{1+\alpha}\cdot a_m\cdot  \left(c_{m-1}+ \ell_2\right),
    \end{align}
    where $\ell_2:=(b+\ell_0)^{1+\alpha} + \ell_0+1$. Therefore, similar to the argument above there exists $\ell_3:= (c_1+\ell_2)/(1+2\lambda_{1+\alpha}\cdot a_1)$ such that for all $m\in \N$,
    \begin{align*}
    c_m +\ell_2 &\leq  \ell_3\cdot \prod_{j=1}^{m} (1+2\lambda_{1+\alpha}\cdot a_j)
    \end{align*}
     which implies that the sequence $\{c_m\}_{m\in \N}$ is bounded by some $c>0$. So by (\ref{Hormander3})
     \begin{equation}\label{eq: appendix-3}
         \|\Phi_m - \Phi_{m-1}\|_{\cC^{1+\alpha}} \leq 2\lambda_{1+\alpha} (c+\ell_2) \cdot \|\phi_m\|_{\cC^{1+\alpha}}.
     \end{equation}
     Hence, $\{\Phi_m\}_{m\in \N}$ is Cauchy by \eqref{eq: Appendic-phi_m norm} and we obtain the $\cC^{1+\alpha}$ convergence. 
     The sequence $DG_m$ converges to a non-singular matrix $DG$ in the $\cC^1$ norm which implies that $G$ has a $\cC^1$ inverse equal to $\lim_{m\rightarrow \infty}G_m^{-1}$. Then using \cite[Theorem A.9]{Hormander} with similar calculations as above it can be proved that
    $|DG_m^{-1}-DG^{-1}|_{\alpha}$ tends to zero as $m$ goes to infinity 
    which implies that $\lim_{m\rightarrow \infty}G_m = G$
    is in $\diffloc{1+\alpha}(\R^d)$. This completes the proof of first part of the lemma.
    
    To prove the continuity of limit with respect to the sequence $\{g_m\}_{m\in \N}$, let 
    $$\tilde G_m: =\tilde g_m\circ \tilde g_{m-1} \circ \cdots \circ \tilde g_1,\quad \tilde \Phi_m := \tilde G_m - \id, \quad \tilde \phi_m: = \tilde g_m -\id$$
    with $\sum_{m\in \N} d_{\cC^{1+\alpha}}(g_m, \tilde g_m) < \delta$ for some $\delta>0$ to be determined later.
    By (\ref{Hormander})
     \begin{align}\label{Hormander4}
     \|\Phi_m-\tilde \Phi_m\|_{\cC^{1}}   
     = &~\|(\Phi_{m-1}-\tilde \Phi_{m-1}) + (\phi_m-\tilde \phi_m)\circ (\id+\tilde \Phi_{m-1})\notag\\
     & +\phi_m\circ (\id+\Phi_{m-1})-\phi_m \circ (\id+\tilde \Phi_{m-1})\|_{\cC^{1}}\notag\\
    \leq &~ \|\Phi_{m-1}-\tilde \Phi_{m-1}\|_{\cC^{1}}\notag  \\
    &+\lambda_{1}\cdot \|\phi_m-\tilde \phi_m\|_{\cC^{1}}\cdot  
    \left( 2\|\id+\Phi_{m-1}\|_{\cC^{1}} + 1\right)\notag \\
    &+\|\phi_m\circ (\id+\Phi_{m-1})-\phi_m \circ (\id+\tilde \Phi_{m-1})\|_{\cC^{1}}\notag\\
    \leq &~ \|\Phi_{m-1}-\tilde \Phi_{m-1}\|_{\cC^{1}}+ \lambda_{1}\cdot \ell_4\cdot \|\phi_m-\tilde \phi_m\|_{\cC^{1}}\notag\\
    & +\|\phi_m\circ (\id+\Phi_{m-1})-\phi_m \circ (\id+\tilde \Phi_{m-1})\|_{\cC^{1}},
         \end{align}
         where $\ell_4:=1+2(\ell_0+b)$. Using the estimate \eqref{eq: appendix holder estimate g-h} one has
         \begin{align} 
         \|\phi_m\circ &(\id +\Phi_{m-1})-\phi_m \circ (\id+\tilde \Phi_{m-1})\|_{\cC^{1}}\notag \\ 
         &\leq\|\phi_m\|_{\cC^{1+\alpha}} \cdot (2+\|\Phi_{m-1}\|_{\cC^1})\cdot \left(\|\Phi_{m-1}-\tilde \Phi_{m-1}\|_{\cC^1}+\|\Phi_{m-1}-\tilde \Phi_{m-1}\|_{\cC^1}^{\alpha}\right). \label{eq: last term in Hormander4}
         \end{align}
         Let $s_m := \|\Phi_m-\tilde \Phi_m\|_{\cC^{1+\alpha}}$ and $t_m := \|\phi_m-\tilde \phi_m\|_{\cC^{1+\alpha}}$.
         By \eqref{Hormander4} and \eqref{eq: last term in Hormander4} we have
         \begin{align}\label{Hormander5}
             s_m &\leq s_{m-1} + \lambda_{1} \ell_4\cdot t_m +(2+b)\cdot a_m \cdot (s_{m-1}+s_{m-1}^{\alpha})\notag\\
             &\leq s_{m-1}(1+\ell_4\cdot a_m)+ s_{m-1}^{\alpha}\cdot \ell_4\cdot a_m + \lambda_{1} \ell_4\cdot t_m,
         \end{align}  
    Since $\sum_{m=1}^{\infty} a_m<\infty$, there exists $m_0\in \N$ such that $\sum_{j\geq m_0} a_j<(8\ell_4)^{-1}$. Thus, by $\sum_{m\in \N} t_m<\delta$, if $\delta<(4\lambda_1\ell_4)^{-1}$ we have
    $\sum_{j\geq m_0} (2\ell_4 a_j+\lambda_1\ell_4 t_j)<1/2$. On the other hand, by continuity of finite composition $g_{m_0}\circ\cdots\circ g_1$ with respect to $g_1,\dots,g_{m_0}$, there is a small $\delta_0>0$ such that for $\delta<\delta_0$ we have
    $s_1,s_2,\dots,s_{m_0}<1/2$. Arguing by induction we show that for $m\geq m_0$ one has 
    $$s_m\leq 1- \sum_{j\geq m+1} (2\ell_4a_j+\lambda_1\ell_4t_j).$$ 
    Since $s_{m-1}<1$, we have $s_{m-1}^{\alpha}<1$. So by \eqref{Hormander5} we can write
    \begin{align*}
    s_m&\leq s_{m-1}+(2\ell_4a_m+\lambda_1\ell_4t_m)\notag \\
    &\leq 1-\sum_{j\geq m} (2\ell_4a_j+\lambda_1\ell_4t_j) + (2\ell_4a_m+\lambda_1\ell_4t_m)\notag \\
    & = 1-\sum_{j\geq m+1} (2\ell_4a_j+\lambda_1\ell_4t_j),
    \end{align*}
    and the induction step follows. Thus, $s_m<1$ for all $m\in \N$. So by \eqref{Hormander5} we have 
    \begin{align}\label{eq: s_m linear recursive estimation}
        s_m&\leq s_{m-1}(1+\ell_4\cdot a_m)+ s_{m-1}^{\alpha}\cdot \ell_4\cdot a_m + \lambda_{1} \ell_4\cdot t_m \notag \\
        & \leq s_{m-1} +\ell_5\cdot u_m,
    \end{align}
    where $\ell_5 := \ell_4(1+\lambda_1)$ and $u_m: = 2a_m+t_m$. Thus, for all $m,N\in \N$ we have
    \begin{align}\label{eq: s_m final estimation}
        s_{m+N}&\leq s_m + \ell_5\sum_{j=1}^{N} u_{m+j}\notag \\
        &\leq s_m + \ell_5\sum_{j> m}u_j.
    \end{align}
    So, given $\varepsilon>0$ there is a large enough $m_1$ such that $\sum_{j\geq m_1}a_j <\varepsilon(8\ell_5)^{-1}$, which implies that for $\delta<\varepsilon(4\ell_5)^{-1}$ and all $N\in \N$,
    $s_{m_1+N}<s_{m_1}+ \varepsilon/2$. By continuity of the finite composition $g_{m_1}\circ\cdots\circ g_1$ with respect to $g_1,\dots,g_{m_1}$, there is a small $\delta_1>0$ such that for $\delta<\delta_1$ we have $s_{m_1}< \varepsilon/2$.
    Therefore, for all $N\in \N$ when $\delta:= \min\{\delta_0,\delta_1,\varepsilon(4\ell_5)^{-1},(4\lambda_1\ell_4)^{-1}\}$ we have $s_{m_1+N}<\varepsilon$. Note that $m_1,m_0$ are depended only on the sequence $\{g_m\}_{m\in \N}$. This yields the second part of the lemma.  
\end{proof}

\begin{remark} Continuity of the limit $G= \lim_{n\to \infty} G_m$ also holds in $\cC^{1+\alpha'}$ topology for any positive $\alpha'<\alpha$ rather than $\cC^1$ topology as stated in Lemma \ref{lem: appendix- convergence lemma}.
    More generally, same calculations can be done to prove $\cC^{\ell}$-convergence of maps $G_m$ for any $\ell\geq 1$ provided that
    $\{g_m\}_{m\in \N}$ are maps in $\diffloc{\ell}(\R^d)$.
    Moreover, when $\ell>1$ the limit is continuous in 
    $\cC^{\ell'}$ topology for any $1\leq \ell'<\ell$.
\end{remark}
\begin{remark}\label{remark: appendix speed of convergence}
    If we have $d_{\cC^{1+\alpha}}(g_m, \id)\leq C\cdot \gamma ^m$ for some positive constants $C$ and $\gamma<1$ then the speed of convergence $G_m \rightarrow H$ is about $\gamma^{m}$. In fact there is a constant $C_1$ such that $d_{\cC^{1+\alpha}}(G_m, H) \leq C_1\cdot \gamma^{m}$.
\end{remark}

From the triangle inequality 
  $$
  \big|d_{\cC^{1+\alpha}}(g_m, \tilde g_m)-d_{\cC^{1+\alpha}}(\tilde g_m, \id)\big| < d_{\cC^{1+\alpha}}(g_m, \id),
    $$
  we conclude that whenever $\sum_m d_{\cC^{1+\alpha}}(g_m, \id) < \infty$, then
  $\sum_m d_{\cC^{1+\alpha}}(\tilde g_m, \id) < \infty$ is equivalent to $\sum_m d_{\cC^{1+\alpha}}(\tilde g_m, g_m) < \infty$. Therefore, according to the proof of continuity of the limit $H = \lim_{m\to \infty} G_m$ with respect to the sequence $\{g_m\}_{m\in \N}$ in Lemma \ref{lem: appendix- convergence lemma}, specifically the relation \eqref{eq: s_m final estimation}, we have the following continuity lemma which is a variant of the continuity part of Lemma \ref{lem: appendix- convergence lemma}.
  
  \begin{lemma}\label{lem: appendix continuity second version} 
  Let $C,\epsilon>0$ and $0<\gamma <1$ and $\{g_m\}_{m\in \N}$ be a sequence of maps satisfying conditions in Lemma \ref{lem: appendix- convergence lemma} such that $d_{\cC^{1+\alpha}}(g_m, \id)<C\cdot\gamma^m$ for all $m\in \N$. Then there exists $\delta>0$ and  an integer $m_0>0$ such that for any sequence of maps $\{\tilde g_m \}_{m=1}^{\infty}$ such that $d_{\cC^{1+\alpha}}(\tilde g_m, \id)<C\cdot\gamma^m$ for all $m\in \N$, $d_{\cC^{1}}(H, \tilde H)< \epsilon$ if $d_{\cC^{1}}(g_k, \tilde g_k)< \delta$ for $k=1,2,\dots,m_0$.
  \end{lemma}

\section{Affine estimation}

Here, we prove the following general linear estimation.

\begin{lemma}\label{lemma estimation of chosing linear instead of non linear part}
Let $U\subset \R^d$ an open set, $X\subset U$ a set such that it's convex hull is also contained inside $U$, $p\in X$ a given point, $\phi:U\rightarrow \R^d$ be a $\cC^{1+\alpha}$ map with $\cC^{1+\alpha}$ norm $C_{\phi, X}$ on the domain $X$ and $B\in \Aff(d,\R)$ be an affine contracting map such that there are constants $C_{B,X}$ and $C_{\phi}$ where 
$$C_{B,X}^{-1}\leq \dfrac{\|DB\|_{\it op}}{\operatorname{diam}(X)}\leq C_{B,X},\;\; C_{\phi} := \sup_{x\in X} \|(D\phi_x)^{-1}\|_{\it op}.$$
Denote $A_{\phi,p}$ as the affine estimation of $\phi$ which is an affine map with derivative equal to $D\phi_p$ and maps $p$ to $\phi(p)$. Then the following holds on the domain $B^{-1}(X)$
$$\|B^{-1}\circ A_{\phi,p}^{-1}\circ \phi \circ B - \id\|_{\cC^{1+\alpha}}<C'\cdot \|DB^{-1}\|_{\textit{op}}\cdot \|DB\|_{\textit{op}}\cdot \diam(X)^{\alpha},$$
where 
$C'= C_{\phi,X}\cdot (1+2C_{B,X})\cdot C_\phi$.
\end{lemma}

\begin{proof} Denote $F := B^{-1}\circ A_{\phi,p}^{-1}\circ \phi \circ B - \id$.
    Given an affine map $H\in \Aff(d,\R)$ and $\psi$ a map on $\R^d$, then one has
    \begin{equation}\label{eq: affine lemma computation}
        H\circ \psi-\id = DH \circ (\psi - H^{-1}).
    \end{equation}
    Taking $H=B^{-1}\circ A_{\phi,p}^{-1}$ and $\psi =\phi\circ B $  in \eqref{eq: affine lemma computation} implies that
    $$F = 
    DB^{-1}\circ DA_{\phi,p}^{-1} \circ (\phi - A_{\phi,p})\circ B.
    $$
    On the other hand, using the H\"older regularity of $\phi$ on domain $X$ we have 
    \begin{equation}\label{eq: holder regularity of phi}
        \|\phi-A_{\phi,p}\|_{\cC^0}\leq C_{\phi,X} \cdot \operatorname{diam}(X)^{1+\alpha},\;\; \|D\phi_x-D\phi_p\|_{\it op}\leq C_{\phi,X} \cdot \operatorname{diam}(X)^{\alpha}.
    \end{equation}
    Then using above relations we have the following estimates on $\cC^{1+\alpha}$ norm of $F$.
    \begin{align*}
        \|F\|_{\cC^0}&\leq  \|DB^{-1}\|_{\it op}\cdot 
        \|(D\phi_p)^{-1}\|_{\it op}\cdot \|\phi-A_{\phi,p}\|_{\cC^0}\\
        & \leq C_{\phi,X}\cdot C_{B,X}\cdot C_{\phi}\cdot  \|DB^{-1}\|_{\it op}\cdot \|DB\|_{\it op} \cdot \operatorname{diam}(X)^{\alpha},\\
        \|DF\|_{\cC^0}&\leq \|DB^{-1}\|_{\it op}\cdot \|(D\phi_p)^{-1}\|_{\it op}\cdot \sup_{x\in X}\|D\phi_x-D\phi_p\|_{\it op} \cdot \|DB\|_{\it op}\\
        & \leq C_{\phi,X}\cdot C_{\phi}\cdot  \|DB^{-1}\|_{\it op}\cdot \|DB\|_{\it op} \cdot \operatorname{diam}(X)^{\alpha},\\
        |DF|_{\alpha} &
        =\sup _{x,y \in B^{-1}(X)} \dfrac{\|DB^{-1}\circ (D\phi_p)^{-1}\circ (D\phi_{B(x)}-D\phi_{B(y)})\circ DB\|_{\it op}}{|x-y|^{\alpha}}\\
        &\leq C_{\phi,X}\cdot C_{\phi}\cdot C_{B,X}\cdot 
        \|DB^{-1}\|_{\it op}\cdot \|DB\|_{\it op} \cdot \operatorname{diam}(X)^{\alpha}.
    \end{align*}
    Thus, since $\|F\|_{\cC^{1+\alpha}}\leq \|F\|_{\cC^0}+\|DF\|_{\cC^0} + |DF|_{\alpha}$ we conclude the proof.
\end{proof}

\section{Operations on the space of affine maps}\label{appendix affine actions}

 Here, we study the structure of left and right action of $\Aff(d,\bF)$ on itself. Denote  
$$\SL^{\pm}(d, \bR):= \{A \in \SL(d, \bR) ~|~ \det(A) = \pm 1
\}.$$
Let $\SL^\star(d,\bF)$ be a notation for $\SL(d,\bF)$ in the cases that $(d,\bF)=(2k+1,\bR)$ or 
$\bF = \bC$ and for $\SL^\pm(d,\bR)$ in the case of $(d,\bF)=(2k,\bR)$. Then we have the correspondence $\Aff(d,\bR) \cong \R^d \times \R^*\times \SL^\star(d,\R)$ via the homeomorphism
$$(v,\lambda ,A)\mapsto [x\mapsto \lambda \cdot Ax+v]$$
with the inverse map 
$$[x\mapsto Ax+v] \mapsto \left(v,s_A,s_A^{-1} A\right)$$
where $s_A : = \frac{\operatorname{det}(A)}{|\operatorname{det}(A)|}\sqrt[d]{|\operatorname{det}(A)|}$. Similarly, 
$\Aff(d,\bC) \cong \bC^d\times  \bC^*\times \SL^\star(d,\bC)$
via the homeomorphism 
$$(v,\lambda e^{i\theta},A)\mapsto [z\mapsto \lambda e^{i\frac{\theta}{d}}\cdot Az+v] $$
where $\lambda>0$ is a real number, with the inverse map
$$
[z\mapsto Az+v] \mapsto \left(v,s_A,s_A^{-1} A\right),$$
where 
$s_A :=e^{i\operatorname{arg}(\operatorname{det}(A))}  \sqrt[d]{|\operatorname{det}(A)|}$.
Let $[x\mapsto Px+w]$ be an element of $\Aff(d,\bF)$. Then for any affine map $[x\mapsto Ax+v] \in \Aff(d,\bF)$ we have 
$$[x\mapsto Px+w] \circ [x\mapsto Ax+v] = [x\mapsto PAx+ (Pv+w)],$$
$$[x\mapsto Ax+v] \circ [x\mapsto Px+w] = [x\mapsto APx +(v+Aw)].$$
So, we have the following interpretation of the action of the group $\Aff(d,\bF)$ on itself.

\begin{lemma}\label{lem: affine left and right action details}
    The above correspondence is a group homomorphism with the following group operation on $  \bF^d \rtimes (\bF^* \times \SL^\star(d,\bF))$, 
    $$(v,a,A)*(w,b,B) = (a\cdot Aw+v,ab,AB ).$$
    This implies that $\bF^d\rtimes \bF^*$ is a subgroup of $\Aff(d,\bF)$ under the injection map $(v,a)\hookrightarrow (v,a,\id)$ with group operation
    \begin{equation}\label{eq: Wd is a subgroup of Aff group operation}
        (v,a)*(w,b) = (aw+v,ab).
    \end{equation}  
\end{lemma}

 We denote the subgroup $\bF^d \rtimes \bF^*$ in above definition as $\Aff_{\id}(d,\bF)$, i.e. subgroup of affine maps that are composition of a homothety and translation.


\providecommand{\bysame}{\leavevmode\hbox to3em{\hrulefill}\thinspace}

\end{document}